\documentclass{scrartcl}
\pdfoutput=1

\usepackage{amsmath,
amssymb,
amsfonts,
amsthm,
amscd,
amsbsy,
mathtools,
thmtools,
listings,
rotating,
stmaryrd,
enumitem,
shuffle,
mleftright,
float,
microtype,
mlmodern
}

\usepackage[dvipsnames]{xcolor}
\newcommand\diff{\mathop{}\!\mathrm{d}}

\usepackage{ifdraft}

\usepackage{ytableau}
\definecolor{todocolor}{HTML}{fa73ff}
\usepackage[linecolor=white,backgroundcolor=white,bordercolor=white,textsize=tiny,textcolor=todocolor]{todonotes}

\usepackage[normalem]{ulem}

\usepackage{imakeidx}
\makeindex[intoc,title=Symbol index]

\usepackage{tikz,tkz-tab, tikz-cd}
\usepackage{caption}

\usepackage{caption}
\usepackage{subcaption}
\usepackage{bookmark}
\usepackage[nameinlink]{cleveref}
\usepackage{mhequ}
\usepackage{longtable}
\usepackage{booktabs}
\usepackage{authblk}
\usepackage[ocgcolorlinks]{ocgx2}
\hypersetup{linkcolor=BrickRed, citecolor=RoyalBlue}

\newcommand{\frontstick}{\,\raisebox{-1pt}{\begin{tikzpicture}
\draw [line width=1pt,] (0,0)--(0,0.25);
\end{tikzpicture}}\kern+2pt}
\newcommand\freecmseriesh{\overline{\mathfrak{g}_{-1}}}
\newcommand\freecmseriesg{\overline{\mathfrak{g}_{0}}}

\newcommand\alphafree{{\alpha^X}}
\newcommand\betafree{{\beta^X}}

\newcommand\R{\mathbb{R}}
\newcommand\N{\mathbb{N}}

\newcommand\gray[1]{{\color{gray}#1}}

\newcommand\unit{\mathsf{e}} 
\usepackage{graphicx}
\newcommand{\mirrorLsh}{\scalebox{-1}[1]{$\Lsh$}}

\newcommand\rightuparrow{\rotatebox[origin=c]{270}{$\Lsh$}}
\newcommand\uprightarrow{\rotatebox[origin=c]{0}{$\mirrorLsh$}}

\newcommand\GL{\operatorname{GL}}

\newcommand\Aut{\operatorname{Aut}}
\newcommand\Der{\operatorname{Der}}
\newcommand\f{\mathfrak{f}}
\newcommand\g{\mathfrak{g}}
\newcommand\h{\mathfrak{h}}

\newcommand\FL{\mathsf{L}}

\newcommand\im{\operatorname{im}}
\newcommand\DEF[1]{\textbf{\textup{#1}}}
\newcommand\eps\epsilon
\newcommand\ZZ{\mathbb{Z}}
\newcommand\BCH{\mathsf{BCH}}
\newcommand\Z[1]{{\color{cyan}\mathtt{Z}}_{#1}}

\newcommand\A{\alpha}

\newcommand\IIS{S} 
\newcommand\Area{\mathsf{Area}}

\newcommand\LL{\mathcal{L}}
\newcommand\Ad{\operatorname{Ad}}
\newcommand\ad{{\mathsf{ad}}}

\renewcommand\mag{\omega} 

\newcommand\coker{\operatorname{coker}}

\newcommand\feedbackGG{\mathfrak T}        
\newcommand\feedbackAA{\mathfrak t}        
\newcommand\actionGG{m}             
\ifdraft
{
\newcommand\actionGA{{\dot m}}      
\newcommand\actionAA{{\ddot m}}     
\newcommand\ddotn{{\ddot n}}
}{
\newcommand\actionGA{%
	\overset{\text{\tikz[baseline]{\draw[fill=white] (0,0) circle (1pt);\fill (0.12,0) circle (1pt);}}}{m}%
}
\newcommand\actionAA{%
	\overset{\text{\tikz[baseline]{\fill (0,0) circle (1pt);\fill (0.12,0) circle (1pt);}}}{m}%
}
\newcommand\ddotn{%
	\overset{\text{\tikz[baseline]{\fill (0,0) circle (1pt);\fill (0.12,0) circle (1pt);}}}{n}%
}
}

\newcommand\restr[2]{\ensuremath{\left.#1\right|_{#2}}}

\newcommand\sab{{\mathsf{sab}}}
\newcommand\ab{{\mathsf{ab}}}
\newcommand\freecm[1]{\mathfrak g_{#1}}
\newcommand\feedbackFreeCmGG{\mathfrak T}       
\newcommand\feedbackFreeCmAA{\tau}
\ifdraft
{
 \newcommand\actionFreeCmGA{\dot \triangleright}
 \newcommand\actionFreeCmAA{\ddot \triangleright}
}
{
	\newcommand\actionFreeCmGA{
		\overset{\text{\tikz[baseline]{\draw[fill=white] (0,0) circle (1pt);\fill (0.12,0) circle (1pt);}}}{\triangleright}%
	}
	\newcommand\actionFreeCmAA{
		\overset{\text{\tikz[baseline]{\fill (0,0) circle (1pt);\fill (0.12,0) circle (1pt);}}}{\triangleright}%
	}
}

\DeclarePairedDelimiter{\dsbr}{[\mkern-3mu[}{]\mkern-3mu]}

\newcommand\BB{\mathcal{B}}
\newcommand\CC{\mathcal{C}}
\newcommand\DD{\mathcal{D}}

\newcommand\CO{\mathcal{O}}
\newcommand\YY{\mathcal{Y}}
\DeclareMathOperator{\diam}{diam}

\newcommand\up{\uparrow}
\newcommand\dn{\downarrow}
\newcommand\rt{\rightarrow}
\newcommand\lt{\leftarrow}
\newcommand\PD{\mathcal{P}} 
\newcommand{\SD}{\mathcal{R}}
\newcommand{\RS}{\mathcal{X}}

\newcommand{\Hol}[1]{#1\textnormal{-H{\"o}l}}

\newcommand\id{\operatorname{id}}
\newcommand\SQUARE[2]{\mathop{{}_{\scalebox{0.5}{$#1$}}\square^{\scalebox{0.5}{$#2$}}}}
\newcommand\HOOK[2]{\mathop{{}_{\scalebox{0.5}{$#1$}}\!\Rsh^{\scalebox{0.5}{$#2$}}}}

\newcommand\DN[2]{\mathop{{}^{\scalebox{0.5}{$#1$}}_{\scalebox{0.5}{$#2$}}\!\!\dn}}
\newcommand\UP[2]{\mathop{{}^{\scalebox{0.5}{$#2$}}_{\scalebox{0.5}{$#1$}}\!\!\up}}
\newcommand\RT[2]{\mathop{{}_{\scalebox{0.5}{$#1$}}\!\!\rt_{\scalebox{0.5}{$#2$}}}}
\newcommand\LT[2]{\mathop{{}_{\scalebox{0.5}{$#2$}}\!\!\lt_{\scalebox{0.5}{$#1$}}}}

\newcommand\ENW{\raisebox{-0.5ex}{\includegraphics[width=.8em]{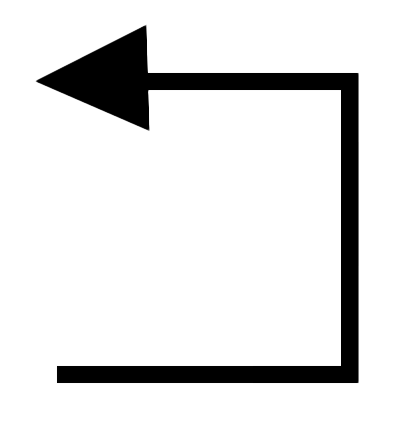}}}
\newcommand\NES{\raisebox{-0.5ex}{\includegraphics[width=.8em]{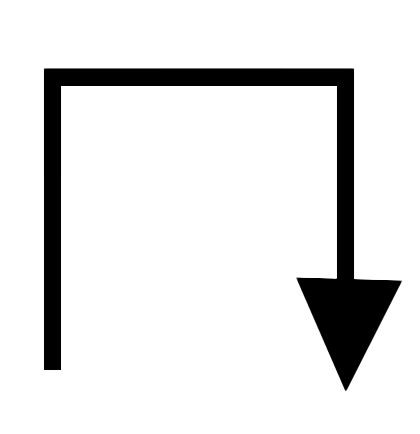}}}

\newcommand\UPRT[2]{{}_{\scalebox{0.5}{$#1$}}\uprightarrow^{\scalebox{0.5}{$#2$}}}
\newcommand\RTUP[2]{{}_{\scalebox{0.5}{$#1$}}\!\rightuparrow^{\, \scalebox{0.5}{$#2$}}}
	
\newcommand\PATHS{\mathsf{Paths}}
\newcommand\RECT{\mathsf{Rect}}

\newcommand\bottomboundary[1]{\stackrel{\scalebox{0.6}{#1}}{\rightarrow}}
\newcommand\leftboundary[1]{%
  \mathchoice
    {\uparrow\kern-0.4em\scalebox{0.6}{#1}} 
    {\uparrow\kern-0.4em\scalebox{0.6}{#1}} 
    {\uparrow\kern-0.0em\scalebox{0.6}{#1}} 
    {\uparrow\kern-0.1em\scalebox{0.6}{#1}} 
}

\newcommand\Paths{\mathsf{Paths}}
\newcommand\Rect{\mathsf{Rect}}
\newcommand\Lasso{\mathsf{Lasso}}

\newcommand{\bfh}{\mathbf{h}}
\newcommand{\bfw}{\mathbf{w}}
\newcommand{\bfa}{\mathbf{a}}
\newcommand{\bfb}{\mathbf{b}}
\newcommand{\bfe}{\mathbf{e}}

\newcommand{\bfH}{\mathbf{H}}
\newcommand{\bfW}{\mathbf{W}}

\newcommand{\freeG}{\mathrm{G}}
\newcommand{\freeH}{\mathrm{H}}

\newcommand{\proj}{\mathbf{p}}

\usepackage{forest}
\forestset{
  decor/.style = {
    label/.expanded = {[inner sep = 0.2ex, font=\unexpanded{\tiny}]right:{$#1$}}
  },
  root/.style = {minimum size = 0.1ex},
  decorated/.style = {
    for tree = {
      circle, fill, inner sep = 0.3ex, minimum size = 1.ex,
      grow' = south, l = 0, l sep = 1.2ex, s sep = 0.7em,
      fit = tight, parent anchor = center, child anchor = center,
      delay = {decor/.option = content, content =}
    }
  },
  default preamble = {decorated, root},
  begin draw/.code={\begin{tikzpicture}[baseline={([yshift=-0.5ex]current bounding box.center)}]},
}

\declaretheorem[numberwithin=section]{theorem}
\declaretheorem[sibling=theorem]{proposition, corollary, lemma}
\declaretheorem[style=definition, sibling=theorem]{definition, convention, example, remark}

\title{A multiplicative surface signature through its Magnus expansion}

\author[1,2]{Ilya Chevyrev}
\author[3]{Joscha Diehl}
\author[4]{Kurusch Ebrahimi-Fard}
\author[5,6]{Nikolas Tapia}

\affil[1]{{\footnotesize Institut f\"ur Mathematik, TU Berlin,
10623 Berlin, Germany}}
\affil[2]{{\footnotesize School of Mathematics, University of Edinburgh, EH9 3FD, UK}}
\affil[3]{{\footnotesize Institute of Mathematics and Computer Science, University of Greifswald, 17489 Greifswald, Germany}}
\affil[4]{{\footnotesize Department of Mathematical Sciences, Norwegian University of Science and Technology,\linebreak 7034 Trondheim, Norway}}
\affil[5]{{\footnotesize Weierstrass Institute, 10117 Berlin, Germany}}
\affil[6]{{\footnotesize Institut für Mathematik, HU Berlin, 12489 Berlin, Germany}}

\mleftright
\begin{document}

\maketitle

\date{\today}

\begin{abstract}
  \small
In the last decade, the concept of path signature has achieved significant success in data science applications. It offers a powerful set of features that effectively capture and describe the characteristics of paths or sequential data. This is partly explained by the fact that the signature of a path can be computed in linear time, using a dynamic programming principle based on Chen’s identity. The path signature can be viewed as a specific example of a product or time-/path-ordered integral. In other words, it represents a one-parameter object built on iterated integrals over a path.

Defining a signature over surfaces requires considering iterated integrals over these surfaces, effectively introducing an additional parameter, resulting in a two-parameter signature. This extended signature is intrinsically connected to a non-commutative generalization of Stokes’ theorem, which is fundamentally connected to the concept of crossed modules of groups. The latter provides a well-established framework in higher gauge theory, where crossed modules with feedback maps exhibiting non-trivial kernels, combined with multiparameter iterated integrals, play a pivotal role.

Building on Kapranov’s work, we explore the surface analog of the log-signature for paths by introducing a Magnus-type formula for the logarithm of the surface signature. This expression takes values in a free crossed module of Lie algebras, defined over a free Lie algebra. We furthermore prove a non-commutative sewing lemma applicable to the crossed module setting and give a definition of rough surface in the so-called Young--H\"older regularity regime along with a corresponding continuous extension theorem. This approach enables the analysis and computation of surface features that go beyond what can be expressed by computing line integrals along the boundary of a surface. 

\end{abstract}

\newpage

\tableofcontents

\newpage

%
%
%
%
%
%
%
%
%
%
%
%
%
%


\section{Introduction}
\label{sec:intro}

The notion of path signature originated in K.-T.~Chen's work on iterated integrals and homotopy theory (see, for example, \cite{chen1957integration,chen1977iterated}), and rose to prominence in the context of stochastic analysis through the work of T.~Lyons on rough paths \cite{Lyons98}.
The path signature consists of a sequence of values computed through iterated integrals along a path. It provides a powerful tool for extracting meaningful information about the path. 

Recently, several authors have promoted the practical utility of the concept of path signature in data science, where it offers features that encapsulate information encoded in sequential data, i.e., chronologically ordered sequences of data points such as, for instance, time-series in finance \cite{cass2024lecture,chevyrev2016primer,lyons2024signature}. Its computational effectiveness derives, in part, from the ability to calculate the path signature in linear time, which is made possible by a dynamic programming approach rooted in Chen's identity (see \eqref{eq:Chen1} below as well as \cite{DEFT22}). 

The main aim of the work at hand is to extend the concept of signature from paths to surfaces which amounts to considering some form of iterated integrals along surfaces. Several attempts at this challenging task have been made recently \cite{diehl2024signature,Lee_Oberhauser_23_surfaces}.  Our approach is based on the categorical framework of higher gauge theory where iterated integrals along surfaces play an important role \cite{BaHu2011,BaSchr07}. More precisely, we utilize the concept of crossed modules of groups and Lie algebras, as detailed in Wagemann's textbook \cite{Wage2021}.

Our work draws from the seminal contribution of Kapranov \cite{kapranov2015membranes}, where a comprehensive description
of the algebraic structure involved, i.e., free crossed modules of Lie algebras, is given.
The latter replaces the free Lie algebra underlying the path signature case, thereby providing the Lie algebraic framework for describing the logarithm of the surface signature. We will motivate this structure in \Cref{ss:why} and provide precise description of this analog of the free Lie algebra in \Cref{s:free_crossed_module}.
This can be seen as a multiparameter generalization of the Magnus expansion. We compute its first few orders explicitly to better understand its components and their relevance for describing features of the underlying surface. Regarding the classical Magnus expansion and its relation to the path signature, we refer to \Cref{sec:prelim} for a brief recap. 

\medskip

Our first definition of surface development (or its Magnus series respectively log-signature) is based on classical path development and is close, to a large extent, to the methods of \cite{SoncZucc15,Lee_Oberhauser_23_surfaces}.
While this construction relies only on classical tools from controlled differential equations,
it has the drawback that it seemingly breaks the natural coordinate invariance of the problem (essentially by treating the $x$-axis and $y$-axis differently), and only a posteriori one shows some form of coordinate invariance.

One of our main contributions is to provide another construction of the surface development that does not break coordinate invariance.
This construction starts from a ``germ'' and applies a 2D non-commutative sewing lemma inspired by rough path theory that we state and prove in \Cref{sec:sewing}.
We emphasize that this construction is independent of the first one (via controlled differential equations), but the two agree in the case when both are well-defined.
However, our construction via sewing
applies in weaker regularity regimes.
This fact naturally leads us to the definition of a ``rough surface'' (\Cref{def:RS})
and to a corresponding continuous extension theorem (\Cref{thm:RS_ext}).

In \Cref{subsec:Young} we show how a surface $X\colon[0,1]^2\to\R^n$ leads to a rough surface in the so-called Young--H\"older regularity regime.
This regime agrees with the Young regime of \cite{Lee_Oberhauser_23_surfaces} when measuring regularity in terms of rectangular increments.

\medskip

\textbf{Contributions}

\begin{enumerate}

	\item
	In \Cref{ss:kapranov}, we provide a detailed elaboration of Kapranov's algebraic structure, including several examples.
	Although probably known to experts, we also provide a proof of the freeness of his crossed module of Lie algebras in \Cref{ss:freeness}.

	\item
	In \Cref{sec:logsigsuf}, we define the 
	surface signature by its Magnus expansion.
	For calculations in the Lie algebra, we spell out Chen's identity on that level, \Cref{thm:log_signature}.
	In \Cref{ex:general_case}, we give explicit expressions of the first few levels of the surface Magnus expansion, 
	in addition to multiple consistency checks.
	In particular, we emphasize that our perspective generalizes Magnus’ expansion to the setting of two-parameter objects.
	We provide a universality statement: every surface development factors through the surface log-signature, \Cref{thm:universality}
	(this can be considered a special form of the log-ODE method
	\cite{morrill2020neural}, for linear ODEs).
	

	\item
	In \Cref{sec:sewing}, we provide a 2D non-commutative sewing lemma which we use to define surface developments.
	We furthermore give a definition of a ``rough surface'' along with a corresponding continuous extension theorem 
	(i.e.~the existence of a (log-)signature). 

	\item
	In \Cref{ss:why}, we provide a novel heuristic motivation for the use of a crossed module structure within the context of surface signatures.
\end{enumerate}

\textbf{Related works}

Recently, there have been several proposed generalizations of the path signature to surfaces:
Reference \cite{diehl2022two} extrapolates the time-warping invariance of the iterated-sums signature \cite{diehl2020time}
to two parameters;
reference \cite{giusti2022topological} extrapolates the mapping space viewpoint of Chen to multiple parameters;
reference \cite{zhang2022two} uses a specific set of integrals for texture classification;
reference \cite{diehl2024signature} extrapolates the ``integration over ordered tuples''-viewpoint of Chen to surfaces.
All these approaches have the shortcoming of \emph{not} satisfying a meaningful Chen's identity
and hence pose computational challenges.


Surfaces can be considered as ``paths in the space of paths'',
so they are naturally related to higher category theory.
The two, in many aspects equivalent, main approaches 
are formulated in terms of $2$-category, or globular, and double category, or cubical, concepts.

The $2$-category approach to surface development, and its relation to higher gauge theory, was constructed in
\cite{BaSchr07,BaHu2011,schreiber2011smooth}, see also \cite{girelli2004higher,SoncZucc15,parzygnat2019two}.
The algebraic structure underlying our work is laid out in Kapranov's 2015 preprint
\cite{kapranov2015membranes}.
We note also the work \cite{breen2005differential}, 
where the term \emph{fake curvature} seems to originate from,
and the use of a formal homology connection in \cite{kohno2020higher}.

\medskip

Regarding the algebraic and categorical foundation of the double categorical approach to surface development,
much work was done by Brown and Higgins.
We mention \cite{brown1981equivalence} for an early reference and \cite{brown2011nonabelian} for a more recent textbook treatment.
See also \cite{F2016}.
For the surface development itself, see \cite{martins2008cubical,SoncZucc15},
and more recently \cite{Lee_Oberhauser_23_surfaces} (which also covers certain non-smooth surfaces).

\medskip

There have been more ``hands-on'' approaches regarding the Lie algebraic formulation of a non-commutative Stokes theorem,
\cite{chacon1991stokes},
\cite{radford1998local},
\cite{ZadMans2022},
but they only consider the so-called ``Schlesinger case'' in which the surface integrals --at the Lie algebra level-- are in fact determined by boundary integrals.
The exception is the work of Yekutieli 
\cite{yekutieli2016nonabelian}, which defines ``multiplicative surface integrals''
in a general (crossed module) setting.


During the preparation of this paper, we learned of an independent work \cite{lee2024_2Dsigs} that achieves several similar results as this article, but from a different and complementary viewpoint. Algebraically, \cite{lee2024_2Dsigs} constructs and works primarily in the free crossed module of \emph{associative algebras} and embeds Kapranov's construction of the free crossed module of Lie algebras into this structure.
In contrast, we work directly in the framework of Lie algebras and define surface developments via their Magnus expansion.
Therefore, \cite{lee2024_2Dsigs} constructs a more direct analogue of the \emph{signature}, while our method construct directly the surface analogue of the \emph{log-signature}.

Analytically, both works provide an extension theorem via sewing-type methods for (potentially rough) surfaces.
The main difference between our results is that \cite{lee2024_2Dsigs} considers decay estimates for the signature terms, which is similar to the original extension theorem of \cite{Lyons98}.
This decay estimate implies that the signature is `universal' in the sense that the surface development\footnote{Sometimes called surface holonomy.} for affine connections factors through the signature.
In contrast, we do not aim to control the size of terms as we send the truncation level $N\to\infty$,
and this is anyway not possible in the Lie algebra setting (log-signature terms do not decay factorially in general),
but consequently our analysis appears simpler.
We also provide a `universality' statement for the surface log-signature, but which requires either a small domain or nilpotency assumption.

\bigskip


{\textbf{Acknowledgements}}: 

IC gratefully acknowledges support from the EPSRC via the New Investigator Award EP/X015688/1
and the Deutsche Forschungsgemeinschaft (DFG, German Research  Foundation) CRC/TRR 388 ``Rough Analysis, Stochastic Dynamics and Related Fields" through a Mercator Fellowship held at TU Berlin.
JD was supported by a Trilateral ANR-DFG-JST AI program, DFG Project 442590525. KEF is supported by the Research Council of Norway through project 302831 ``Computational Dynamics and Stochastics on Manifolds'' (CODYSMA). 
IC, JD, and KEF were supported by the project ``Signatures for Images'' at the Centre for Advanced Study at the Norwegian Academy of Science and Letters in Oslo, Norway, during the academic year 2023/24.
NT was funded by the DFG under Germany's Excellence Strategy – The Berlin Mathematics  Research Center MATH+ (EXC-2046/1, project ID: 390685689) and CRC/TRR 388, Project B01.
We thank W. Steven Gray for his help obtaining a copy of the reference \cite{reutenauer1990dimensions}.


\section{Preliminaries: log-signature and Magnus expansion}
\label{sec:prelim}

In this paper, we study the Magnus expansion of the ``signature of a surface''.
\phantomsection
\label{first-appearance-n}
A surface for us will be a (smooth enough) function
$X\colon [0,1]\to \R^n$ or $X\colon \R^2 \to \R^n$.
To put our main contribution in perspective, we first recall that the algebraic framework underlying the concept of path signature is based on the free Lie algebra $\FL(\mathbb{R}^{n})$ generated over $\mathbb{R}^{n}$ and its corresponding universal enveloping algebra, which is understood to be the tensor algebra
\begin{equation*}
	\text{T}(\mathbb{R}^{n})\coloneqq\bigoplus_{k=0}^{\infty}(\mathbb{R}^{n})^{\otimes k},
\end{equation*}
with its usual associative and non-commutative algebraic tensor product $\otimes$ defined over $\mathbb{R}$ and with the empty word $\mathbf{1}$ as its unit ($(\mathbb{R}^{n})^{\otimes 0}\coloneqq\mathbb{R}\mathbf{1}$) -- it is, in fact, a graded connected cocommutative Hopf algebra \cite{reutenauer1993free}. Its elements are denoted as words $w={i_1} \cdots {i_m} \in \text{T}(\mathbb{R}^{n})$, of length $m\eqqcolon|w|$ (the empty word has length $|\mathbf{1}|=0$). 

Elements in the space of tensor series
\begin{equation*}
	\text{T}(( \mathbb{R}^n))\coloneqq \prod_{k \ge 0} (\mathbb{R}^n)^{\otimes k}
\end{equation*}
can be seen as linear maps from $\mathrm{T}(\mathbb{R}^{n})$ to $\mathbb{R}$ written as formal series
\begin{equation}
\label{tensorseries}
	F = \sum_{k \ge 0} \sum_{1 \le j_1, \ldots, j_k \le n} \alpha^{(j_1 \cdots j_k)} \Z{j_1} \cdots \Z{j_k}
\end{equation}
in monomials of the canonical basis $\{\Z{j}\}_{j=1}^n$ of $\mathbb{R}^n$. The natural pairing between $\text{T}(( \mathbb{R}^n))$ and $\text{T}(\mathbb{R}^{n})$ is used such that for $F \in \text{T}(( \mathbb{R}^n))$ and a word ${i_1} \cdots {i_m} \in \text{T}(\mathbb{R}^{n})$ one finds
$$
	F({i_1} \cdots {i_m}) 
	=\langle F, {i_1} \cdots {i_m}\rangle 
	= \alpha^{(j_1 \cdots j_m)}  \in \mathbb{R}.
$$ 
Note that the tensor product in $ \text{T}(\mathbb{R}^{n})$ can be extended to $\text{T}((\mathbb{R}^{n}))$ making it a unital, associative and non-commutative algebra, with unit  $\mathbf{1}$.

Consider an $ \mathbb{R}^n$-valued Lipschitz path defined over the interval $[0,T]$
\begin{equation*}
	\eta_t = \sum_{i=1}^n \eta^{(i)}_t \Z{i}
\end{equation*}
and the linear initial value problem in $\text{T}(( \mathbb{R}^n))$ 
\begin{equation}
\label{SigODE}
	dS_{s,t} = S_{s,t} \otimes d\eta_t, \quad S_{s,s} = \mathbf{1},
\end{equation}
driven by the pullback 
\begin{equation}
\label{const1form}
	d\eta_t = \eta_t^* a = \sum_{i =1}^n \Z{i}d\eta^{(i)}_{t}  = \sum_{i =1}^n \dot{\eta}^{(i)}_t \Z{i} dt,
\end{equation}
of the Lie algebra valued constant one-form on $\mathbb{R}^{n}$
\begin{equation*}
	a = \sum_{i=1}^n \Z{i} \diff{x}^i \in \Omega^{1}(\mathbb{R}^{n}, \FL(\mathbb{R}^{n})).
\end{equation*}
The solution $S=S(\eta)$ of \eqref{SigODE} is called path signature and defines a 1-parameter family of functionals defined with respect to paths in $\mathbb{R}^n$ 
\begin{equation}
\label{eq:pathsignature}
	S_{s,\cdot}(\eta)\colon [0,T] \to \text{T}(( \mathbb{R}^n)).
\end{equation}

Seen as a function sending words $w \in \text{T}(\mathbb{R}^{n})$ to $\langle S_{s,t}(\eta), w \rangle \in \mathbb{R}$, the path signature has a representation as a sequence of numbers \eqref{tensorseries} in $\mathrm{T}(( \mathbb{R}^n))$ by setting
\begin{equation}
\label{SigSeries}
	S_{s,t}(\eta) = \mathbf{1} + \sum_{k>0} \sum_{1 \le i_1, \ldots, i_k \le n} \langle 
					S_{s,t}(\eta), i_1 \cdots i_k \rangle\, {\Z{i_1}} \cdots {\Z{i_k}},
\end{equation}
where the coefficients in \eqref{SigSeries} are given in terms of Chen iterated integrals. See \cite{StFlour07} for a review. 
The path signature satisfies Chen's identity, that is, it has the $1$-cocycle property\footnote{We follow terminology used by Soncini and Zucchini \cite{SoncZucc15} in their formulation of higher parallel transport, which suggests to identify the path signature as a particular example of the concept of time- or path-ordered exponential characterised by the 1-cocycle property -- in the same spirit the surface signature, as we will see later, can be seen as a 2-cocycle.} with respect to concatenation of paths
\begin{equation}
\label{eq:Chen1}
	S(\eta_1 \eta_2) = S(\eta_1) S(\eta_2).
\end{equation}

As the unique solution to a linear initial value problem, it can be shown that the path signature \eqref{SigSeries} is a 1-parameter family of elements in the particular subset $\freeG((\mathbb{R}^{n})) \subset \text{T}((\mathbb{R}^{n}))$ of so-called group-like elements. In fact, $\freeG((\mathbb{R}^{n}))$ forms a non-commutative (infinite dimensional) group with respect to the tensor product in $\text{T}((\mathbb{R}^{n}))$. Its elements can be expressed as tensor-exponentials
\begin{equation*}
	\freeG((\mathbb{R}^{n})) = \exp^\otimes\bigl(\FL((\mathbb{R}^n))\bigr).
\end{equation*}
The linear space $\FL((\mathbb{R}^n))$ consists of elements in $\text{T}((\mathbb{R}^{n}))$ which can be written as an infinite sum $\sum_{k \geq 1} p_k$ of Lie polynomials $p_k \in \FL(\mathbb{R}^{n})$
(of increasing minimal degree).

With the group structure in $\freeG((\mathbb{R}^{n}))$ available, we may return to \eqref{SigODE} and define the (left) logarithmic derivation
\begin{equation}
\label{logSigODE}
 	\mathcal{L}_t(S_{s,\cdot}) \coloneqq S^{-1}_{s,t} \otimes dS_{s,t} =  d\eta_t.
\end{equation}
Following \cite{BCOR2008,Magnus54}, the group-valued solution to \eqref{logSigODE} can be expressed in terms of  the Lie series $\mag_{s,u}$ resulting as the solution to the Magnus ODE  
\allowdisplaybreaks
\begin{align}
\label{MagnusODE}
\begin{aligned}
	d\mag_{s,u}
	= 	\frac{ \text{ad}_{\mag_{s,u}}}{1 - e^{-\text{ad}_{\mag_{s,u}}}}( d\eta_u )
	=	\sum_{k \ge 0} \frac{B_k}{k!} \text{ad}^{(k)}_{ \mag_{s,u}} d\eta_{u}.
\end{aligned}
\end{align}

Note that, to avoid confusion, though diverging from standard notation, we use $\mag=\mag(\eta)$ for denoting the Magnus expansion in the usual, i.e., path setting, while reserving $\Omega$ for its application in the surface setting.

The second equality in \eqref{MagnusODE} derives from the well-known generating series 
$$
	\frac{z}{1-e^{-z}}=\sum_{k \ge 0} B_k\frac{z^k}{k!}, 
$$
where $\{B_m\}_{m\ge 0}$ are the Bernoulli numbers: $B_0=1$, $B_1=\frac{1}{2}$, $B_2=\frac{1}{6}$, \ldots, and $B_{2k+1}=0$, for $k>0$. Integrating \eqref{MagnusODE} defines a particular Lie series known as the Magnus expansion\footnote{At this point we do not worry about existence and convergence of the Magnus expansion; but note that no issues arise in the free case considered here.}
\allowdisplaybreaks
\begin{align*}
	\mag_{s,t}
	&= \int_s^t d\eta_{u} 
		+ \frac{1}{2}  \int_s^t \int_s^{u_1} [d\eta_{u_2},d\eta_{u_1}]
	+ \frac{1}{4} \int_s^t \int_s^{u_1} \int_s^{u_2} \big[ [d\eta_{u_3},d\eta_{u_2}] ,d\eta_{u_1}\big] \\
	&\qquad + \frac{1}{12} \int_s^t \int_s^{u_1} \int_s^{u_1} \big[d\eta_{u_3},[d\eta_{u_2},d\eta_{u_1}]\big] + \cdots.
\end{align*}
Using \eqref{const1form}, we obtain the following expansion\footnote{To simplify notation, we employ Einstein's summation convention implying summation over repeated (upper and lower) indices.}
\allowdisplaybreaks
\begin{align*}
	\mag_{s,t}
	&= \int_s^t d\eta^{(i)}_{u} \Z{i}
		+ \frac{1}{2}  \int_s^t \int_s^{u_1} d\eta^{(i_2)}_{u_2} d\eta^{(i_1)}_{u_1} [\Z{i_2},\Z{i_1}]
	+ \frac{1}{4} \int_s^t \int_s^{u_1} \int_s^{u_2} d\eta^{(i_3)}_{u_3} d\eta^{(i_2)}_{u_2}d\eta^{(i_1)}_{u_1}	
	\big[ [\Z{i_3},\Z{i_2}\big] ,\Z{i_1}] \\
	&\qquad\  + \frac{1}{12} \int_s^t \int_s^{u_1} \int_s^{u_1}d\eta^{(i_3)}_{u_3} d\eta^{(i_2)}_{u_2} d\eta^{(i_1)}_{u_1}	
	\big[ \Z{i_3}, [ \Z{i_2},\Z{i_1} ]\big] + \cdots.
\end{align*}

In summary, seen as an element in $\freeG((\mathbb{R}^{n})) \subset \text{T}((\mathbb{R}^{n}))$, the path signature can be expressed as a tensor exponential of the Magnus expansion, more commonly known as the log-signature 
\begin{align}
	S_{s,t}(\eta) 
	&= \mathbf{1} + \sum_{k>0} \sum_{1 \le j_1, \ldots, j_k \le n} \langle 
					S_{s,t}(\eta), j_1 \cdots j_k \rangle\, {\Z{j_1}} \cdots {\Z{j_k}} 	\label{SigSeries2}  \\
	&=\exp^\otimes\big( \mag_{s,t}(\eta) \big) 									\label{Log-MagSeries2} .
\end{align}

According to \cite{cass2024lecture}, the coefficients of the path signature \eqref{SigSeries2} can be viewed as a ``library" of features that describe the path $\eta$ at the group level. Uniqueness results from \cite{BGLY2016} indicate that these features, in fact, characterize the path (up to thin homotopy, or tree-like, equivalence, and translation). From \eqref{Log-MagSeries2}, on the other hand, we see that this information is equally captured by the description of its logarithm represented as an element of $\FL((\mathbb{R}^{n}))$ in terms of the Magnus expansion. Indeed, at the Lie algebra level path features are identified with the coefficients in front of the iterated Lie brackets in the Lie series $\mag_{s,t}(\eta)$. We note, however, that there is redundancy in the feature set \eqref{SigSeries2} coming from, for instance, integration by parts, e.g., $\langle S_{s,t}(\eta), i_1  \rangle\, \langle S_{s,t}(\eta), i_2  \rangle = \langle S_{s,t}(\eta), i_1 i_2 \rangle + \langle S_{s,t}(\eta), i_2 i_1 \rangle $. At the Lie algebra level \eqref{Log-MagSeries2} this redundancy is eliminated in Magnus' Lie series $\mag_{s,t}(\eta)$.
This comes with the drawback that Chen's indentity now has to be formulated in terms of the Baker--Campbell--Hausdorff series which poses a significant computation challenge.
In summary, choosing to work with either representation involves considering a trade-off between high dimensionality (path signature) versus computational complexity (Magnus series).

{\textbf{Outlook}} From this perspective, the paper at hand aims to extend the concept of log-signature seen as a ``library'' of (redundancy-free) path features to a library of surface features. To achieve this, we need to ``enlarge our library'' by considering iterated Lie brackets in the letters $\{\Z{i}\}_{i=1}^n \cup \{\Z{ij}\}_{i,j=1}^n$, where the $\Z{ij}$s should be interpreted as unit surface elements. For comparison we list the surface analogue of the Magnus series up to order three:
\begin{align*}
		\Omega^{(1)}_{s_1, s_2;t_1,t_2}  &= 0 \\[0.2cm]
		\Omega^{(2)}_{s_1, s_2;t_1,t_2}  &= \sum_{i<j} \int_{s_1}^{t_1} dr \int_{s_2}^{t_2} dp\ J^{(ij)}_{r,p} \Z{ij} \\[0.2cm]
		\Omega^{(3)}_{s_1, s_2;t_1,t_2}  &= \sum_i \sum_{j<k} \int_{s_1}^{t_1} dr \int_{s_2}^{t_2} dp\ (X^{(i)}_{r,p} - X^{(i)}_{s_1,s_2} )J^{(jk)}_{r,p} [\Z{i}, \Z{jk}],
\end{align*}
where $X\colon \mathbb{R}^2 \to \mathbb{R}^n$ is a smooth enough surface and the Jacobi minors 
\begin{align*}
	J^{(jk)}_{r,p}
	\coloneqq
	\partial_1 X^{(j)}_{r,p}
	\partial_2 X^{(k)}_{r,p}
	-
	\partial_1 X^{(k)}_{r,p}
	\partial_2 X^{(j)}_{r,p}.
\end{align*}
The multi-parameter index $({s_1, s_2;t_1,t_2})$ describes the lower left and upper right corners of a rectangle inside $\mathbb{R}^2$. For details, we refer to the later sections in the paper.


%


\section{Cocycles}
\label{s:cocycles}

We first introduce notation and terminology regarding paths and rectangles, which will be ``indexing objects'' for cocycles.
%
In \Cref{ss:cocycles_lie_group} we start by recalling the (one-parameter) notion of cocycles in Lie groups.
These correspond to path developments.
Then we introduce the notion of $2$-cocycles and give two constructions as ``surface developments''
of them: an ODE construction from the literature and a novel ``sewing'' construction.
In \Cref{ss:cocycles_lie_algebra}
we expand the ODE solutions to both $1$-cocycles and $2$-cocycles
in Magnus expansions in the Lie algebra.
This approach enables us to directly apply Kapranov's construction of the free crossed module of Lie algebras (\Cref{s:free_crossed_module}) to build the logarithm of a surface signature, effectively circumventing the need to construct the group itself.

\begin{remark}
The term “cocycle” originates from the geometric perspective of holonomy in principal bundles, a viewpoint we adopt here to align with the literature \cite{SoncZucc15}. However, this perspective will not be used in this paper. Readers may instead interpret “cocycle” as synonymous with “multiplicative map” within the context of either 1D or 2D frameworks.
\end{remark}


Recall (\cite{barrett1991holonomy,caetano1994axiomatic}) that two $\R^n$-valued paths $\gamma, \eta\colon [0,1] \to \R^n$ are \DEF{thinly homotopic} if they are related under the equivalence relation generated by the following relation: 

$\gamma$ is related to $\eta$ if there exists a piecewise smooth map $H\colon [0,1]^2 \to \R^n$ such that
	\begin{align*}
		H(t,0) = \gamma_t, \quad H(t,1) = \eta_t, \quad
		H(0,s) = H(1,s) = \gamma_0 = \eta_0,
	\end{align*}
	and $H([0,1]^2)$ is contained in the 
	 image of the union of the images of $\gamma$ and $\eta$.
For piecewise smooth paths, thin homotopic equivalence is equivalent to \textbf{tree-like equivalence} \cite{HL10, BGLY2016}.

Concatenation of (equivalence classes of) paths written as $\gamma \eta = \gamma \sqcup \eta$, is well-defined as soon as the endpoint of $\gamma$ coincides with the starting point of $\eta$. As a result, the set $\PATHS_n$ of thin-homotopy classes of piecewise smooth paths in $\R^n$ forms a groupoid. 

For later use, we define the following elements of $\PATHS_2$:
		\begin{itemize}
			\item 
			$\UP{v_1,v_2}{v_1,q}$: the linear path from $(v_1,v_2)$ to $(v_1,q)$ (usually $v_2 \le q$),

			\item
			$\DN{v_1,q}{v_1,v_2}$: the linear path from $(v_1,q)$ to $(v_1,v_2)$ (usually $v_2 \le q$),

			\item
			$\RT{v_1,v_2}{r,v_2}$: the linear path from $(v_1,v_2)$ to $(r,v_2)$ (usually $v_1 \le r$),

			\item 
			$\LT{r,v_2}{v_1,v_2}$: the linear path from $(r,v_2)$ to $(v_1,v_2)$ (usually $v_1 \le r$),

			\item
			$\UPRT{v_1,v_2}{r,q} \coloneqq \UP{v_1,v_2}{v_1,q} \sqcup \RT{v_1,q}{r,q}$,
			\item
			$\RTUP{v_1,v_2}{r,q} \coloneqq \RT{v_1,v_2}{r,v_2} \sqcup \UP{r,v_2}{r,q} $.
		\end{itemize}
\newcommand\INTERVAL[2]{\RT{#1}{#2}}
For $a,b\in\R$ we denote the unique (up to thin homotopy) element
of $\PATHS_1$ starting at $a$ and ending at $b$ by
	$\INTERVAL{a}{b}$.

We will also need the concept of ``rectangles'' in $\R^2$. Define the set 
	\begin{align*}
		\RECT \coloneqq
		\{ (e^-,e^+) \in \R^2\times\R^2 \mid e^-_1 \le e^+_1, e^-_2 \le e^+_2 \},
	\end{align*}
where $(e^-,e^+)$ is considered the bottom-left and top-right corner of a rectangle. General elements of $\RECT$ will sometimes be denoted by the symbol $\square$. If we want to specify the corners of a rectangle, we will also write $\SQUARE{e^-}{e^+}$. The area of a rectangle $\square$ is denoted by~$|\square|$.

Let $\square_A, \square_B \in \RECT$ be two rectangles. If the right edge of $\square_A$ coincides with the left edge of $\square_B$, then we can concatenate them horizontally, to obtain a new rectangle $\square = \square_A \square_B$, whose bottom-left corner is that of $\square_A$ and whose top-right corner is that of $\square_B$. Analogously, if the top edge of $\square_A$ coincides with the bottom edge of $\square_B$, then we can concatenate them vertically, to obtain a new rectangle $\square = \genfrac{}{}{0pt}{1}{\square_B}{\square_A}$, whose bottom-left corner is that of $\square_A$ and whose top-right corner is that of $\square_B$.

As boundaries of rectangles are crucial for our computations with crossed modules, we define the boundary map
	\begin{align}\label{eq:boundary_rect_paths}
		\partial\colon \RECT \to	\PATHS_2,
	\end{align}
which sends a rectangle to (the thin-homotopy class of) its boundary path, starting at the bottom-left corner and going counter-clockwise. For example, if $\square = (e^-,e^+)$, then 
$$
	\partial(\square) = \RT{e^-_1,e^-_2}{e^+_1,e^-_2} \sqcup \UP{e^+_1,e^-_2}{e^+_1,e^+_2} \sqcup \LT{e^+_1,e^+_2}{e^-_1,e^+_2} \sqcup \DN{e^-_1,e^+_2}{e^-_1,e^-_2}.
$$ 
%


\subsection{Cocycles in (Lie) groups}
\label{ss:cocycles_lie_group}

\begin{definition}
	\label{def:1cocycle}
	Let $G$ be a group.
	A \textbf{1-cocycle} is a functor $\PD\colon \PATHS_n \to G$.
	Spelled out, this means that $\PD$ is a function from equivalence classes of paths to $G$ and
	for any composable (equivalence classes of) paths $\gamma_1, \gamma_2 \in \PATHS_n$,
	\begin{align}
	\label{eq:1cocycle}
		\PD(\gamma_1 \gamma_2) = \PD(\gamma_1) \PD(\gamma_2).
	\end{align}
\end{definition}

Let $G$ be a Lie group. It is well known that (smooth enough) 1-cocycles are in one-to-one correspondence with (smooth enough) $1$-forms. For general results see \cite{schreiber2009parallel}, but we mention the following statement.

\begin{theorem}[{\cite[Proposition 2.4]{SoncZucc15}, Path development}]
	\label{thm:1cocycleIntegration}
	Let $a \in \Omega^1(\R^n, \mathfrak g)$ be a $\mathfrak g$-valued 1-form where $\mathfrak g$ is the Lie algebra of a Lie group $G$.
	Let $\gamma\colon [0,1] \to \R^n$ be a continuous, piecewise smooth path.
	Then the ODE
	\begin{align*}
		\LL_t( Y ) = (\gamma^* a)_t, \qquad Y_0 = \mathbf{1}_G,
	\end{align*}
	has a unique solution $Y$. Here $\LL_t$ is the \DEF{(left) logarithmic derivative}, defined by 
	\begin{align*}
		\LL_t(Y) \coloneqq (L_{(Y_t^{-1})})_* dY_t,
	\end{align*}
	where $L_g$ is left-multiplication by $g \in G$.
	$Y_1$ is constant on thin-homotopy classes of paths. 
	One then defines the \DEF{path development}
	$$
		\PD^a(\gamma) \coloneqq Y_1,
	$$
	which is a 1-cocycle.
\end{theorem}

We also have the following alternative construction of $1$-cocycles, which will
follow from the sewing constructions in \Cref{sec:sewing} (in particular
\Cref{lem:1D_sewing}).

\begin{theorem}
	Admit the assumptions of \Cref{thm:1cocycleIntegration},
	parametrize $\gamma \in \PATHS_n$ by $[0,1]$,
	and define
	\begin{align*}
		\widehat\PD([s,t])
		\coloneqq \exp( \int_s^t \gamma^* a ).
	\end{align*}
	Then there exists a ``sewed'' map $\PD$, which is sufficiently close to \(\widehat\PD\),
	and we have $\PD^a(\gamma) = \PD([0,1])$.
\end{theorem}


\begin{example}
	\label{ex:pathsignature}
A natural example of 1-cocycle is the path signature $\IIS$, \eqref{eq:pathsignature}, which satisfies Chen's property \eqref{eq:Chen1}.
Indeed, even though it does not exactly fit the previous two theorems (since $\mathrm{G}((\R^d))$ is not an honest Lie group),
it is a 1-cocycle that satisfies an analogous ODE.
\end{example}

The notion of 2-cocycles is supposed to generalize path development and \eqref{eq:1cocycle} to surface development.
To support this viewpoint, we sketch a concrete computation of the path signature over a closed loop $\eta$ in $\R^n$
and how it relates - via a Stokes'-type formula - to the integral over a rectangle.
The loop is identified with the ``boundary'' of a surface patch $X\colon [0,1]^2 \to \R^n$.\footnote{By a choice of interpolation, any closed loop can be written as the boundary of a surface. Also, note that we use the terms surface and membrane interchangeably.}

Recall the notation of axis-parallel paths from \Cref{s:cocycles} and let
\begin{align*}
	\gamma \coloneqq \gamma^{1,1} \coloneqq \RT{0,0}{1,0} \sqcup \UP{1,0}{1,1} \sqcup \LT{1,1}{0,1} \sqcup \DN{1,1}{0,1},
\end{align*}
which can be thought of as the boundary of the unit square.
Define $\eta \coloneqq X_* \gamma$, the image of $\gamma$ under $X$.
As before, we consider the Lie algebra-valued, constant $1$-form
\begin{align*}
  a = \sum_{i=1}^n \Z{i}\diff{x}^i,
\end{align*}
and its pullback, $\A\coloneqq X^* a = \sum_{i=1}^n \Z{i} \partial_1 X^i \diff{t} + \sum_{i=1}^n \Z{i} \partial_2 X^i \diff{s}$.
Even though $\freeG((\R^n))$ is not a Lie group,
it is well-known (see \Cref{ex:pathsignature})
that the path development, \Cref{thm:1cocycleIntegration},
$\PD^{\A}(\gamma)$,
is well-defined (and coincides with the path signature of $\eta$).
In this case of a closed loop, using physics terminology, it amounts to computing a so-called Wilson-loop \cite{Karpetal1999} in the group $\freeG((\mathbb{R}^n))$ of invertible tensor series.

We now turn $\gamma$ into a 2-parameter family of curves, parametrized by the top-right corner of the rectangle $((0,0),(t_1,t_2))$,
\begin{align*}
	\gamma^{t_1,t_2} = \RT{0,0}{t_1,0} \sqcup \UP{t_1,0}{t_1,t_2} \sqcup \LT{t_1,t_2}{0,t_2} \sqcup \DN{0,t_2}{0,0}.
\end{align*}
A lengthy calculation similar to the Wilson-loop computation presented in \cite{Karpetal1999}  
gives
\begin{align}
	\label{eq:ll_wilson}
	\LL_{t_2}( \PD^{\A}(\gamma^{t_1,\bullet}) )
	&= K_{t_1,t_2} \diff{t}_2,
\end{align}
where $K_{t_1,t_2}=K_{0,0; t_1,t_2}$ is the function
\begin{align}
\label{eq:Kmap}
	K_{t_1, t_2}
  \coloneqq
	\int_0^{t_1}
	\PD^{\A}( \UPRT{0,0}{r,t_2} )
	\left( \partial_1 \A^{(2)}_{r,t_2} - \partial_{2} \A^{(1)}_{r,t_2} + [\A^{(1)}_{r,t_2}, \A^{(2)}_{r,t_2}] \right)
	\PD^{\A}( \UPRT{0,0}{r,t_2} )^{-1} dr.
\end{align}

\begin{remark}\label{rmk:wedge}
\begin{enumerate}
\item
Recall that the 2-form\footnote{In the literature this is also denoted $\frac12 [\alpha, \alpha']$ or $\frac12 [\alpha \wedge \alpha']$.}
$$
			\alpha \wedge \alpha'(X,Y) \coloneqq \frac12([\alpha(X),\alpha'(Y)] - [\alpha(Y),\alpha'(X)]).
$$
		Explicitly, for $\alpha = \sum_{i=1}^k A_i \alpha^i$ and $\alpha' = \sum_{j=1}^l A'_j \alpha'^j$ in $\Omega^1(M, \mathfrak g)$, with $A_i, A'_j \in \mathfrak g$, and $\alpha^i, \alpha'^j \in \Omega^1(M)$, we have $\alpha \wedge \alpha'
			= \frac12 \sum_{i,j} [A_i, A'_j] \alpha^i \wedge \alpha'^j 
			= \sum_{i < j} [A_i, A'_j] \alpha^i \wedge \alpha'^j.$

\item
The element
$$
	F_{12} \coloneqq (d\A + \A\wedge \A)(\partial_{t_1}, \partial_{t_2}) =\partial_1 \A^{(2)}_{r,t_2} - \partial_{2} \A^{(1)}_{r,t_2} + [\A^{(1)}_{r,t_2}, \A^{(2)}_{r,t_2}]
$$
is a coefficient of the 2-form $F$ which could be interpreted as curvature form  associated to the 1-form $\A$ (\cite{Karpetal1999}).  

\item
Owing to the fact that $a$ is a constant $1$-form, \eqref{eq:Kmap} simplifies to
	\begin{align*}
		K_{t_1, t_2}
		=
		\int_0^{t_1}
		\PD^{\A}( \UPRT{0,0}{r,t_2} )
		\,[\A^{(1)}_{r,t_2}, \A^{(2)}_{r,t_2}]\,
		\PD^{\A}( \UPRT{0,0}{r,t_2} )^{-1}
		dr.
	\end{align*}
However, we shall carry the more general expression along.
\end{enumerate}
\end{remark}

If we want to stay in the Lie algebra throughout, we can use the Magnus expansion to write 
$$
	\mag^{r,t} \coloneqq
	\log \PD^{\A}(\UPRT{0,0}{r,t}).
$$
Then, using the fact that $\exp(x)y\exp(-x)=e^{\ad_x}y$, we obtain
\begin{align*}
	K_{t_1, t_2}
	\coloneqq
	\int_0^{t_1}
	e^{\ad_{\mag^{r,t_2}}}
	\left[ \partial_1 \A^{(2)}_{r,t_2} - \partial_{2} \A^{(1)}_{r,t_2} + [\A^{(1)}_{r,t_2}, \A^{(2)}_{r,t_2}] \right]
	 dr,
\end{align*}
and \eqref{eq:ll_wilson} can then also be solved via another use of Magnus expansion
\begin{equation}
\label{eq:2magnus}
	\mag^{t_1,t_2}
		= \int_0^{t_2}\frac{\text{ad}_{\mag^{t_1,r}}}{1 - e^{-\text{ad}_{\mag^{t_1,r}}}}( K_{t_1,r} )dr,
\end{equation}
which defines a 2-parameter Magnus series that exponentiates to the following identity in the group $G$ as a result of computing the path signature along a closed path $\gamma^{t_1,t_2}$
$$
	\exp(\mag^{t_1,t_2}) = \PD^{\A}(\gamma^{t_1,t_2}).
$$
Introducing notation for the left-hand side, we find 
\begin{align}
\label{Schlesinger}
\begin{aligned}
	\SD(\SQUARE{0,0}{t_1,t_2})
	\coloneqq\exp(\mag^{t_1,t_2}) 
	&=\PD^{\A}(\gamma^{t_1,t_2}) \\
	&= \PD^{\A}(\RT{0,0}{t_1,0}) \PD^{\A}(\UP{t_1,0}{t_1,t_2})
		\PD^{\A}(\LT{t_1,t_2}{0,t_2})\PD^{\A}(\DN{0,t_2}{0,0}),
\end{aligned}
\end{align}   
which can be summarized in the ``non-commutative Stokes' identity''
\begin{align}
	\label{eq:first_stokes}
	\SD(\square) = \PD^{\A}( \partial \square ).
\end{align}

It is not hard to generalize the above computation to a rectangle $\SQUARE{s_1,s_2}{t_1,t_2}$ with arbitrary corners corners $(s_1,s_2)$, $(s_1,t_2)$, $(t_1,t_2)$, $(t_1,s_2)$. The non-commutative Stokes identity \eqref{eq:first_stokes} permits to show that $\SD$ satisfies the following ``horizontal Chen-identity'', $s_1 < u_1 < t_1$, $s_2 < t_2$, 
\begin{align}
	\label{eq:schlesinger_chenHorizontal}
	\SD( \SQUARE{s_1,s_2}{t_1,t_2} )
	=
	\actionGG_{\PD^{\A}(\RT{s_1,s_2}{t_1,s_2})}\left(\SD(\SQUARE{u_1,s_2}{u_1,t_2}) \right)\, \SD(\SQUARE{s_1,s_2}{u_1,t_2}),
\end{align}
as well as the ``vertical Chen's identity'', $s_1 < t_1$, $s_2 < u_2 < t_2$,
\begin{align}
	\label{eq:schlesinger_chenVertical}
	\SD( \SQUARE{s_1,s_2}{t_1,t_2} )
	=
	\SD( \SQUARE{s_1,s_2}{t_1,u_2} )\, \actionGG_{\PD^{\A}(\UP{s_1,s_2\;\,}{s_1,u_2\;\,} )}(\SD(\SQUARE{s_1,u_2}{t_1,t_2}) ),
\end{align}
where we define $\actionGG_g(h) \coloneqq g h g^{-1}$. For example, using \eqref{eq:first_stokes} and the multiplicative property of the path development, we have
\begin{align*}
	\SD( \SQUARE{s_1,s_2}{u_1,t_2} )
	&=
	\PD^{\A}( \partial \SQUARE{s_1,s_2}{u_1,t_2} ) \\
	&=
	\PD^{\A}( \RT{s_1,s_2}{t_1,s_2}
	\sqcup \partial \SQUARE{t_1,s_2}{u_1,t_2}
	\sqcup \LT{t_1,s_2}{s_1,s_2}
	\sqcup \partial \SQUARE{s_1,s_2}{t_1,t_2} ) \\
	&=
	\actionGG_{\PD^{\A}( \RT{s_1,s_2}{t_1,s_2})}
	\left( \PD^{\A}( \partial \SQUARE{t_1,s_2}{u_1,t_2} ) \right)\,
	\PD^{\A}( \partial \SQUARE{s_1,s_2}{t_1,t_2} ) \\
	&=
	\actionGG_{\PD^{\A}(\RT{s_1,s_2}{t_1,s_2} )}
	\left(\SD( \SQUARE{t_1,s_2}{u_1,t_2} ) \right)\,
	\SD( \SQUARE{s_1,s_2}{t_1,t_2} ).
\end{align*}

Identity \eqref{eq:first_stokes} says that the surface development on the left-hand side is completely determined by the path development along the boundary of the surface $X$, which makes the former less interesting.
However, the pertinent feature here is that identity \eqref{eq:first_stokes}
can be used to ``extrapolate'', by interpreting the $\id$-map on the group $\mathrm{G}((\mathbb{R}^n))$
as a group morphism sending elements in $\mathrm{G}((\mathbb{R}^n))$ specified by surface development along $X$ to elements in $\mathrm{G}((\mathbb{R}^n))$ determined by path development along the boundary $X$:
\begin{equation*}
	\id_{\mathrm{G}((\mathbb{R}^n))}\left( \SD(\square) \right) = \PD^{\A}( \partial \square ).
\end{equation*}
It is also clear that the group $\mathrm{G}((\mathbb{R}^n))$ acts on itself via inner automorphism $\actionGG\colon \mathrm{G}((\mathbb{R}^n)) \to \Aut(\mathrm{G}((\mathbb{R}^n)))$.
This interpretation is summarized in the fact that every group trivially defines a \emph{crossed module} over itself.
It motivates the following definition of 2-cocycles on a general crossed module of groups.

The concept of crossed modules of groups appeared in Whitehead's work \cite{whitehead1949combinatorial}, alongside the notion of crossed morphisms. The simplest example of a crossed module consists of a group $G$ and its group of automorphisms, $\Aut(G)$, along with the canonical homomorphism from the former to the latter. For a comprehensive historical account, we refer to Huebschmann's recent articles \cite{Huebschmann2021, Huebschmann2023}. According to \cite{Huebschmann2021}, the notion of crossed module of a Lie algebra, also known as differential crossed modules, first appeared in the 1982 work by Kassel and Loday  \cite{KasLod1982}.

In this subsection, we briefly recall the definition of crossed modules of groups and Lie algebras. We derive the latter from the former, analogous to the classical relationship between (Lie) groups and Lie algebras (see \Cref{prop:crossedmodule} for the particular example of going from a crossed module of nilpotent Lie algebras to a crossed module of Lie groups). For a comprehensive exposition of the modern theory of crossed modules and its applications in theoretical physics, particularly in the context of higher gauge theory, we refer the reader to Wagemann's textbook \cite{Wage2021} and additional works \cite{BaSchr07,BaHu2011,girelli2004higher,parzygnat2019two,SoncZucc15}.

\begin{definition}\label{def:xmod}
    Let $H$ and $G$ be (Lie) groups.  		
    A \textbf{crossed module of (Lie) groups} $(\feedbackGG\colon H \to G, \actionGG)$
    is a diagram
    \begin{align*}
        H \stackrel{\feedbackGG}\to G \stackrel{\actionGG}\to \Aut(H),
    \end{align*}
    where $\feedbackGG$ and $\actionGG$ are (Lie) group morphisms, satisfying the identities
	(the second is also known as the Peiffer identity)
    \begin{align}
        \feedbackGG( \actionGG_g( h' ) ) &= g \feedbackGG(h') g^{-1} && h' \in H, g \in G \label{eq:equivariant} \tag{EQUI}\\
        \actionGG_{\feedbackGG(h)}( h' ) &= h h' h^{-1}    && h, h' \in H \label{eq:peiffer}. \tag{PEIF}
    \end{align}
\end{definition}

\begin{example}
	\label{ex:group_example}
	Let $G$ be any group,
	$H$ a normal subgroup of $G$, $\feedbackGG$ the inclusion map, and $\actionGG_g(h) \coloneqq g h g^{-1}$.
\end{example}

Taking derivatives in the Lie group case, we obtain the following definition of crossed module of Lie algebras.
\begin{definition}\label{def:diffxmod}
	Let $ \mathfrak h$ and $\mathfrak g$ be Lie algebras. 
    A \DEF{crossed module of Lie algebras}%
	\footnote{Also called a \emph{differential crossed module}.}
	$(\feedbackAA\colon \mathfrak h \to \mathfrak g, \actionAA)$
    is a diagram
    \begin{align*}
        \mathfrak h \stackrel{\feedbackAA}\to \mathfrak g \stackrel{\actionAA}\to \Der(\mathfrak h),
    \end{align*}
    where $\feedbackAA$ and $\actionAA$ are Lie algebra morphisms, satisfying
    \begin{align}
        \feedbackAA( \actionAA_x( y' ) ) &= \left[ x, \feedbackAA(y') \right] && y' \in \mathfrak h, x \in \mathfrak g \label{eq:equiDiff} \tag{EQUID} \\
        \actionAA_{\feedbackAA(y)}( y' ) &= \left[ y, y' \right] && y, y' \in \mathfrak h. \label{eq:peifferDiff} \tag{PEIFD}
    \end{align}
\end{definition}

\begin{example}
	Let $\mathfrak g$ be any Lie algebra, 
	$\mathfrak h \subset \mathfrak g$ a Lie ideal, $\feedbackAA$ the inclusion map, and $\actionAA_x(y) \coloneqq [x,y]$.
\end{example}

\begin{lemma}
	\label{lem:action_exp}
	Let $(\feedbackGG\colon H\to G,\actionGG)$ be a crossed module of Lie groups.
	Let $\mathfrak h,\mathfrak g$ be the Lie algebras of $H,G$ respectively.
	Let $\feedbackAA\colon \mathfrak h \to \mathfrak g$ be the differential of $\feedbackGG$.
	Let $\actionGA_g\colon \mathfrak h \to \mathfrak h$ be the differential of $\actionGG_g$,
	which is an element of $\Aut(\mathfrak h)$.
	Finally, let $\actionAA\colon \mathfrak g \to \Der(\mathfrak h)$ be the differential of $\actionGA$.

	Then $(\feedbackAA\colon \mathfrak h \to \mathfrak g, \actionAA)$
	is a crossed module of Lie algebras.
	Further, for $x \in \mathfrak g$
	\begin{align*}
		\actionGA_{\exp_G( x ) } = \exp( \actionAA_x ),
	\end{align*}
	as automorphisms of $\mathfrak h$.
\end{lemma}

\begin{remark}
Our notation $\actionGG, \actionGA, \actionAA$ makes sense in light of the last lemma.
\end{remark}

\begin{proof}[Proof of \Cref{lem:action_exp}]
	Both statements are elementary, we prove the second one.
	Define the following paths in $\Aut(\mathfrak h)$,
	\begin{align*}
		f_t \coloneqq \exp( t \actionAA_x ), \qquad
		g_t \coloneqq \actionGA_{\exp_G( t x ) }.
	\end{align*}
	Then $f_0 = g_0 = \id$, and $\dot f_t = f_t \actionAA_x$.
	Further
	\begin{align*}
		\dot g_t
		&=
		 \restr{\frac{d}{dh}}{h=0} g_{t+h}
		=
		\actionGA_{\exp_G( t x ) }
		 \restr{\frac{d}{dh}}{h=0}
		\actionGA_{\exp_G( h x ) } \\
		&=
		\actionGA_{\exp_G( t x ) }
		\actionAA_x = g_t \actionAA_x.
	\end{align*}
	Hence, by uniqueness of solutions of linear ODEs, $f_t = g_t$, $t \in \R$.
	Taking $t=1$ gives the result.
\end{proof}

\begin{definition}\label{def:diffxmorph}
	A \textbf{morphism of two crossed modules of Lie algebras}
    $(\feedbackAA\colon \mathfrak h \to \mathfrak g, \actionAA)$ and
    $(\feedbackAA'\colon \mathfrak h' \to \mathfrak g', \actionAA')$,
	is a tuple $(\psi,\phi)$, where
	\begin{align*}
		\psi\colon \mathfrak h \to \mathfrak h', \quad
		\phi\colon \mathfrak g \to \mathfrak g',
	\end{align*}
	are Lie algebra morphisms, satisfying
	\begin{align*}
		\phi \circ \feedbackAA &= \feedbackAA' \circ \psi\\
		\psi( \actionAA_{g}( h ) )
		&=
		\actionAA'_{\phi(g)}( \psi(h) ) \qquad g \in \mathfrak g,\ h \in \mathfrak h.
	\end{align*}
\end{definition}


\smallskip

\begin{definition}\label{def:Stokes_Chen}
	Let $(\feedbackGG\colon H\to G,\actionGG)$
	be a crossed module of groups.
	Given a functor $\PD\colon \PATHS_2 \to G$ on paths (later this will be a path development),
	a \DEF{$2$-cocycle} is a map $\SD\colon \RECT \to H$ such that:
	\begin{itemize}
		\item It satisfies ``Stokes' '' theorem,
		\begin{align}
			\label{eq:stokes}
			\tag{Stokes}
			\feedbackGG( \SD(\square) ) = \PD(\partial \square), \qquad \square \in \RECT.
		\end{align}

		\item It satisfies a \DEF{horizontal Chen's identity}
			\begin{align*}
				\label{eq:chenHorizontal}
				\SD( \square_A \square_B ) = \actionGG_{\PD(\bottomboundary{A})}( \SD (\square_B) )\, \SD (\square_A),
				\tag{Chen-H}
			\end{align*}
			and it satisfies a \DEF{vertical Chen's identity}
			\begin{align*}
				\label{eq:chenVertical}
				\SD ( \genfrac{}{}{0pt}{1}{\square_B}{\square_A} )
				=
				\SD (\square_A)\, \actionGG_{\PD(\leftboundary{A})}( \SD(\square_B) ). \tag{Chen-V}
			\end{align*}
			(We also just say that $\SD$ satisfies \DEF{Chen's identity}.)

			Here, $\bottomboundary{A}$ is the (equivalence class of the) horizontal path from the bottom-left corner of $\square_A$ to the bottom-right corner of $\square_A$,
			and $\leftboundary{A}$ is the (equivalence class of the) vertical path from the bottom-left corner of $\square_A$ to the top-left corner of $\square_A$.

%
	\end{itemize}
\end{definition}
\begin{remark}
	If $\SD$ and $\PD$
	satisfy
	\eqref{eq:stokes}, then
	\begin{align*}
		\SD(\square_A) \actionGG_{\PD(\NES)}( \SD(\square_B) )
		&\stackrel{\text{\eqref{eq:peiffer}}}=
		\actionGG_{\feedbackGG( \SD(\square_A) )}(\actionGG_{\PD(\NES)}( \SD(\square_B) ))\,
		\SD(\square_A) \\
		&\stackrel{\hphantom{\text{\eqref{eq:peiffer}}}}=
		\actionGG_{\PD(\partial \square_A) \PD(\NES)}( \SD(\square_B) )\,
		\SD(\square_A) \\
		&\stackrel{\hphantom{\text{\eqref{eq:peiffer}}}}=
		\actionGG_{\PD(\bottomboundary{A})}( \SD(\square_B) )\,
		\SD(\square_A).
	\end{align*}
	Here $\NES$ is the boundary path of $\square_A$, \emph{clockwise}
	and leaving out the last linear piece at the bottom.
	This then leads to an alternative right-hand side of \eqref{eq:chenHorizontal}.
	
	Analogously,
	\begin{align*}
		\actionGG_{\PD(\ENW)}( \SD(\square_B) )\,
		\SD(\square_A)
		=
		\SD(\square_A)\, \actionGG_{\PD(\leftboundary{A})}( \SD(\square_B) ).
	\end{align*}
	Here $\ENW$ is the boundary path of $\square_A$, counter-clockwise 
	and leaving out the last linear piece on the left.
\end{remark}

It is well-known in the higher category theory literature, that a $2$-cocycle can be constructed from a compatible pair of forms
(in fact, there is an equivalence of categories, see \cite{schreiber2011smooth}, but note that they work in the globular, i.e.~$2$-category setting whereas our formulation is closer to that of a double category).

\begin{theorem}[{$2$-cocycles from forms / Surface development; \cite{SoncZucc15}}]
	\label{thm:2cocycleIntegration}
	Let $(\feedbackGG\colon H\to G,\actionGG)$ be a crossed module of Lie groups,
	and $(\feedbackAA\colon \mathfrak h \to \mathfrak g, \actionAA)$ be the corresponding crossed module of Lie algebras.

	Let $\alpha \in \Omega^1(\R^2, \mathfrak g)$ and 
	$\beta \in \Omega^2(\R^2, \mathfrak h)$ 
	and assume that they are related in terms of 
	the \DEF{vanishing fake curvature} condition
	\begin{align}
	\label{eq:vanishfakecurvature}
		\feedbackAA( \beta ) = d\alpha + \alpha \wedge \alpha.
	\end{align}

	Then, letting $\PD^\alpha$ be the path development of $\alpha$,
	define the following $1$-form $K^{\alpha,\beta}_{v_1,v_2;s_1,\bullet} \in \Omega^1(\R,\mathfrak h)$,
	\begin{align}
	\label{eq:1formK}
		K^{\alpha,\beta}_{v_1,v_2;s_1,t}\ dt
		&\coloneqq
		\int_{v_1}^{s_1}
			\actionGA_{\PD^\alpha(
				\UPRT{v_1,v_2}{r,t}
				)}
			\left( \beta_{r,t} \right) dr\ dt,
	\end{align}
	and define the \DEF{surface development}\footnote{Recall, see at beginning of \Cref{s:cocycles}, that we write $[a,b]$ for the unique, up to thin homotopy, path in $\R$ from $a$ to $b$.}
	\begin{align}\label{eq:SD_def}
		\SD^{\alpha,\beta}( {}_{v_1,v_2}\square^{s_1,s_2} )
		=
		\PD^{K^{\alpha,\beta}_{v_1,v_2;s_1,\bullet}}([v_2,s_2]).
	\end{align}

	Then: $\SD^{\alpha,\beta}$, together with $\PD^\alpha$, is a $2$-cocycle.

	Moreover, ``switching the axes'', defining
	\begin{align*}
		\bar K^{\alpha,\beta}_{v_1,v_2;t,s_2}\ dt
		&=
		\int_{v_2}^{s_2}
			\actionGA_{
			\PD^\alpha(
				\RTUP{v_1,v_2}{t,q}
				)}
			\left( \beta_{t,q} \right) dq\ dt,
	\end{align*}
	we have
	\begin{align}\label{eq:SD_def2}
		\SD^{\alpha,\beta}( {}_{v_1,v_2}\square^{s_1,s_2} )
		=
		\bar{\PD}^{\bar K^{\alpha,\beta}_{v_1,v_2;\bullet,s_2}}([v_1,s_1]),
	\end{align}
	where $\bar{\PD}$ is as stated in the path development 
	\Cref{thm:1cocycleIntegration}, but
	using the \emph{right}-logarithmic derivative instead.
\end{theorem}

\begin{remark}
	We are mostly interested in the case where $\alpha = X^* a$ and $\beta = X^* b$ are the pullbacks
	of some (constant) 1- respectively 2-forms on $\R^n$ under a smooth enough map $X\colon \R^2 \to \R^n$.
	Alternatively, and closer to \Cref{def:1cocycle}, one could consider
	$2$-cocycles as \emph{double functors} on equivalence classes of such maps
	(``surfaces'' or ``membranes''). In this setting one can also formulate the thin-homotopy invariance of $2$-cocycles.
	We want to avoid (double) categorical language and stick to the concrete formulation in terms of rectangles.
\end{remark}

Note the ostensible breaking of symmetry here: each of the two equivalent definitions of $\SD$ treat the two axes of the plane differently.
In \Cref{sec:sewing}, we present another method to obtain $2$-cocycles. In addition to working in regularity regimes where the above theorem does not apply, this method \emph{is} symmetric in the two axes of the plane.

\begin{theorem}
	\label{thm:2cocycle_sewing}
	Admit the assumptions of \Cref{thm:2cocycleIntegration}.
	Consider the ``germ''
	\begin{align*}
		\widehat{\SD}(\square) = \exp\Big( \int_\square \beta \Big)
	\;, \qquad \square \in \RECT\;.
	\end{align*}
	Then there exists a unique $2$-cocycle $\SD = \SD^{\alpha,\beta}$
	such that\footnote{Recall that \(\lesssim\) means up to a constant factor.}
	\begin{equ}[eq:weak_bound]
	\rho(\SD(\square),\widehat\SD(\square)) \lesssim \bfa^{\lambda}
	\end{equ}
	for any $\lambda \in (2,3]$,
	where $\bfa$ is the length of the longer side of $\square$
	and $\rho$ is a Riemannian metric on $H$.
	Furthermore, $\SD$ satisfies the stronger bound
	\begin{equ}[eq:strong_bound]
	\rho(\SD(\square),\widehat\SD(\square) \lesssim \bfa^{2}\bfb
	\end{equ}
	where $\bfb$ is the other side length of $\square$ with $\bfb\leq \bfa$.
\end{theorem}

The proof of the theorem is given at the end of \Cref{subsec:2D_sewing}.

In case of smooth forms (as is assumed in the above theorems), the two methods coincide:
\begin{corollary}
	Let $\SD$ be the $2$-cocycle from \Cref{thm:2cocycleIntegration}.
	Then $\SD$ agrees with the $2$-cocycle provided by \Cref{thm:2cocycle_sewing}, i.e.~$\SD$ satisfies \eqref{eq:weak_bound}.
\end{corollary}

\begin{proof}
Since $\SD$ is a $2$-cocycle, by the uniqueness part of \Cref{thm:2cocycle_sewing} it suffices to prove that $\SD$ satisfies the bound~\eqref{eq:weak_bound} for some $\lambda>2$.
To this end, observe that $\actionFreeCmGA_{\PD^\alpha(\UPRT{v_1,v_2}{r,t})} = \id + O(\bfa)$.
Therefore $K_{v_1,v_2;s_1,t} = \int_{v_1}^{s_1} \beta_{r,t}\,dr + O(\bfa^2)$
and thus
\begin{equ}
	\SD(\square)
	= \PD^{K_{v_1,v_2;s_1,\bullet}}([v_2,s_2]) 
	= \exp\Big( \int_{v_2}^{s_2}\int_{v_1}^{s_1} \beta_{r,t}\,dr\,dt + O(\bfa^3) \Big) \;.
\end{equ}
Since the double integral above is precisely $\int_\square \beta$, it follows that $\rho(\SD(\square), \widehat\SD(\square))\lesssim \bfa^3$.
\end{proof}


\subsection{Cocycles in the Lie algebra}
\label{ss:cocycles_lie_algebra}

Although the sewing perspective of \Cref{thm:2cocycle_sewing} is far more flexible (and, arguably, more symmetric), we use the ODE perspective of \Cref{thm:2cocycleIntegration} to obtain a Magnus-type expansion.

\newcommand\magnus{{\mathsf{m}}}

We recall the classical setting first, proven for example in \cite[p.263]{blanes1998magnus},\cite[Theorem 4.1]{iserles2000lie}.
As a side-note: in the (graded) nilpotent case of \Cref{Log-MagSeries2}, the Magnus expansion exists without any smallness assumption.

\begin{theorem}
	\label{thm:magnus_expansion_1cocycle}

	There exists $C_\magnus > 0$ such that the following is true.

	Let $a \in \Omega^1(\R^n, \mathfrak g)$ where $\mathfrak g$ is the matrix Lie algebra of a matrix Lie group $G$
	(and we endow $\mathfrak g$ with the norm induced by a sub-multiplicative matrix norm). Let $\PD^a$ be the 1-cocycle of \Cref{thm:1cocycleIntegration}.
	Then, for $\gamma \in \PATHS_n$ satisfying\footnote{$\gamma$ is an equivalence class, so it is enough for one representative to satisfy this.}
	\begin{align*}
		\int_0^1 ||(\gamma^* a)_{t}( \partial_t )|| dt
		=
		\int_0^1 ||a_{\gamma(t)}( \dot \gamma(t) )|| dt < C_{\magnus},
	\end{align*}
	the logarithm 
	\begin{align*}
		\mag \coloneqq \log \PD^a(\gamma)
	\end{align*}
	is well-defined, satisfies
	the \DEF{Magnus ODE} \eqref{MagnusODE}
	and is given
	by a converging \DEF{Magnus series}
	$\mag = \sum_{n=1}^\infty \mag^{(n)}$, where  ($A_t \coloneqq (\gamma^* a)(\partial_t))$
\begin{align}
	\label{Magrec}
\begin{aligned}
	\mag^{(1)}_t &= \int_0^t A_s ds \\
	\mag^{(i)}_t
	&=\int_0^t \sum_{m > 0}\frac{B_m}{m!}\sum_{i_1+\cdots+i_m=i-1}\ad_{\mag^{(i_1)}_s}\cdots \ad_{\mag^{(i_m)}_s}A_s  ds,
	\quad i>1.
\end{aligned}
\end{align}
The first terms are given by:
\allowdisplaybreaks
	\begin{align}
		\label{eq:first_terms_magnus}
		\begin{split}
		\mag^{(2)}_t &= \frac12 \int_0^t [\mag^{(1)}_s, A_s] ds \\
		\mag^{(3)}_t &=
				\frac12 \int_0^t [\mag^{(2)}_s, A_s] ds
					+ \frac1{12} \int_0^t \big[\mag^{(1)}_s, [\mag^{(1)}_s, A_s]\big] ds \\
		\mag^{(4)}_t &=  \frac12 \int_0^t [\mag^{(3)}_s, A_s] ds	
					+ \frac1{12} \int_0^t \big[\mag^{(1)}_s, [\mag^{(2)}_s, A_s]\big] ds 
					+ \frac1{12} \int_0^t \big[\mag^{(2)}_s, [\mag^{(1)}_s, A_s]\big] ds .
		\end{split}
	\end{align}
\end{theorem}

\begin{theorem}
	\label{thm:magnus_expansion_general}
	Let $(\feedbackGG\colon H\to G,\actionGG)$,
	$\alpha \in \Omega^1(\R^2,\mathfrak g)$, $\beta \in \Omega^2(\R^2,\mathfrak h)$ 
		be as in \Cref{thm:2cocycleIntegration}.
	Let $v_1 \le s_1$, $v_2 \le s_2$.
	Then, for
	\begin{enumerate}
		\item either small $|s_1-v_1|$, $|s_2-v_2|$, or
		\item nilpotent Lie algebras $\mathfrak h$ and $\mathfrak g$ (and arbitrary $v_1<s_1$, $v_2<s_2$),
	\end{enumerate}
	for $r \in [v_1, s_1]$, $t \in [v_2, s_2]$,
	\begin{align*}
		\mag^{r,t} \coloneqq \log\left( \PD^{\alpha}( \HOOK{v_1,v_2}{r,t} ) \right),
	\end{align*}
	is well-defined and given by a converging Magnus series.
	$K^{\alpha,\beta}$ can be written as
	\begin{align*}
		K^{\alpha,\beta}_{v_1,v_2;s_1,t}\ dt
		&=
		\int_{v_1}^{s_1}
			\exp(\actionAA_{\mag^{r,t}})
			\left( \beta_{r,t} \right) dr\ dt,
	\end{align*}
	where the exponential is the operator exponential,
	and $\Omega^{\alpha,\beta}$ defined by the Magnus series of $K^{\alpha,\beta}$, corresponding to the Magnus ODE
  \begin{equation}\label{eq:magnus_general}
		\frac{d}{dt} \Omega^{\alpha,\beta}_{v_1,v_2; s_1, t}
		=
		\frac{\ad_{\Omega^{\alpha,\beta}_{v_1,v_2;s_1,t}}}{1 - \exp(-\ad_{\Omega^{\alpha,\beta}_{v_1,v_2;s_1,t}})}( K^{\alpha,\beta}_{v_1,v_2;s_1,t} ), \quad \Omega^{\alpha,\beta}_{v_1,v_2;s_1,v_2} = 0,
	\end{equation}
	is well-defined and
	\begin{align*}
		\SD^{\alpha,\beta}(\SQUARE{v_1,v_2}{s_1,s_2}) = \exp( \Omega^{\alpha,\beta}_{v_1,v_2;s_1,s_2} ),
	\end{align*}
	where $\SD$ is the $2$-cocycle from \Cref{thm:2cocycleIntegration}.

	Moreover, ``switching the axes'' and defining
	for $v_1\le t \le s_1$,
	\begin{align*}
		\bar\mag^{t,q} \coloneqq \log\left( \bar\PD^\alpha( \RTUP{v_1,v_2}{t,q} ) \right),
	\end{align*}
	defining $\bar K^{\alpha,\beta}$
	\begin{align*}
		\bar K^{\alpha,\beta}_{v_1,v_2;s_1,t}\ dt
		&=
		\int_{v_1}^{s_1}
			\exp({\actionAA_{\bar\mag^{t,q}}})
			\left( \beta_{t,q} \right) dq\ dt,
	\end{align*}
	then $\bar \Omega^{\alpha,\beta}$ defined by the Magnus series of $\bar K$, corresponding to the Magnus ODE
	\begin{align*}
		\frac{d}{dt} \bar\Omega^{\alpha,\beta}_{v_1,v_2; t, s_2}
		=
		\frac{\ad_{\bar\Omega^{\alpha,\beta}_{v_1,v_2;t,s_2}}}{\exp(\ad_{\bar\Omega^{\alpha,\beta}_{v_1,v_2;t,s_2}})-1}( \bar K^{\alpha,\beta}_{v_1,v_2;t,s_2} ), \quad \bar\Omega^{\alpha,\beta}_{v_1,v_2;s_1,s_2} = 0,
	\end{align*}
	is well-defined and $\Omega^{\alpha,\beta} = \bar \Omega^{\alpha,\beta}$.
\end{theorem}

\begin{remark} 
The work \cite{ZadMans2022} discusses the Magnus expansion as a method for solving linear systems of PDEs involving two independent parameters that arise in variational problems with Lie group symmetries. The paper shows that, provided certain compatibility conditions are satisfied (think fake curvature condition), the Magnus expansion can be consistently applied sequentially along each parameter to solve these PDEs while ensuring compatibility conditions in Lie group settings are maintained. 
\end{remark}

\begin{remark}
	\label{rem:magnus_expansion_general}
	Take $\alpha \leadsto \eps \alpha, \beta \leadsto \eps^2 \beta$.
	Then the ansatz
	\begin{align*}
		\mag       = \sum_{k=1}^\infty \eps^k \mag^{(k)}, \qquad
		\Omega    = \sum_{k=1}^\infty \eps^k \Omega^{(k)},
	\end{align*}
	leads to $\Omega^{(1)} \equiv 0$ and
	\begin{align*}
		\Omega^{(2)}_{v_1,v_2;s_1,s_2} 
			&= \int_{[v_1,s_1]\times[v_2,s_2]} \beta,
			\qquad
		\Omega^{(3)}_{v_1,v_2;s_1,s_2} = \int_{[v_1,s_1]\times[v_2,s_2]} \actionAA_{\mag^{(1);r,q}}( \beta_{r,q} ) \\
		\Omega^{(4)}_{v_1,v_2;s_1,s_2}
		&=
		\int_{[v_1,s_1]\times[v_2,s_2]} \actionAA_{\mag^{(2);r,q}}(\beta_{r,q})
		+
		\frac12
		\int_{[v_1,s_1]\times[v_2,s_2]} \actionAA_{\mag^{(1);r,q}}( \actionAA_{\mag^{(1);r,q}}(\beta_{r,q}) ) \\
		&\qquad
		+
		\frac12 \int_{v_2 < q^1 < q^2 < s_2}
			\left[ \int_{[v_1,s_1]} \beta_{r,q^1}, \int_{[v_1,s_1]} \beta_{r,q^2} \right].
	\end{align*}
	At order five, we find 
	\begin{align*}
		\lefteqn{\Omega^{(5)}_{v_1,v_2;s_1,s_2}
		=
		\int_{[v_1,s_1]\times[v_2,s_2]} \actionAA_{\mag^{(3);r,q}}(\beta_{r,q})} \\
		&\qquad
		+
		\frac12
		\int_{[v_1,s_1]\times[v_2,s_2]} \actionAA_{\mag^{(1);r,q}}( \actionAA_{\mag^{(2);r,q}}(\beta_{r,q}) )
		+
		\frac12
		\int_{[v_1,s_1]\times[v_2,s_2]} \actionAA_{\mag^{(2);r,q}}( \actionAA_{\mag^{(1);r,q}}(\beta_{r,q}) ) \\
		&\qquad
		+
		\frac16
		\int_{[v_1,s_1]\times[v_2,s_2]} \actionAA_{\mag^{(1);r,q}}( \actionAA_{\mag^{(1);r,q}}( \actionAA_{\mag^{(1);r,q}}(\beta_{r,q}) ) ) \\
		&\qquad
		+
		\frac12 \int_{v_2 < q^1 < q^2 < s_2}
			\left[ \int_{[v_1,s_1]} \actionAA_{\mag^{(1);r,q^1}}( \beta_{r,q^1} ), \int_{[v_1,s_1]} \beta_{r,q^2} \right] \\
		&\qquad
		+
		\frac12 \int_{v_2 < q^1 < q^2 < s_2}
			\left[ \int_{[v_1,s_1]} \beta_{r,q^1}, \int_{[v_1,s_1]} \actionAA_{\mag^{(1);r,q^2}}( \beta_{r,q^2} ) \right].
	\end{align*}
\end{remark}
\begin{proof}
	We argue the non-nilpotent case.
	Embed the finite-dimensional Lie algebra $\mathfrak h$
	into a Lie algebra of matrices, which comes endowed with
	a sub-multiplicative norm $||.||$.
	We argue in this embedding.
	By \Cref{thm:magnus_expansion_1cocycle}
	there is a constant $C_\magnus$
	such that if a curve $\gamma$ satisfies
	\begin{align*}
		\int_0^1 ||\alpha_{\gamma_t}( \dot \gamma_t )|| dt < C_{\magnus},
	\end{align*}
	then $\mag \coloneqq \log \PD^\alpha(\gamma)$
	exists and is given by a converging Magnus series.
	Furthermore
	\begin{align*}
		||\mag|| \le  \psi( \int_{\gamma} ||\alpha|| ),
	\end{align*}
	for a continuous $\psi$ with $\psi(0) = 0$.

	Choose $s_1,s_2$ with $|s_1-v_1|$, $|s_2-v_2|$ small enough such that
	for every $v_1 \le u_1 \le s_1$, $v_2 \le u_2 \le s_2$ the path
	$\gamma^{u_1,u_2} \coloneqq \HOOK{v_1,v_2}{u_1,u_2}$ satisfies
	\begin{align*}
		\int_0^1 ||\alpha_{\gamma^{u_1,u_2}_t}( \dot \gamma^{u_1,u_2}_t )|| dt < C_{\magnus}.
	\end{align*}
	Then, for all such $u_1,u_2$
	\begin{align*}
		\mag^{u_1,u_2} \coloneqq \log\left( \PD^\alpha( \gamma^{u_1,u_2} ) \right)
	\end{align*}
	exists and is given by a converging Magnus series.
	Now, taking the operator norm on the linear endomorphisms of $\mathfrak h$, we get
	\begin{align*}
		||\exp( \actionAA_{\mag^{r,t}} )||
		\le
		\exp( ||\actionAA_{\mag^{r,t}}|| )
		\le
		\exp( C ||\mag^{r,t}|| ),
	\end{align*}
	for some constant $C$.
	Using \Cref{lem:action_exp},
	$K$ from \Cref{thm:2cocycleIntegration} can be written as
	\begin{align*}
		K_{v_1,v_2;s_1,t}\ dt
		&\coloneqq
		\int_{v_1}^{s_1}
			\actionGA_{\PD^\alpha( \UPRT{v_1,v_2}{r,t})}
			(\beta_{r,t}) dr\ dt \\
		&=
		\int_{v_1}^{s_1}
			\exp\left( \actionAA_{\mag^{r,t}} \right)(\beta_{r,t} ) dr\ dt.
	\end{align*}
	By continuity of $\psi$ we can pick $s_1,s_2$ small enough such that
	\begin{align*}
		\int_0^{s_1} dr ||K_{r,t}||
		\le
		\int_0^{s_1} dr \int_0^{s_2} dq \exp( C ||\mag^{r,t}|| ) |\beta_{r,t}|
		\le
		C_{\magnus}.
	\end{align*}
	The first statement then follows.

	Now, decreasing $s_1,s_2$ if necessary, $\bar K$ and $\bar \Omega$ are also well-defined
	and by \Cref{thm:2cocycleIntegration} the claimed identity follows.
\end{proof}

We also encode how Chen's identities translate to the Lie algebra level.

\begin{lemma}
	\label{lem:chen_in_lie_algebra}
	Let $(\feedbackGG\colon H\to G,\actionGG)$,
	$\alpha \in \Omega^1(\R^2,\mathfrak g)$, 
	$\beta \in \Omega^2(\R^2,\mathfrak h)$ 
	be as in \Cref{thm:2cocycleIntegration}
	and $\PD^\alpha$, $\SD$ be the corresponding $1$- resp.~$2$-cocycles.
	Let $v_1 \le s_1$, $v_2 \le s_2$.
	Then, for
	\begin{enumerate}
		\item either small $|s_1-v_1|$, $|s_2-v_2|$, or
		\item nilpotent Lie algebras $\mathfrak h$ and $\mathfrak g$ (and arbitrary $v_1<s_1$, $v_2<s_2$),
	\end{enumerate}
	the following
	\begin{align*}
		\mag(\gamma) 
		\coloneqq \log_G\left( \PD^\alpha( \gamma ) \right), \qquad
		\Omega^{\alpha,\beta}_{v_1,v_2;s_1,s_2}
						\coloneqq \log_H \SD^{\alpha,\beta}( \SQUARE{v_1,v_2}{s_1,s_2} ),
	\end{align*}
	are well-defined, and using temporarily the more suggestive notation
	$\Omega^{\alpha,\beta}(\SQUARE{v_1,v_2}{s_1,s_2}) \coloneqq \Omega^{\alpha,\beta}_{v_1,v_2;s_1,s_2}$,
	\begin{align}
		\label{eq:stokes-LA}
		\tag{Stokes-LA}
		\feedbackAA( \Omega(\square) ) = \mag(\partial \square),
	\end{align}
	and for composable rectangles contained in $[v_1,s_1]\times[v_2,s_2]$ 
	(recall the notation $\bottomboundary{A},\leftboundary{A}$ from \Cref{def:Stokes_Chen}),
	\begin{align}
		\label{eq:chenHorizontal-LA}
		\tag{Chen-H-LA}
		\Omega(\square_A \square_B)
		&=
		\BCH_{\freecm{n}^{-1}}\left( \exp( \actionFreeCmAA_{\mag(\bottomboundary{A})} ) \Omega(\square_B) ,\Omega(\square_A) \right),
	\end{align}
	and
	\begin{align}
		\label{eq:chenVertical-LA}
		\tag{Chen-V-LA}
		\Omega(\genfrac{}{}{0pt}{1}{\square_B}{\square_A})
		&=
		\BCH_{\freecm{n}^{-1}}\left( \Omega(\square_A), \exp( \actionFreeCmAA_{\mag(\leftboundary{A})} ) \Omega(\square_B) \right).
	\end{align}
\end{lemma}
\begin{proof}
	For the logarithms, and then $\BCH$, to be well-defined,
	one needs to show, in the non-nilpotent case,
	that the group elements are close to the identity.
	Calculations very similar to those in the proof of
	\Cref{thm:magnus_expansion_general} show that this is the case
	on small enough rectangles.
\end{proof}


\section{A free crossed module over the free Lie algebra}
\label{s:free_crossed_module}

We introduce a \emph{free} crossed module, which underpins the surface signature, much like the free Lie algebra underpins the path signature
in \Cref{ss:kapranov}.
In \Cref{ss:truncated}, we consider quotients of crossed modules of Lie algebras,
and show how from a nilpotent Lie algebra one can always construct a crossed module of Lie groups.
In \Cref{ss:freeness} we show that the constructed crossed module is indeed free,
and in \Cref{ex:free_nilpotent} we construct the truncated, nilpotent versions of it.


\subsection{Kapranov's construction of a free crossed module of Lie algebras} 
\label{ss:kapranov}

A differential graded Lie algebra
(\textbf{dg Lie algebra})%
\footnote{
A historical reference is Quillen's work on rational homotopy theory \cite{quillen1969rational},
and a modern reference is the book \cite[Section 5.6]{manetti2022lie}.
}
is a $\ZZ$-graded vector space
\begin{align*}
    L = \bigoplus_i L_i,
\end{align*}
with a bilinear, graded map $[.,.]$
(i.e.~$[.,.] \colon L_i \otimes L_j \to L_{i+j}$)
satisfying
\begin{itemize}

    \item 
    (graded anti-symmetry)
    \begin{align*}
        [x,y] = (-1)^{|x||y|+1} [y,x]
    \end{align*}

    \item
    (graded Jacobi identity)
	\begin{align*}
        (-1)^{|x||z|} [x,[y,z]]
        +
        (-1)^{|y||z|} [y,[z,x]]
        +
        (-1)^{|z||x|} [z,[x,y]]
        =
        0
    \end{align*}
    
\end{itemize}
and maps (cohomologically graded, as used by Kapranov \cite{kapranov2015membranes})
\begin{align*}
    d^i\colon L^i \to L^{i+1},
\end{align*}
satisfying
\begin{itemize}
    \item (chain complex)
    \begin{align*}
        d^2 = 0
    \end{align*}
    \item (derivative)
    \begin{align*}
        d [x,y] = [dx,y] + (-1)^{|x|} [x,dy].
    \end{align*}
\end{itemize}

A \DEF{morphism} of dg Lie algebras
\begin{align*}
	\phi\colon L \to L'
\end{align*}
is a chain map (i.e.$\phi_m\colon L_m \to L'_m$ is a family
of linear maps commuting with the differentials:
$d' \phi = \phi d$)
respecting the Lie bracket: $\phi([x,y]) = [\phi(x), \phi(y)], x,y \in L$.


%
%
%
%
%
%


Most results in this section are taken from \cite{kapranov2015membranes} and \cite{reutenauer1990dimensions}.
The main results are \Cref{thm:free_differential_crossed_module}, \Cref{thm:free_generators} and \Cref{thm:kernel_basis}.

Recall that $n \in \N$ is the dimension of the ambient space.
%
Consider the symbols, $p \in \N$,
\begin{align*}
    \Z{I} = \Z{i_1,\ldots,i_p}, \quad I = \{ i_1 < \cdots < i_p \} \subset [n],
\end{align*}
and assign the degrees
$\deg( \Z{I} ) \coloneqq -p + 1$.
Consider the free, graded Lie algebra (see for example \cite[Section B.2]{quillen1969rational}) generated by it, $\f^\bullet(\R^n)$.

\begin{remark}~
	\begin{enumerate}
		\item $\f^0(\R^n)$ can be identified with the free Lie algebra $\FL(\R^n)$.

		\item This is a \emph{graded} object, so that, for example $[ \Z{ij}, \Z{ij} ] \not= 0$	.
	\end{enumerate}
\end{remark}

Define the codifferential $d$, on generators
as
\begin{align*}
    d\Z{I} \coloneqq \frac12 \sum_{J \sqcup K = I} \mathsf{sign}(J,K) [\Z{J}, \Z{K}],
\end{align*}
and extend uniquely using the Leibniz rule.
Here, $\mathsf{sign}$ is the sign of the shuffle permutation corresponding to $J,K$.
For example
\begin{align*}
    d\Z{i} = 0, \quad d\Z{ij} &= [\Z{i}, \Z{j}], \quad d\Z{ijk} = [\Z{i}, \Z{jk}] - [\Z{j}, \Z{ik}] + [\Z{k}, \Z{ij}] \\
	&d[\Z{ij}, \Z{kl}] = [[\Z{i},\Z{j}], \Z{kl}] - [\Z{ij}, [\Z{k}, \Z{l}]].
\end{align*}

The following proposition is a direct consequence of the definition of $d$.
\begin{proposition}
	\label{prop:dgLie}
	$(\f^\bullet(\R^n), d)$ is a dg Lie algebra.
	Moreover,\footnote{In fact it is exact everywhere except at $0$, \cite[Proposition 1.1.4]{kapranov2015membranes}.}
    \begin{align*}
        \im d^{-1} = [\f^0,\f^0].
    \end{align*}
\end{proposition}

The following lemma is implicit in \cite{kapranov2015membranes},
we spell it out for completeness.
\begin{lemma}
	\label{lem:basis_f_minus_1}
	A linear basis for $\f^{-1}(\R^n)$ is given by
	\begin{align*}
		[\Z{i_1}, [\Z{i_2}, \dots, [\Z{i_{p-2}}, \Z{i_{p-1}i_p}]\dots]], \quad p \in \N_{\ge2}, i_j \in [n], j=1,...,p,\ i_{p-1} < i_p.
	\end{align*}
\end{lemma}
\begin{proof}
	Let $\f^0(\R^n)_\ell$ be the subspace of $\f^0(\R^n)$
	generated by brackets with $\ell$ symbols $\Z{i}$.
	Let $\f^{-1}(\R^n)_\ell$ be the subspace of $\f^{-1}(\R^n)$
	generated by brackets with $\ell-1$ symbols $\Z{i}$
	and one symbol $\Z{jk}$.
	Then,
	\begin{align*}
		\dim \f^{-1}(\R^n)_\ell = n^{\ell-1} \cdot \binom{n}{2}.
	\end{align*}
	Indeed, we only have to count the number of Lyndon words
	over the alphabet $\{ \Z{i}, \Z{jk} \}$ of length $\ell$
	that contain exactly one letter from $\{ \Z{jk} \}$.
	For this it is enough to count the number of Lyndon words of length $\ell$
	containing exactly one letter from $\{ \Z{jk} \}$ which is equal to $\Z{12}$,
	and multiply the result by $\binom{n}{2}$.

	Without loss of generalization we can assume that in the order
	on the alphabet, the letter $\Z{12}$ is minimal.
	Then, every such Lyndon word
	is of the form
	\begin{align*}
		\Z{12} \Z{i_1} \dots \Z{i_{\ell-1}},
	\end{align*}
	with arbitrary $i_1, \dots, i_{\ell-1} \in [n]$. This gives the claimed dimension.

	Now it remains to show that an arbitrary element of 
	$\f^{-1}(\R^n)$ can we written in terms of ``right-combed'' brackets,
	with the letters $\Z{ij}$ appearing at the end.
	This is a routine check using the Jacobi identity and anti-symmetry.
	For $\ell=2$ it is immediate.
	Let it be true up to $\ell-1$.
	Then, let $x \in \f^{-1}(\R^n)_\ell$,
	\begin{align*}
		x = [ P, Q ],
	\end{align*}
	with $P \in \f^{-1}(\R^n)_{j}$ and $Q \in \f^{0}(\R^n)_k$ with $j+k = \ell$.
	By induction hypothesis, $P$ can be assumed to be right-combed.
	For the free Lie algebra it is well-known that $Q$ can be assumed right-combed.
	If $k=1$, we can re-arrange $x$ to be right-combed by using anti-symmetry.
	Else write
	\begin{align*}
		P &= [\Z{i_1}, [\Z{i_2}, \dots, [\Z{i_{q-2}}, \Z{i_{q-1}i_q}] \dots ]] \\
		Q &= [\Z{w_1}, [\Z{w_2}, \dots, [\Z{w_{k-2}}, [\Z{w_{k-1}},\Z{w_k}]] \dots ]].
	\end{align*}
	Then, by the Jacobi identity,
	\begin{align*}
		[P,Q]
		&=
		\left[ \left[ P, \Z{w_1}\right], \left[\Z{w_2}, \dots, [\Z{w_{k-2}}, [\Z{w_{k-1}},\Z{w_k}]] \dots \right]\right] \\
		&\qquad
		+
		\left[\Z{w_1}, \left[P, [\Z{w_2}, \dots, [\Z{w_{k-2}}, [\Z{w_{k-1}},\Z{w_k}]] \dots ]\right]\right].
	\end{align*}
	The second term is $[\Z{w_1}, \tilde P]$, where $\tilde P \in \f^{-1}(\R^n)_{\ell-1}$
	which can be right-combed by induction hypothesis. Hence, again using anti-symmetry,
	$[\Z{w_1}, \tilde P]$ can be right-combed.

	Regarding
	\begin{align*}
		\left[ [ P, \Z{w_1}], [\Z{w_2}, \dots, [\Z{w_{k-2}}, [\Z{w_{k-1}},\Z{w_k}]] \dots ]\right],
	\end{align*}
	we can apply the Jacobi identity again and continue in this fashion: one term
	is covered by the induction hypothesis and the other term has moves $P$ one bracketing-level deeper.
\end{proof}

A dg Lie algebra $\h^\bullet$ is \textbf{semiabelian}
if $[\h^{\le -1},\h^{\le -1}] = 0$.
Let $\h$ be any dg Lie algebra. Define its \textbf{semiabelianization} as
\begin{align*}
    \h^\bullet_\sab \coloneqq \h^\bullet / ( [\h^{\le -1},\h^{\le -1}] + d[\h^{\le -1},\h^{\le -1}]).
\end{align*}

The main algebraic object used in this paper is the following.
\begin{theorem}
	\label{thm:free_differential_crossed_module}
	Consider the dg Lie algebra $\f^\bullet(\R^n)$ from the beginning of this subsection,
	perform its semiabelianization
	and restrict to degrees $-1$ and $0$, i.e.
	\begin{align*}
		\freecm{n}^{-1} &\coloneqq \f^{-1}_\sab(\R^n) = \f^{-1}(\R^n) / d [\f^{-1}(\R^n),\f^{-1}(\R^n)] \\
		   \freecm{n}^0 &\coloneqq \f^0_\sab(\R^n) = \f^0(\R^n).
	\end{align*}
	Define the \DEF{derived bracket} on $\freecm{n}^{-1}$ as
		\index{derived backet@\([x, y]_{\freecm{n}^{-1}}\), derived bracket on \(\freecm{n}^{-1}\)}
	\begin{align*}
		[x, y]_{\freecm{n}^{-1}} \coloneqq [d x, y].
	\end{align*}
	The differential $d$ on $\f^\bullet(\R^n)$ induces a map
	$\feedbackFreeCmAA\colon \freecm{n}^{-1} \to \freecm{n}^{0}$.
	The adjoint action of $\f^0(\R^n)$ on $\f^{-1}(\R^n)$
	induces an action $\actionFreeCmAA$ of $\freecm{n}^0$ on $\freecm{n}^{-1}$.
	We use the derived bracket $[.,.]_{\freecm{n}^{-1}}$ on $\freecm{n}^{-1}$.
	Then:
	\begin{enumerate}
		\item 
		\index{freecm@\(\freecm{n}^{\ge -1}\), free differential crossed module over $\R^n$}
		\index{feedback freecm@\(\feedbackFreeCmAA\), feedback of the differential crossed module over $\freecm{n}^{\ge -1}$}
		\index{action freecm@\(\actionFreeCmAA\), action of the differential crossed module over $\freecm{n}^{\ge -1}$}
		$(\feedbackFreeCmAA\colon \freecm{n}^{-1} \to \freecm{n}^0, \actionFreeCmAA)$ is a differential crossed module,
		the free differential crossed module over the free Lie algebra
		$\freecm{n}^0$ and the map $\Z{ij}\mapsto[\Z{i},\Z{j}]$.

		\item
		$\freecm{n}^0$ is the free Lie algebra over $\R^n$
		and, as Lie algebras and as $\GL(\R^n)$-representations
		\begin{align*}
			\freecm{n}^{-1}
			\cong
			\ker \feedbackFreeCmAA \oplus [\freecm{n}^0, \freecm{n}^0].
		\end{align*}

	\end{enumerate}
\end{theorem}
\begin{remark}
	When we endow $\freecm{n}^{-1}$ with the derived bracket,
	$\freecm{n}^{\ge -1}$ is \emph{not} a dg Lie algebra.
	When we write $\freecm{n}^\bullet$ we mean the crossed module (i.e.~using the derived bracket).
	When we write $\f^\bullet_\sab(\R^n)$ we mean the dg Lie algebra.
\end{remark}
\begin{proof}~
	\begin{enumerate}

		\item
			This is not proven in \cite{kapranov2015membranes}
			and we prove it in \Cref{ss:freeness}.

		\item
			By \Cref{prop:dgLie},
			$\feedbackFreeCmAA$ is onto $W = [\freecm{n}^0,\freecm{n}^0]$.
			Since $W$ is a free Lie algebra (\Cref{thm:free_generators})
			there exists a Lie morphism $s\colon W \to \freecm{n}^{-1}$
			such that $\feedbackFreeCmAA \circ s = \id_W$.
			Since $\ker \feedbackFreeCmAA$ is in the center of $\freecm{n}^{-1}$,
			we get, as Lie algebra,
			\begin{align*}
				\freecm{n}^{-1}
				\cong
				\ker \feedbackFreeCmAA \oplus [\freecm{n}^0, \freecm{n}^0].
			\end{align*}

			Now, $\feedbackFreeCmAA$ is $\GL(\R^n)$-equivariant, so that $\ker \feedbackFreeCmAA$ is a $\GL(\R^n)$-representation.
			Moreover, $[\freecm{n}^0, \freecm{n}^0]$ is clearly a $\GL(\R^n)$-representation.

	\end{enumerate}
\end{proof}

Following \cite{kapranov2015membranes}
we now describe,
a set of (free) generators for the Lie algebra $[\freecm{n}^0, \freecm{n}^0]$, \Cref{thm:free_generators},
(this goes back to \cite{reutenauer1990dimensions}),
and
a linear basis for the (commutative) Lie algebra $\ker \feedbackFreeCmAA$, \Cref{thm:kernel_basis},
the proofs of which are deferred to
\Cref{app:kapranov}.

We shall need the sets
\newcommand\I{\mathsf{I}}
\newcommand\J{\mathsf{J}}
\begin{align*}
	\I_{n,p} &\coloneqq \left\{ i_1 \ge i_2 \ge \dots \ge i_{p-1} < i_p \mid i_j \in [n] \right\} \\
	\J_{n,p} &\coloneqq \left\{ i_1 \ge i_2 \ge \dots \ge i_{p-2} < i_{p-1} < i_p \mid i_j \in [n] \right\}.
\end{align*}

\newcommand\cl{{\operatorname{cl}}}
\begin{theorem}
	\label{thm:free_generators}
	The Lie algebra $W\coloneqq[\freecm{n}^0, \freecm{n}^0]$ is free over the generators
	\begin{align*}
		[\Z{i_1}, [\Z{i_2}, \dots, [\Z{i_{p-1}}, \Z{i_p}]\dots]], \quad (i_1, \dots, i_p) \in \I_{n,p},\ p \ge 2.
	\end{align*}

	A basis  of $W$ compatible with these generators is given as follows.
	Consider $H_1$ the following subset of binary trees with leaves indexed by $[n]$:
	\begin{align*}
		( i_1, (i_2, \dots, (i_{p-2}, (i_{p-1}, i_p) ) \dots ) ), \quad (i_1, \dots, i_p) \in \I_{n,p},\ p \ge 2.
	\end{align*}
	(In particular: the Lie bracketings of $H_1$ provides the above mentioned generators.)
	Order $H_1$ totally (in \Cref{ex:dimensions} we use the lexicographic order on the tree's foliation).
	Define the following sets of trees recursively:
	\begin{align*}
		H_{n+1} \coloneqq \{ (t_1, (t_2, \dots, (t_{p-1}, t_p) ) \dots ) \mid t_1, \dots, t_p \in H_n, t_1 \ge t_2 \ge \dots \ge t_{p-1} < t_p \},
	\end{align*}
	and order $H_{n+1}$ totally.
	Then a basis of $W$ is given by the Lie bracketing of the trees in $\bigcup_{n \ge 1} H_n$.

	\medskip

	As a $\GL(\R^n)$-representation, $(W/[W,W])_\ell$ is the irrep corresponding to the shape $(\ell-1,1)$.
	(On $W/[W,W]$ we consider the grading induced from the grading of $W$, which is supported in degrees $\ge 2$.)
\end{theorem}
\begin{proof}
	Every Lie subalgebra of a free Lie algebra is free, \cite[Theorem 2.5]{reutenauer1993free},
	hence $[\freecm{n}^0, \freecm{n}^0]$ is free.
	%

	The statement regarding the free generators and the compatible basis is proven
	in \cite[Section 5.3]{reutenauer1993free}.
	The $\GL(\R^n)$-representation structure follows
	from \cite[Theorem 2]{reutenauer1990dimensions}
	(and is made more explicit with \Cref{thm:reutenauers_map} \Cref{it:rho0iso} together with \Cref{lem:cochain} \Cref{it:gamma_is_irrep}).
\end{proof}

\begin{theorem}
	\label{thm:kernel_basis}
	The abelian Lie algebra $\ker \feedbackFreeCmAA$ is isomorphic,
	as a $\GL(\R^n)$-representation, to a direct
	sum of the representations corresponding to the shapes $(p-1,1,1)$, $p \ge 3$.

	In particular, a basis for it is indexed by $\J_{n,p}, p \ge 3$.
	A concrete indexing is given by
	\begin{align*}
		J_{n,p} \ni (i_1, \dots, i_{p-2}, i_{p-1}, i_p)
		&\mapsto
		[\Z{i_1}, [\Z{i_2}, \dots, [\Z{i_{p-2}}, \Z{i_{p-1}i_p}]\dots]]
		-
		[\Z{i_1}, [\Z{i_2}, \dots, [\Z{i_{p-1}}, \Z{i_{p-2}i_p}]\dots]] \\
		&\qquad\qquad
		+
		[\Z{i_1}, [\Z{i_2}, \dots, [\Z{i_{p}}, \Z{i_{p-2}i_{p-1}}]\dots]] \in \ker \feedbackFreeCmAA.
	\end{align*}
\end{theorem}
\Cref{thm:kernel_basis} follows from
\Cref{lem:rho_minus_1} and goes back to Kapranov.
We provide an explicit indexing of a basis for convenience,
although it is immediate for experts in representation theory.
%

\medskip

We finish with some low-dimensional dimension counting and examples.
Note that (see for example \cite[Theorem 6.3 (1)]{fulton2013representation})
\begin{align*}
	\# \I_{n,p} = \prod_{1 \le i < j \le n} \frac{ \lambda_i - \lambda_j + j - i }{j-i},\qquad
	\# \J_{n,p} = \prod_{1 \le i < j \le n} \frac{ \lambda_i - \lambda_j + j - i }{j-i}.
\end{align*}

\begin{example}
	Consider ambient dimension $n=2$.
	Then $\#\J_{2,p} = 0$ for all $p$. Hence $\ker \feedbackFreeCmAA = 0$
	and
	\begin{align*}
		\freecm{2}^{-1} \cong [\freecm{2}^0,\freecm{2}^0].
	\end{align*}

	\begin{center}
	\begin{tabular}{ c||c|c|c|c|c }
		$p$ (level) & $\dim \FL_p(\R^2) = \dim\mathfrak f^0_p$ & \gray{$\# \I_{2,p}$} & $\# \J_{2,p}$ & $\dim (d [\mathfrak f^{-1},\mathfrak f^{-1}])_p$ & $\dim \mathfrak f^{-1}_p$ \\ 
		\hline
		2 &  1 & \gray{ 1} & 0 & 0 &  1 \\  
		3 &  2 & \gray{ 2} & 0 & 0 &  2 \\
		4 &  3 & \gray{ 3} & 0 & 1 &  4 \\
		5 &  6 & \gray{ 4} & 0 & 2 &  8
	\end{tabular}
	\end{center}


\end{example}

\begin{example}
	\label{ex:dimensions}
	Consider ambient dimension $n=3$.

	\begin{center}
	\begin{tabular}{ c||c|c|c|c|c }
		$p$ (level) & $\dim \FL_p(\R^3) = \dim\mathfrak f^0_p$ & \gray{$\# \I_{3,p}$} & $\# \J_{3,p}$ & $\dim (d [\mathfrak f^{ -1},\mathfrak f^{-1}])_p$ & $\dim \mathfrak f^{-1}_p$ \\ 
		\hline
		2 &  3 & \gray{ 3} & 0 &  0 &  3 \\  
		3 &  8 & \gray{ 8} & 1 &  0 &  9 \\
		4 & 18 & \gray{15} & 3 &  6 & 27 \\
		5 & 48 & \gray{24} & 6 & 27 & 81
	\end{tabular}
	\end{center}

	Here $\dim (d [\mathfrak f^{ -1},\mathfrak f^{-1}])_p$ is the dimension of level $p$ of the image $d [\mathfrak f^{ -1},\mathfrak f^{-1}]$,
	the space that is quotiened by to obtain $\mathfrak f^{-1}_{\mathsf{sab}}(\R^3) = \mathfrak g^{-1}$.
	Their values 
	must be equal to
	\begin{align*}
		\dim \mathfrak f^{-1}_p - \dim\mathfrak f^0_p - \# \J_{3,p},
	\end{align*}
	but at the moment it is unclear to us how to obtain these values from first principles.

\bigskip
A basis for $\FL(\R^3)$ up to level $5$, compatible with the derived series (\Cref{thm:free_generators} with lexicographic order on the $H_k$):
\begin{center}
\begin{minipage}{.8\linewidth}
\begin{tiny}
\begin{verbatim}
------- LEVEL= 1 (# =  3) ------- LEVEL= 4 (# =  18) ------- LEVEL= 5  (# =  48)
1                         (1, (1, (1, 2)))           (1, (1, (1, (1, 2))))    ((1, (1, 2)), (1, 2))
2                         (2, (1, (1, 2)))           (2, (1, (1, (1, 2))))    ((1, (1, 3)), (1, 2))
3                         (2, (2, (1, 2)))           (2, (2, (1, (1, 2))))    ((1, (1, 2)), (1, 3))
------- LEVEL= 2 (# =  3) (3, (1, (1, 2)))           (2, (2, (2, (1, 2))))    ((1, (1, 3)), (1, 3))
(1, 2)                    (3, (2, (1, 2)))           (3, (1, (1, (1, 2))))    ((1, 2), (2, (1, 2)))
(1, 3)                    (3, (3, (1, 2)))           (3, (2, (1, (1, 2))))    ((1, 3), (2, (1, 2)))
(2, 3)                    (1, (1, (1, 3)))           (3, (2, (2, (1, 2))))    ((1, 2), (2, (1, 3)))
------- LEVEL= 3 (# =  8) (2, (1, (1, 3)))           (3, (3, (1, (1, 2))))    ((1, 3), (2, (1, 3)))
(1, (1, 2))               (2, (2, (1, 3)))           (3, (3, (2, (1, 2))))    ((1, 2), (2, (2, 3)))
(2, (1, 2))               (3, (1, (1, 3)))           (3, (3, (3, (1, 2))))    ((1, 3), (2, (2, 3)))
(3, (1, 2))               (3, (2, (1, 3)))           (1, (1, (1, (1, 3))))    ((1, (1, 2)), (2, 3))
(1, (1, 3))               (3, (3, (1, 3)))           (2, (1, (1, (1, 3))))    ((1, (1, 3)), (2, 3))
(2, (1, 3))               (2, (2, (2, 3)))           (2, (2, (1, (1, 3))))    ((2, (1, 2)), (2, 3))
(3, (1, 3))               (3, (2, (2, 3)))           (2, (2, (2, (1, 3))))    ((2, (1, 3)), (2, 3))
(2, (2, 3))               (3, (3, (2, 3)))           (3, (1, (1, (1, 3))))    ((2, (2, 3)), (2, 3))
(3, (2, 3))               ((1, 2), (1, 3))           (3, (2, (1, (1, 3))))    ((1, 2), (3, (1, 2)))
                          ((1, 2), (2, 3))           (3, (2, (2, (1, 3))))    ((1, 3), (3, (1, 2)))
                          ((1, 3), (2, 3))           (3, (3, (1, (1, 3))))    ((2, 3), (3, (1, 2)))
                                                     (3, (3, (2, (1, 3))))    ((1, 2), (3, (1, 3)))
                                                     (3, (3, (3, (1, 3))))    ((1, 3), (3, (1, 3)))
                                                     (2, (2, (2, (2, 3))))    ((2, 3), (3, (1, 3)))
                                                     (3, (2, (2, (2, 3))))    ((1, 2), (3, (2, 3)))
                                                     (3, (3, (2, (2, 3))))    ((1, 3), (3, (2, 3)))
                                                     (3, (3, (3, (2, 3))))    ((2, 3), (3, (2, 3)))
\end{verbatim}
\end{tiny}
\end{minipage}
\end{center}

A corresponding basis for a complement to $\ker d$ inside of $\freecm{3}^{-1}$, the ``non-abelian part'', is given by 
replacing in each expression the last bracket $[i,j]$ by $\Z{ij}$.

A basis for $\Gamma^\cl_2(\R^3)$ up to level $5$, according to \Cref{thm:kernel_basis}
is given as follows:
\begin{align*}
	&[\Z1,\Z{23}] - [\Z2,\Z{13}] + [\Z3,\Z{12}] \\
 	&[\Z1,[\Z1,\Z{23}]] - [\Z1,[\Z2,\Z{13}]] + [\Z1,[\Z3,\Z{12}]], \\
 	&[\Z2,[\Z1,\Z{23}]] - [\Z2,[\Z2,\Z{13}]] + [\Z2,[\Z3,\Z{12}]], \\
 	&[\Z3,[\Z1,\Z{23}]] - [\Z3,[\Z2,\Z{13}]] + [\Z3,[\Z3,\Z{12}]].
\end{align*}

\end{example}


\subsection{Quotients of crossed modules and nilpotent crossed modules}
\label{ss:truncated}

For a Lie algebra $(\mathfrak g, [.,.])$, recall the notion of \DEF{Lie ideal} which is a subset $V \subset \mathfrak g$ such that for all $v \in V$ and all $g \in \mathfrak g$ the Lie bracket $[v,g] \subset V$. It follows that a Lie ideal is a Lie subalgebra.

\begin{definition}
  Let \((\feedbackAA\colon\mathfrak{h}\to\mathfrak{g},\actionAA)\) be a crossed module of Lie algebras.
  A \DEF{crossed ideal} is a pair \((I,J)\) such that
  \begin{enumerate}
    \item \(I\subset\mathfrak{h}\) and \(J\subset\mathfrak{g}\) are Lie ideals,
    \item \(\feedbackAA(I)\subset J\), 
    \item \(\actionAA_g (I)\subset I\) for all \(g \in\mathfrak{g}\),
    \item \(\actionAA_x(\mathfrak{h})\subset I\) for all \(x\in J\).
  \end{enumerate}
\end{definition}

The proofs of the next three statements are in \Cref{ape:proofs}.
\begin{proposition}
  \label{prop:quotient}
  Let \((\feedbackAA\colon\mathfrak{h}\to\mathfrak{g},\actionAA)\) be a crossed module of Lie algebras and \((I, J)\) a crossed ideal.
  There exists a unique crossed module \((\feedbackAA'\colon\mathfrak{h}/I\to\mathfrak{g}/J,\actionAA')\).
  Moreover, the pair \((\psi,\phi)\), where \(\psi\colon\mathfrak{h}\to\mathfrak{h}/I\) and \(\phi\colon\mathfrak{g}\to\mathfrak{g}/J\) are the canonical projections, is a morphism of crossed modules of Lie algebras.
\end{proposition}

We now restrict to the case where both \(\mathfrak{g}\) and \(\mathfrak{h}\) are nilpotent and realize the group laws via the Baker--Campbell--Hausdorff formula.
Let \(\mathfrak{g}\) be a nilpotent Lie algebra and \(G=\exp(\mathfrak{g})\) the corresponding group.
Given \(x,y\in\mathfrak{g}\) there exists a unique element \(x\star y\in\mathfrak{g}\) such that \(\exp(x)\exp(y)=\exp(x\star y)\), that is, \(x\star y\coloneqq \operatorname{BCH}(x,y)\). The nilpotency assumption implies that this operation is well-defined, in the sense that the BCH series terminates and so it yields an element of \(\mathfrak{g}\), no completion is needed. 
Multiplication in \(G\) implies that this operation is associative.
Indeed,
\begin{align*}
  \exp((x\star y)\star z) &= \exp(x\star y)\exp(z)\\
                          &= \exp(x)\exp(y)\exp(z)\\
                          &= \exp(x)\exp(y\star z)\\
                          &= \exp(x\star(y\star z)).
\end{align*}
Inverses for this operation are simply given by \(x^{\star -1}=-x\).
Therefore, \((\mathfrak{g},\operatorname{BCH}_\mathfrak{g})\) is a group, where \(\mathfrak{g}\) is regarded as a set.

\begin{lemma}
  \label{lem:exp}
  Let \(\mathfrak{h}\) be a nilpotent Lie algebra and \(\delta\in \operatorname{Der}(\mathfrak{h})\).
  Then \(\exp(\delta)\in \operatorname{Aut}(\mathfrak{h})\).
  In particular, if \(H=(\mathfrak{h},\mathrm{BCH}_{\mathfrak{h}})\) is the corresponding group, then \(\exp(\delta)\in \operatorname{Aut}(H)\).
\end{lemma}

\begin{proposition}
  \label{prop:crossedmodule}
  Let \((\feedbackAA\colon\mathfrak{h}\to\mathfrak{g},\actionAA)\) be a crossed module of nilpotent Lie algebras.
  Define \(G\coloneqq(\mathfrak{g},\mathrm{BCH}_{\mathfrak{g}})\), \(H=(\mathfrak{h},\mathrm{BCH}_{\mathfrak{h}})\), an action \(\actionGG\colon G\to\operatorname{Aut}(H)\) by
  \[
    \actionGG_g\coloneqq\exp(\actionAA_g)
  \]
  and a feedback \(\feedbackGG\colon H\to G\) by
  \[
    \feedbackGG(h) \coloneqq \mathfrak{t}(h).
  \]
  Then, \((\feedbackGG\colon H\to G,\actionGG)\) is a crossed module of Lie groups.
\end{proposition}

Our main example of this construction is given in \Cref{ex:free_nilpotent}.

\subsection{Freeness}
\label{ss:freeness}

A crossed module satisfying \Cref{def:diffxmod} except for \cref{eq:peifferDiff} is called a \DEF{precrossed module}.
We recall the definition of a free differential crossed module \cite[p. 69]{Martins2016}, \cite[p.455]{whitehead1949combinatorial}, \cite{ratcliffe1980free}, \cite{hog1993two}.

\begin{definition}\label{def:freepxmod}
  A precrossed module
  \((\feedbackAA\colon\mathfrak{e}\to\mathfrak{g},\actionAA)\), together with an
  injection $\iota\colon E \to \mathfrak e$, is a \DEF{free precrossed module} over
  \(\feedbackAA_0\colon E\to\mathfrak{g}\) if given any precrossed module
  \((\feedbackAA'\colon\mathfrak{e}'\to\mathfrak{g}',\actionAA')\) and maps
  \(\phi\colon\mathfrak{g}\to\mathfrak{g}'\), \(\psi_0\colon E\to\mathfrak{h}'\) such
  that \(\phi\circ\feedbackAA_0 = \feedbackAA'\circ \psi_0\), there exists a unique
  \(\psi\colon\mathfrak{e}\to\mathfrak{e}'\) with \(\psi\circ\iota=\psi_0\), such that
  the pair \((\psi,\phi)\) is a morphism of differential precrossed modules.
  In other words, the following diagram commutes:
  \begin{center}
	\begin{tikzcd}
	\mathfrak e \arrow[r, "t"] \arrow[dd, "\psi", dashed] & \mathfrak g \arrow[dd, "\phi", bend left]                    \\
														  & E \arrow[u, "t_0"] \arrow[lu, "\iota"] \arrow[ld, "\psi_0"'] \\
	\mathfrak e' \arrow[r, "t'"]                          & \mathfrak g'                                                
	\end{tikzcd}
\end{center}

A crossed module is a \DEF{free crossed module} if it satisfies above universal property
with $(\mathfrak t\colon \mathfrak h \to \mathfrak g, \actionAA)$ and $(\mathfrak t'\colon \mathfrak h' \to \mathfrak g', \actionAA')$ being crossed modules.
\end{definition}

Once we have constructed the free precrossed module, the free crossed module can be obtained by taking an appropriate quotient.

\begin{definition}[Peiffer commutator]
  \label{def:Peifcomm}
  Let \((\feedbackAA\colon\mathfrak{e}\to\mathfrak{g},\actionAA)\) be a precrossed module.
  The \DEF{Peiffer commutator} is the bilinear map \(\dsbr{-,-}\colon\mathfrak{e}\otimes\mathfrak{e}\to\mathfrak{e}\) given by
  \[
    \dsbr{x,y}\coloneqq\actionAA_{\feedbackAA(x)}(y)-[x,y].
  \]
\end{definition}

\begin{lemma}\label{lem:Peifcomm}
  The subspace 
  \begin{align*}
	\dsbr{\mathfrak{e},\mathfrak{e}} \coloneqq \{ \dsbr{x,y} \mid x,y \in \mathfrak{e} \},
  \end{align*}
  is invariant under the action of $\actionAA_x$ for all $x \in \mathfrak{g}$.
  Moreover, \(\dsbr{\mathfrak{e},\mathfrak{e}}\subset\ker\feedbackAA\).
\end{lemma}

\begin{proof}
  Let \(z\in\mathfrak{g}\) and \(x,y\in\mathfrak{e}\). Since \(\actionAA_z\in \operatorname{Der}(\mathfrak{e})\) we have
  \begin{align*}
    \actionAA_z(\dsbr{x,y}) &= \actionAA_z(\actionAA_{\feedbackAA(x)}(y)-[x,y]) \\
                            &= \actionAA_z\circ\actionAA_{\feedbackAA(x)}(y) - [\actionAA_z(x),y] - [x,\actionAA_z(y)].
  \end{align*}
  Since the action is a Lie morphism we may write
  \[
    \actionAA_z\circ\actionAA_{\feedbackAA(x)}=\actionAA_{[z,\feedbackAA(x)]} + \actionAA_{\feedbackAA(x)}\circ\actionAA_z,
  \]
  so that
  \begin{align*}
    \actionAA_z(\dsbr{x,y}) &= \actionAA_{\feedbackAA(x)}(\actionAA_z(y))-[x,\actionAA_z(y)]+\actionAA_{[z,\feedbackAA(x)]}(y)-[\actionAA_z(x),y]\\
                            &= \dsbr{x,\actionAA_z(y)}  + \actionAA_{[z,\feedbackAA(x)]}(y)- [\actionAA_z(x),y].
  \end{align*}
  Finally, \cref{eq:equiDiff} we see that \(\actionAA_{[z,\feedbackAA(x)]} = \actionAA_{\feedbackAA(\actionAA_z(x))}\), hence
  \[
    \actionAA_z(\dsbr{x,y}) = \dsbr{x,\actionAA_z(y)} + \dsbr{\actionAA_z(x), y}.
  \]

  To verify the inclusion, recall that \(\feedbackAA\) is a Lie morphism and that \cref{eq:equiDiff} holds, so that
  \begin{align*}
    \feedbackAA(\dsbr{x,y}) &=  \feedbackAA(\actionAA_{\feedbackAA(x)}(y))-[\feedbackAA(x),\feedbackAA(y)]\\
                            &=  [\feedbackAA(x),\feedbackAA(y)] - [\feedbackAA(x),\feedbackAA(y)]\\
                            &= 0.\qedhere
  \end{align*}
\end{proof}
Let us consider the Lie ideal \(\mathfrak{p}\subset\mathfrak{e}\) generated by \(\dsbr{\mathfrak{e},\mathfrak{e}}\).
\Cref{lem:Peifcomm} implies that \(\mathfrak{p}\subset\ker\feedbackAA\) and
$\actionAA_x( \mathfrak p ) \subset \mathfrak p$ for all $x \in \mathfrak g$.

\begin{lemma}\label{lem:xmodquot}
  Let \((\feedbackAA\colon\mathfrak{e}\to\mathfrak{g},\actionAA)\) be a differential precrossed module and consider the Lie algebra \(\mathfrak{h}\coloneqq\mathfrak{e}/\mathfrak{p}\).
  The action and feedback descend to unique maps \(\ddotn\colon\mathfrak{g}\to\mathfrak{h}\) and \(\widetilde{\feedbackAA}\colon\mathfrak{h}\to\mathfrak{g}\), such that \((\widetilde{\feedbackAA}\colon\mathfrak{h}\to\mathfrak{g},\ddotn)\) is a differential crossed module.
\end{lemma}
\begin{proof}
  Since \(\mathfrak{p}\subset\ker\feedbackAA\) there exists a unique map \(\widetilde{\feedbackAA}\colon\mathfrak{h}\to\mathfrak{g}\) such that \(\widetilde{\feedbackAA}\circ\pi=\feedbackAA\), where \(\pi\colon\mathfrak{e}\to\mathfrak{h}\) is the canonical projection.
  Moreover, since \(\actionAA_z(\mathfrak{p})\subset\mathfrak{p}\) for every \(z\in\mathfrak{g}\), the maps \(\ddotn_z\colon\mathfrak{h}\to\mathfrak{h}\), \(\ddotn_z(\pi(x))\coloneqq \pi(\actionAA_z(x))\) are well-defined and give derivations on \(\mathfrak{h}\).

  Next, we note that \cref{eq:equiDiff} is preserved. Indeed, given \(z\in\mathfrak{g}\) and \(h=\pi(x)\in\mathfrak{h}\),
  \begin{align*}
  \widetilde{\feedbackAA}(\ddotn_z(h)) &= \widetilde{\feedbackAA}(\pi(\actionAA_z(x)))\\
  &= \feedbackAA(\actionAA_z(x))\\
  &= [z,\feedbackAA(x)]\\
  &= [z,\widetilde{\feedbackAA}(h)].
  \end{align*}

  Hence, \((\widetilde{\feedbackAA}\colon\mathfrak{h}\to\mathfrak{g},\ddotn)\) is a well-defined differential precrossed module.
  Furthermore, \cref{eq:peifferDiff} is satisfied, since if \(h=\pi(x),h'=\pi(x')\in\mathfrak{h}\),
  \begin{align*}
  \ddotn_{\widetilde{\feedbackAA}(h)}(h') &= \pi(\actionAA_{\feedbackAA(x)}(x')) \\
  &= \pi([x,x'])+\pi(\dsbr{x,x'}) \\
  &= [h,h'].\qedhere
  \end{align*}
\end{proof}

In particular, we obtain the following result.
\begin{corollary}\label{crl:xmodfree}
  Let \((\feedbackAA\colon\mathfrak{e}\to\mathfrak{g},\actionAA)\) be a free differential precrossed module.
  Then the crossed module \((\widetilde{\feedbackAA}\colon\mathfrak{h}\to\mathfrak{g},\ddotn)\) given by \Cref{lem:xmodquot} is a free differential crossed module.
\end{corollary}
 \begin{proof}
   Let \((\feedbackAA\colon\mathfrak{h}'\to\mathfrak{g}',\actionAA')\) be any differential crossed module.
   Suppose we are given \(\feedbackAA_0\colon E\to\mathfrak{g}\) and maps \(\phi\colon\mathfrak{g}\to\mathfrak{g}'\), \(\psi_0\colon E\to\mathfrak{h}'\) such that \(\phi\circ\feedbackAA_0=\feedbackAA'\circ \psi_0\).
   By freeness there exists a unique extension \(\psi\colon\mathfrak{e}\to\mathfrak{h}'\) of \(\psi_0\) such that \((\phi,\psi)\) is a morphism of differential precrossed modules.
 
   We note that \(\mathfrak{p}\subset\ker \psi\).
   Indeed, for any \(x,y\in\mathfrak{e}\) we have
   \begin{align*}
     \psi(\dsbr{x,y}) &= \psi(\actionAA_{\feedbackAA(x)}(y)) - [\psi(x), \psi(y)]\\
                   &= \actionAA'_{\phi(\feedbackAA(x))}(\psi(y)) - [\psi(x), \psi(y)]\\
                   &= \actionAA'_{\feedbackAA'(\psi(x))}(\psi(y)) - [\psi(x), \psi(y)].
   \end{align*}
   Since \((\feedbackAA'\colon\mathfrak{h}'\to\mathfrak{g}',\actionAA')\) is a differential crossed module, \cref{eq:peifferDiff} holds so that
   \[
     \psi(\dsbr{x,y}) = [\psi(x), \psi(y)] - [\psi(x), \psi(y)] = 0.
   \]
   Therefore, \(\psi\) descends to a unique Lie morphism \(\widetilde{\psi}\colon\mathfrak{h}\to\mathfrak{h}'\), hence the claim.
 \end{proof}


We now review a realization of the free differential precrossed module over the data \(\mathfrak{g}=\FL(\mathbb{R}^n)\) whose generators we denote by \(\{\Z1,\dotsc,\Z{n}\}\), \(E=\{\Z{ij}:1\le i<j\le n\}\) and \(\feedbackAA_0(\Z{ij}) = [\Z{i}, \Z{j}]\).
By \Cref{crl:xmodfree} it is enough to understand this construction to get a hold of the free differential crossed module.

Consider the real vector space \(\mathcal{E}\) spanned by \(E\), \(V \coloneqq \operatorname{T}(\mathbb{R}^n)\otimes\mathcal{E}\) and \(\mathfrak{e}\coloneqq \FL(V)\).

The Lie algebra \(\mathfrak{g}\) acts on its universal envelope \(\operatorname{T}(\mathbb{R}^n)\) by left multiplication, here denoted as concatenation of words, and
then on $V$ by
 \begin{align*}
 	\actionAA_x(w\otimes \Z{ij})\coloneqq xw\otimes \Z{ij},
 \end{align*}
 extended linearly. It further extends to \(\mathfrak{e}\) as a derivation, 
\cite[Lemma 0.7]{reutenauer1993free}.
 
 \newcommand{\rdd}{\overset{\scalebox{0.5}{$\bullet\bullet$}}{\rho}}
 The adjoint action of \(\mathfrak{g}\) on itself also extends uniquely to an action of \(\operatorname{T}(\mathbb{R}^n)\) on it via
 \[
   \rdd_{\Z{i_1}\dotsm \Z{i_n}}(x) = [\Z{i_1},[\Z{i_2},\dotsc,[\Z{i_n},x]\dotsc]],
 \]
 for \(x\in\mathfrak{g}\).
 
 The feedback map \(\feedbackAA_0\) extends linearly to \(\mathcal{E}\) and then to \(\feedbackAA_0\colon V\to\mathfrak{g}\) by
 \[
   \feedbackAA_0(w\otimes \Z{ij}) = \rdd_w([\Z{i},\Z{j}]).
 \]
 Finally, it extends as a Lie algebra morphism \(\feedbackAA\colon \mathfrak{e}\to \mathfrak{g}\) since \(\mathfrak{e}\) is free as a Lie algebra.
 
 \begin{proposition}
   The maps \((\feedbackAA\colon\mathfrak{e}\to\mathfrak{g},\actionAA)\) defined above form a differential precrossed module.
 \end{proposition}
 \begin{proof}
   Let us verify \Cref{def:diffxmod}. By construction \(\feedbackAA\colon\mathfrak{e}\to\mathfrak{g}\) is a Lie map.
   Given \(x,y\in\mathfrak{g}\), note that
   \begin{align*}
     \actionAA_{[x,y]}(w\otimes \Z{ij}) &= [x,y]w\otimes \Z{ij}\\
     &= xyw\otimes \Z{ij} - yxw\otimes \Z{ij}\\
     &= \actionAA_x(\actionAA_y(w\otimes \Z{ij}))-\actionAA_y(\actionAA_x(w\otimes \Z{ij})).
   \end{align*}
   Therefore, the linear maps \(\actionAA_{[x,y]}\) and \([\actionAA_x,\actionAA_y]_{\operatorname{Der}(\mathfrak{e})}\) coincide on \(V\), hence on \(\mathfrak{e}\). 
 
   For the same reason, it is enough to verify \cref{eq:equiDiff} on \(V\) since for any \(x\in\mathfrak{g}\) the maps
   \[
     \mathfrak{e}\ni a\mapsto \feedbackAA(\actionAA_x(a))\ \text{ and }\ \mathfrak{e}\ni a\to [x,\feedbackAA(a)]
   \]
   are both derivations on \(\mathfrak{e}\).
   In this case we have
   \begin{align*}
     \feedbackAA(\actionAA_x(w\otimes \Z{ij})) &= \feedbackAA(xw\otimes \Z{ij})\\
                                               &= \rdd_{xw}([\Z{i}, \Z{j}])\\
                                               &= [x,\rdd_{w}([\Z{i}, \Z{j}])]\\
                                               &= [x,\feedbackAA(w\otimes \Z{ij})].\qedhere
   \end{align*}
 \end{proof}
 \begin{theorem}\label{thm:freepxmod}
   The precrossed module	\(\mathfrak{X}_n\coloneqq(\feedbackAA\colon\mathfrak{e}\to\mathfrak{g},\actionAA)\) constructed above is the free precrossed module above \(\mathfrak{g}\) and \(\feedbackAA_0(\Z{ij})=[\Z{i},\Z{j}]\) in the sense of \Cref{def:freepxmod}.
 \end{theorem}
\begin{proof}
Let \((\feedbackAA'\colon\mathfrak{e}'\to\mathfrak{g}',\actionAA')\) be an arbitrary precrossed module and consider maps \(\phi\colon\mathfrak{g}\to\mathfrak{g}'\), \(\psi_0\colon E\to\mathfrak{e}'\) such that \(\phi\circ\feedbackAA_0=\feedbackAA'\circ \psi_0\).

Define the extension \(\psi\colon\mathfrak{e}\to\mathfrak{e}'\) by first setting \(\psi(w\otimes \Z{ij})=\actionAA'_{\phi(w)}(\psi_0(\Z{ij}))\) on \(V\), and extending to \(\mathfrak{e}\) as a Lie algebra map.
Here, we are considering the unique extension \(\phi\colon \operatorname{T}(\mathbb{R}^d)\to \mathcal{U}(\mathfrak{g}')\) as an algebra morphism, and \(\actionAA'\colon\mathcal{U}(\mathfrak{g}')\to \operatorname{Der}(\mathfrak{e}')\) is given by
\[
\actionAA'_{x_1\dotsm x_n}(a) = \actionAA'_{x_1}\circ\dotsm\circ\actionAA'_{x_n}(a)
\]
for all \(a\in\mathfrak{e}'\).

Note that the maps \(\phi\circ\feedbackAA\) and \(\feedbackAA'\circ \psi\) are Lie algebra maps on \(\mathfrak{e}\), therefore it suffices to verify they coincide on \(V\).
We have
\begin{align*}
  \phi\circ\feedbackAA(\Z{k_1}\dotsm \Z{k_n}\otimes \Z{ij}) &= \phi(\rdd_{\Z{k_1}\dotsm \Z{k_n}}([\Z{i},\Z{j}]))\\
&= \phi\left( [\Z{k_1},\dotsc[\Z{k_n},[\Z{i},\Z{j}]]\dotsc] \right)\\
&= [\phi(\Z{k_1}),\dotsc[\phi(\Z{k_n}),\feedbackAA'(\psi_0(\Z{ij}))]\dotsc].
\end{align*}
\Cref{eq:equiDiff} implies that \([\phi(\Z{k}),\feedbackAA'(\psi_0(\Z{ij}))]=\feedbackAA'\left(\actionAA'_{\phi(\Z{k})}(\psi_0(\Z{ij}))\right)\) hence
\begin{align*}
\phi\circ\feedbackAA(\Z{k_1}\dotsm \Z{k_n}\otimes \Z{ij}) &= \feedbackAA'\left(\actionAA'_{\phi(\Z{k_1})}\circ\dotsm\circ\actionAA'_{\phi(\Z{k_n})}(\psi_0(\Z{ij}))\right)\\
&= \feedbackAA'\circ\actionAA'_{\phi(\Z{k_1}\dotsm \Z{k_n})}(\psi_0(\Z{ij})) \\
&= \feedbackAA'\circ \psi(\Z{k_1}\dotsm \Z{k_n}\otimes \Z{ij}).
\end{align*}

In order to check the second equality in \Cref{def:diffxmorph}, we note once more that \(\psi\circ\actionAA\) and \(\actionAA'_\phi\circ \psi\) are two Lie morphisms from \(\mathfrak{g}\) to \(\operatorname{Der}(\mathfrak{e})\), hence it suffices to check that they coincide on generators.
We have
\begin{align*}
\psi\circ\actionAA_{\Z{k}}(w\otimes \Z{ij}) &= \psi(\Z{k}w\otimes \Z{ij})\\
&= \actionAA'_{\phi(\Z{k})}\left( \actionAA'_{\phi(w)}(\psi_0(\Z{ij})) \right)  \\
&= \actionAA'_{\phi(\Z{k})}\circ \psi(w\otimes \Z{ij}).\qedhere
\end{align*}
\end{proof}

We now verify that Kapranov's construction, outlined in \Cref{ss:kapranov}, can be understood in these terms.
That is we will first construct an explicit isomorphism between the crossed module obtained from \(\mathfrak{X}_n\) via \Cref{crl:xmodfree} and the crossed module \((\feedbackFreeCmAA\colon\freecm{n}^{-1}\to\freecm{n}^0,\actionFreeCmAA)\) from \Cref{thm:free_differential_crossed_module}

Since \(\freecm{n}^0=\mathfrak{g}=\FL(\R^d)\) we may take $\phi\colon \mathfrak{g} \rightarrow \mathfrak{g}^0$ to be the identity
map. We must find a Lie isomorphism $\psi\colon\mathfrak{e} \rightarrow\freecm{n}^{-1}$ such that the pair $(\phi, \psi)$ is a morphism of (pre)crossed modules. As before, take $\feedbackAA_0(\Z{ij}) = [\Z{i},\Z{j}]$ and
$\psi_0 (\Z{ij}) \coloneqq \Z{ij}$. Note that
\begin{align*}
\phi \circ \feedbackAA_0(\Z{ij}) & = [\Z{i}, \Z{j}]
= d\Z{ij}
= d(\psi_0 (\Z{ij})).
\end{align*}
By the above construction, $\psi_0$ lifts to $\psi : \mathfrak{e} \rightarrow
\freecm{n}^{-1}$, given by
\begin{equation*}
\psi (\Z{i_1} \cdots \Z{i_p}\otimes\Z{i j})  = [\Z{i_1},[\dotsc,[\Z{i_p},\Z{i j}]] .
\end{equation*}
Now we observe that after taking the quotient to obtain the free crossed module from \(\mathfrak{X}_n\), the Lie bracket on \(\mathfrak{e}/\mathfrak{p}\) is simply given by
\begin{align*}
  [u\otimes \Z{ij},v\otimes \Z{kl}] &= \actionAA_{\feedbackAA(u\otimes \Z{ij})}(v\otimes \Z{kl})\\
                                    &= \feedbackAA(u\otimes\Z{ij})v\otimes\Z{kl} \\
                                    &= [\Z{i_1},[\dotsc,[\Z{i_p},[\Z{i},\Z{j}]\dotsc]]v\otimes\Z{kl}
\end{align*}
for \(u=\Z{i_1}\dotsc\Z{i_p}\), which is then mapped to the corresponding element in \(\freecm{n}^{-1}\) by \(\psi\), where as before \(\freecm{n}^{-1}\) is endowed with the derived bracket.
It is easy to check, in view of \Cref{lem:basis_f_minus_1} that $\psi$ is indeed an isomorphism.

\begin{example}
	\label{ex:free_nilpotent}
	We construct now the free nilpotent crossed module by means of \Cref{prop:quotient}.
	For this, grade \(\freecm{n}^0\) and \(\freecm{n}^{-1}\) according to \(|\Z{i}|=1,|\Z{ij}|=2\) and \(|[x,y]|=|x|+|y|\).
	\footnote{Note that this is different than the homological grading in \Cref{ss:kapranov}.}
	We denote the corresponding homogeneous components by \(\freecm{n;k}^0\) and \(\freecm{n;k}^{-1}\), \(k\ge 1\).
	Given a fixed \(N\in\mathbb{N}\), consider the subspaces
	\[
	J_{\ge N+1}\coloneqq\bigoplus_{k=N+1}^\infty \freecm{n;k}^0, \quad I_{\ge N+1}\coloneqq\bigoplus_{k=N+1}^\infty \freecm{n;k}^{-1}.
	\]
	Clearly \(I_{\ge N+1}\) and \(J_{\ge N+1}\) are Lie ideals.
	It is easy to check from the definitions that the feedback is graded,
	in the sense that \(\feedbackAA(\freecm{n;k}^{-1})\subset\freecm{n;k}^0\), and so \(\feedbackAA(I_N)\subset J_N\).
	Then define
	\begin{align*}
		\freecm{n;\le N}^0    &\coloneqq \freecm{n}^0/J_{\ge N+1}, \\
		\freecm{n;\le N}^{-1} &\coloneqq \freecm{n}^{-1}/I_{\ge N+1}.
	\end{align*}
	Note that this truncation is already defined at \cite[p.41]{kapranov2015membranes}.

	There are canonical morphisms of crossed modules of Lie algebras
	$(\freecm{n;\le N}^0, \freecm{n;\le N}^{-1}) \rightarrow (\freecm{n;\le M}^0, \freecm{n;\le M}^{-1})$
	for $M \le N$, and the inverse limit can
	be considered as infinite formal sums of elements of $\freecm{n}^{-1}$, resp. $\freecm{n}^0$,
	where the coefficient on each finite homogeneity stays constant after finite many summands.
	call these spaces $\freecmseriesh$ and $\freecmseriesg$.
	Since the maps $\feedbackFreeCmAA$ and $\actionFreeCmAA$ are graded,
	they extend to $\freecmseriesh$ and $\freecmseriesg$,
	turning the latter into a differential crossed module.
\end{example}


\section{Magnus expansion of the multiplicative surface signature}
\label{sec:logsigsuf}


Let $(\feedbackFreeCmAA\colon \freecm{n}^{-1} \to \freecm{n}^{0}, \actionFreeCmAA)$
be the free differential crossed module of Lie algebras of \Cref{s:free_crossed_module}.
Given a smooth enough surface $X\colon \R^2 \to \R^n$,
define
\begin{align*}
  \alphafree &\coloneqq X^*( \sum_{i=1}^n \Z{i} \diff x^i )
		   = \sum_{i=1}^n \partial_1 X^{(i)}_{t_1,t_2} \Z{i} \ \diff{t}^1 
		   + \sum_{i=1}^n \partial_2 X^{(i)}_{t_1,t_2} \Z{i} \ \diff{t}^2 \\
	\betafree &\coloneqq X^*( \sum_{1\le i<j \le n} \Z{ij} \diff x^i \wedge \diff x^j )
		  = \sum_{i < j} J^{(ij)}_{t_1,t_2} \Z{ij} \ \diff{t}^1 \wedge \diff{t}^2,
\end{align*}
with Jacobian minor 
\begin{align}\label{eq:J_def}
	J^{(ij)}_{t_1,t_2} \coloneqq \partial_1 X^{(i)}_{t_1,t_2} \partial_2 X^{(j)}_{t_1,t_2} - \partial_1 X^{(j)}_{t_1,t_2} \partial_2 X^{(i)}_{t_1,t_2}.
\end{align}
Then -  formally, since we do not consider the group ``above'' the free crossed module -
$\PD^{\alphafree}( \HOOK{v_1,v_2}{r,t} )$
is equal to the path signature (\Cref{sec:intro}) of the path $X_*( \HOOK{v_1,v_2}{r,t} )$.
%
In the spirit of \Cref{thm:magnus_expansion_general}, let%
\footnote
{
	We denote by
	\begin{align*}
		\Area^{ij}(\gamma)
		\coloneqq
		\frac12 \left( \int d\gamma^{(i)} d\gamma^{(j)} - \int d\gamma^{(j)} d\gamma^{(i)} \right),
	\end{align*}
	the signed area. Note the pre-factor $\frac12$, which for example is not used in the nomenclature of \cite{diehl2020areas}. 
}
\begin{align*}
	\mag^{r,t}
	&\coloneqq \log\left( \PD^{\alphafree}( \HOOK{v_1,v_2}{r,t} ) \right) \\
	&=
	\sum_i (X^{(i)}_{r,t} - X^{(i)}_{v_1,v_2}) \Z{i}
	+
	\sum_{i < j} \Area^{ij}(X_*( \HOOK{v_1,v_2}{r,t} ) ) [\Z{i}, \Z{j}]
	+
	\cdots,
\end{align*}
be the log-signature of that path.
Note that $\mag^{r,t}$ is, in general, \emph{not} an element of
$\freecm{n}^{0}$, but rather an element of the projective limit / formal series space (\Cref{ex:free_nilpotent}).
Define
\begin{align*}
	K^{{\alphafree},\betafree}_{v_1,v_2;s_1,t}\ dt
	&\coloneqq
	\int_{v_1}^{s_1}
		e^{\actionFreeCmAA_{\mag^{r,t}}} \left( \betafree_{r,t} \right) dr\ dt.
\end{align*}
Then, $\Omega^{{\alphafree},\betafree}$ defined by the Magnus series of $K$, corresponding to the Magnus ODE
\begin{align*}
	\frac{d}{dt} \Omega^{{\alphafree},\betafree}_{v_1,v_2; s_1, t}
	=
	\frac{\ad_{\Omega^{{\alphafree},\betafree}_{v_1,v_2;s_1,t}}}{1 - \exp(-\ad_{\Omega^{{\alphafree},\betafree}_{v_1,v_2;s_1,t}})}\left( K^{\alphafree,\betafree}_{v_1,v_2;s_1,t} \right), \quad \Omega^{{\alphafree},\betafree}_{v_1,v_2;s_1,v_2} = 0,
\end{align*}
is well-defined as an element of the space of formal series in $\freecm{n}^{-1}$ (\Cref{ex:free_nilpotent}).

Applying \Cref{thm:magnus_expansion_general} to the nilpotent truncations,
we obtain:
\begin{theorem}
  \label{thm:log_signature}
  The maps
	\begin{align*}
		\mag^{\alphafree}\colon    \PATHS_2 \to \freecmseriesg, \quad
		\Omega^{\alphafree,\betafree}\colon \RECT \to \freecmseriesh,
	\end{align*}
	satisfy the following properties:
	\begin{itemize}
		\item $\feedbackFreeCmAA( \Omega^{\alphafree,\betafree}(\square) ) = \mag^\alphafree(\partial \square)$,
		\item
		For composable rectangles, we have
		\begin{align*}
			\Omega^{\alphafree,\betafree}(\square_A \square_B)
			&=
			\BCH_{\freecm{n}^{-1}}( \exp(\actionFreeCmAA_{\mag^\alphafree(\bottomboundary{A})}) \Omega^{\alphafree,\betafree}(\square_B) ,\Omega^{\alphafree,\betafree}(\square_A) ),
		\end{align*}
		and
		\begin{align*}
			\Omega^{\alphafree,\betafree}(\genfrac{}{}{0pt}{1}{\square_B}{\square_A})
			&=
			\BCH_{\freecm{n}^{-1}}( \Omega^{\alphafree,\betafree}(\square_A), \exp(\actionFreeCmAA_{\mag^\alphafree(\leftboundary{A})}) \Omega^{\alphafree,\betafree}(\square_B) ).
		\end{align*}

	\end{itemize}
\end{theorem}

\begin{remark}
	Although this theorem is formulated in the smooth case,	
	we note that using \Cref{thm:RS_ext}, $\Omega$ with the same properties can be 
	defined for surfaces $X$ with less regularity.
	See in particular \Cref{subsec:Young} where we spell out the so-called Young--H\"older regularity regime.
	Outside of this regime, extra data (iterated integrals) need to be defined ``above'' the surface,
	just as for rough paths.
\end{remark}

We now show that surface development, that is, the 2-cocycle \(\SD\) given by the previous theorem is universal, in a sense similar to how the path signature is universal amongst 1-cocycles.
\begin{theorem}
	\label{thm:universality}
	Let $(\feedbackGG\colon H\to G,\actionGG)$
	be a crossed module of Lie groups
	and
	$(\feedbackAA\colon \mathfrak h\to \mathfrak g,\actionAA)$,
	the corresponding crossed module of Lie algebras.

	Let 
	$(\feedbackFreeCmAA\colon \freecm{-1} \to \freecm{0} ,\actionFreeCmAA)$
	be the free crossed module of Lie algebras of \Cref{s:free_crossed_module}.
  In particular, \(\mathfrak{g}_0\) is the free Lie algebra over \(\mathbb{R}^n\).

  Let $\freecmseriesh$ and $\freecmseriesg$ the corresponding infinite-series
	spaces.

	Let
	\begin{align*}
		a = \sum_i A_i \diff x^i \qquad b = \sum_{i<j} B_{ij} \diff x^i \wedge \diff x^j,
	\end{align*}
	be two constant forms on $\R^n$,
	valued in $\mathfrak g$ and $\mathfrak h$ respectively, and let \(\alpha=X^*a,\beta=X^*b\).
	
	Suppose that vanishing of fake curvature \eqref{eq:vanishfakecurvature} and the assumptions of 
	\Cref{thm:magnus_expansion_general} hold, so that
	the Magnus expansion $\Omega^{\alpha,\beta}$ exists.
	Let $\Omega^{\alphafree,\betafree}$ be the Magnus expansion
	in the context of the free case as above.

  There exists a unique morphism of crossed modules \((\Psi,\Phi)\colon(\tau\colon\freecmseriesh\to\freecmseriesg)\to(\mathfrak{t}\colon\mathfrak{h}\to\mathfrak{g})\) such that $\Omega^{\alpha,\beta} = \Psi( \Omega^{\alphafree,\betafree} )$.
\end{theorem}
\begin{proof}
  We begin by constructing a morphism of crossed modules \(\psi\colon\freecm{-1}\to\mathfrak{h}\), \(\phi\colon\freecm{0}\to\mathfrak{g}\). Since \(\freecm{0}\) is the free Lie algebra on \(n\) generators \(\Z1,\dotsc,\Z{n}\) there exists a unique Lie algebra moprhism \(\phi\) such that \(\phi(\Z{i})=A_i\).
  Now, define \(E\coloneqq\{\Z{ij}:i<j\}\) and \(\psi_0\colon E\to\mathfrak{h}\) by \(\psi_0(\Z{ij})=B_{ij}\). Since we have assumed that the fake curvarture vanishes and we are dealing with constant differential forms we have that \(\mathfrak{t}(B_{ij})=[A_i,A_j]\), hence
  \begin{align*}
    \mathfrak{t}(\psi_0(\Z{ij})) &= \mathfrak{t}(B_{ij}) \\
                                 &= [A_i, A_j] \\
                                 &= \phi([\Z{i},\Z{j}]) \\
                                 &= \phi(\tau(\Z{ij})).
  \end{align*}
  Hence, by freeness (see \Cref{def:freepxmod}) there exists a unique Lie algebra morphism \(\psi\colon\freecm{-1}\to\mathfrak{h}\) such that \((\psi,\phi)\) is a morphism of crossed modules.
  Both maps uniquely extend to the corresponding algebraic completion as
  \[
    \Psi(\sum_{n}L_n) = \begin{cases}\sum_n\Psi(L_n)&\text{ if the series converges absolutely,}\\0&\text{else,}\end{cases}
  \]
  and likewise for \(\Phi\).
  Hence we have a morphism of crossed modules \((\Psi,\Phi)\) as in the statement of the Theorem.

  Now it suffices to check that \(\Psi(K^{\alphafree,\betafree}) = K^{\alpha,\beta}\).
  Indeed, once we have shown this identity, it follows that \(\Omega\coloneqq\Psi(\Omega^{\alphafree,\betafree}_{v_1,v_2;s_1,t})\) solves the ODE
  \begin{align*}
    \frac{\mathrm{d}}{\mathrm{d}t}\Omega &= \Psi\left( \frac{\mathrm{d}}{\mathrm{d}t}\Omega^{\alphafree,\betafree}_{v_1,v_2;s_1,t} \right) \\
        &= \Psi\left( \frac{\ad_{\Omega^{\alphafree,\betafree}_{v_1,v_2;s_1,t}}}{1-\exp(-\ad_{\Omega^{\alphafree,\betafree}_{v_1,v_2;s_1,t}})}(K^{\alphafree,\betafree}_{v_1,v_2;s_1,t}) \right)   \\
        &= \frac{\ad_\Omega}{1-\exp(-\ad_\Omega)}(K^{\alpha,\beta}_{v_1,v_2;s_1,t})
  \end{align*}
  which is identical to the one solved by \(\Omega^{\alpha,\beta}\), hence \(\Psi\circ\Omega^{\alphafree,\betafree}=\Omega^{\alpha,\beta}\).

  We first note that by definition \(\Phi\circ\alphafree=\alpha\) and \(\Psi\circ\betafree=\beta\).
  Next, from \Cref{def:diffxmorph} it follows that (denoting by \(\tilde{\omega}^{r,t}=\log\PD^{\alphafree}(\HOOK{v_1,v_2}{r,t})\) the Magnus expansion of the universal 1-cocycle)
  \[
    \Psi\circ\actionFreeCmAA_{\tilde{\omega}^{r,t}} = \actionAA_{\Phi(\tilde{\omega}^{r,t})}\circ\Psi = \actionAA_{\omega^{r,t}}\circ\Psi,
  \]
  where in the last equality we have used that since \(\Phi\) is itself a Lie algebra morphism it clearly holds that \(\Phi(\tilde{\omega}^{r,t})=\omega^{r,t}\eqqcolon\log\PD^\alpha(\HOOK{v_1,v_2}{r,t})\).
  Iterating this identity it is also immediately clear that for any integer \(n\ge 1\) it holds that \(\Psi\circ\actionFreeCmAA_{\tilde{\omega}^{r,t}}^{n}=\actionAA^n_{\omega^{r,t}}\circ\Psi\) and so \(\Psi\circ\exp(\actionFreeCmAA_{\tilde{\omega}^{r,t}})=\exp(\actionAA_{\omega^{r,t}})\circ\Psi\).
  Hence
  \begin{align*}
    \Psi\left(K^{\alphafree,\betafree}_{v_1,v_2;s_1,t}\right)\,\mathrm{d}t &= \int_{v_1}^{s_1}\Psi\left( \exp(\actionFreeCmAA_{\tilde{\omega}^{r,t}})(\betafree_{r,t}) \right)\,\mathrm{d}r\mathrm{d}t \\
    &= \int_{v_1}^{s_1}\exp(\actionAA_{\omega^{r,t}})(\beta_{r,t})\,\mathrm{d}r\mathrm{d}t\\
    &= K^{\alpha,\beta}_{v_1,v_2;s_1,t}.\qedhere
  \end{align*}
\end{proof}


\subsection{The first few terms explicitly}
\label{ssec:explicit}

	\label{ex:general_case}
	Let $X\colon \R^2 \to \R^n$ be a smooth enough surface
	with $X_{0,0} = 0$.
	To streamline notation, we write $K_{s_1,s_2} = K^{{\alphafree},\betafree}_{0,0;s_1,s_2}$,
	$\Omega_{s_1,s_2} = \Omega^{{\alphafree},\betafree}_{0,0;s_1,s_2}$.
	Then
	\begin{align*}
		K_{s_1,t}\ dt
		&=
		\int_0^{s_1} \exp( \actionFreeCmAA_{\mag^{r,t}} )( \betafree_{r,t} ) dr \ dt\\
		&=
		\int_0^{s_1}
		\left(
		\betafree_{r,t}
		+
		[\mag^{r,t}, \betafree_{r,t}]
		+
		\frac12 [\mag^{r,t}, [\mag^{r,t}, \betafree_{r,t}]]
		+
		\dots
		\right) dr \ dt.
	\end{align*}

	Consider $X \leadsto \eps X$
	for a dummy variable $\eps$.
	Then, collecting terms of $\eps^k$ in $\Omega^{(k)}$,
	we get (compare with \Cref{rem:magnus_expansion_general})
	\begin{align*}
		\Omega^{(1)}_{s_1,s_2} &= 0 \\
		\Omega^{(2)}_{s_1,s_2} &= \sum_{i<j} \int_0^{s_2} dq \int_0^{s_1} dr\ J^{(ij)}_{r,q} \Z{ij} \\
		\Omega^{(3)}_{s_1,s_2} &= \sum_i \sum_{j<k} \int_0^{s_2} dq \int_0^{s_1} dr\ X^{(i)}_{r,q} J^{(jk)}_{r,q} [\Z{i}, \Z{jk}] \\
		\Omega^{(4)}_{s_1,s_2} &=
		\frac12 \sum_i \sum_j \sum_{k<\ell} \int_0^{s_2} dq \int_0^{s_1} dr\ X^{(i)}_{r,q} X^{(j)}_{r,q} J^{(k\ell)}_{r,q} [\Z{i}, [\Z{j}, \Z{k\ell}]] \\
		&\quad +
		\sum_{i<j} \sum_{k<\ell}
			\int_0^{s_2} dq \int_0^{s_1} dr\ 
				\Area^{ij}(X_*(\HOOK{0,0}{r,q})) J^{(k\ell)}_{r,q}  [ [\Z{i}, \Z{j}], \Z{k\ell}] \\
		&\quad +
		\frac12 \sum_{i<j} \sum_{k<\ell}
			\int_0^{s_2} dq_2 \int_0^{q_2} dq_1 \int_0^{s_1} dr\ J^{(ij)}_{r,q_1} \int_0^{s_1} dr'\ J^{(k\ell)}_{r',q_2}  [[\Z{i},\Z{j}], \Z{k\ell}].
	\end{align*}


	We now specialize to the case $n=3$.
	\newcommand\XJ[1]{\mathsf{XJ}^{#1}}
	Written in terms of the basis of \Cref{ex:dimensions} we get,
	with $\XJ{i,j,k} \coloneqq \int_0^{s_2} dq \int_0^{s_1} dr\ X^{(i)}_{r,q} J^{(jk)}_{r,q}$,
	and $I \coloneqq \{ (1,1,2), (2,1,2), (1,1,3), (3,1,3), (2,2,3), (3,2,3) \}$.
	\begin{align*}
	\Omega^{(3)}_{s_1,s_2}
		&= 
		\XJ{1,2,3} \left( [\Z{1}, \Z{23}] - [\Z{2}, \Z{13}] + [\Z{3}, \Z{12}] \right)
		+
		\sum_{(i,j,k) \in I }
		\XJ{i,j,k} [\Z{i}, \Z{jk}] \\
		&\qquad
		+
		( \XJ{2,1,3} + \XJ{1,1,3} ) [\Z{2}, \Z{13}]
		+
		( \XJ{3,1,2} - \XJ{1,2,3} ) [\Z{3}, \Z{12}].
	\end{align*}
	Further, $\Omega^{(4)}$ in terms of that basis has coefficients given in \Cref{tab:general_case}.
	\newcommand\XXJ[1]{\mathsf{XXJ}^{#1}}
	\newcommand\AJ[1]{\mathsf{AJ}^{#1}}
	\newcommand\JJ[1]{\mathsf{JJ}^{#1}}
	\label{eq:omega4}
	\begin{table}[ht]
	\centering
	\newcommand{\specialcell}[2][c]{\begin{tabular}[#1]{@{}l@{}}#2\end{tabular}}
	\resizebox{0.85\textwidth}{!}{%
	\begin{tabular}{ll}
	 Lie basis & Coefficient \\
	 \hline
$[\Z{1}, [\Z{1}, \Z{12}]]$ & $\frac{1}{2} \XXJ{1,1,1,2}$  \\ 
$[\Z{2}, [\Z{1}, \Z{12}]]$ & $ \XXJ{2,1,1,2}$  \\ 
$[\Z{2}, [\Z{2}, \Z{12}]]$ & $\frac{1}{2} \XXJ{2,2,1,2}$  \\ 
$[\Z{3}, [\Z{1}, \Z{12}]]$ & $ \XXJ{1,3,1,2}-\frac{1}{2} \XXJ{1,1,2,3}$  \\ 
$[\Z{3}, [\Z{2}, \Z{12}]]$ & $ \XXJ{2,3,1,2}- \XXJ{1,2,2,3}$  \\ 
$[\Z{3}, [\Z{3}, \Z{12}]]$ & $\frac{1}{2} \XXJ{3,3,1,2}- \XXJ{1,3,2,3}$  \\ 
$[\Z{1}, [\Z{1}, \Z{13}]]$ & $\frac{1}{2} \XXJ{1,1,1,3}$  \\ 
$[\Z{2}, [\Z{1}, \Z{13}]]$ & $ \XXJ{1,2,1,3}+\frac{1}{2} \XXJ{1,1,2,3}$  \\ 
$[\Z{2}, [\Z{2}, \Z{13}]]$ & $\frac{1}{2} \XXJ{2,2,1,3}+ \XXJ{1,2,2,3}$  \\ 
$[\Z{3}, [\Z{1}, \Z{13}]]$ & $ \XXJ{3,1,1,3}$  \\ 
$[\Z{3}, [\Z{2}, \Z{13}]]$ & $ \XXJ{2,3,1,3}+ \XXJ{1,3,2,3}$  \\ 
$[\Z{3}, [\Z{3}, \Z{13}]]$ & $\frac{1}{2} \XXJ{3,3,1,3}$  \\ 
$[\Z{2}, [\Z{2}, \Z{23}]]$ & $\frac{1}{2} \XXJ{2,2,2,3}$  \\ 
$[\Z{3}, [\Z{2}, \Z{23}]]$ & $ \XXJ{3,2,2,3}$  \\ 
$[\Z{3}, [\Z{3}, \Z{23}]]$ & $\frac{1}{2} \XXJ{3,3,2,3}$  \\ 
$[[\Z{1}, \Z{2}], \Z{13}]$ & $\frac{1}{2} \XXJ{1,2,1,3}-\frac{1}{2} \XXJ{1,3,1,2}+ \XXJ{1,1,2,3}+ \AJ{1,2,1,3}- \AJ{1,3,1,2}+\frac{1}{2} \JJ{1,2,1,3}-\frac{1}{2} \JJ{1,3,1,2}$  \\ 
$[[\Z{1}, \Z{2}], \Z{23}]$ & $\frac{3}{2} \XXJ{1,2,2,3}-\frac{1}{2} \XXJ{2,3,1,2}+ \AJ{1,2,2,3}- \AJ{2,3,1,2}+\frac{1}{2} \JJ{1,2,2,3}-\frac{1}{2} \JJ{2,3,1,2}$  \\ 
$[[\Z{1}, \Z{3}], \Z{23}]$ & $\frac{1}{2} \XXJ{1,3,2,3}-\frac{1}{2} \XXJ{2,3,1,3}+ \AJ{1,3,2,3}- \AJ{2,3,1,3}+\frac{1}{2} \JJ{1,3,2,3}-\frac{1}{2} \JJ{2,3,1,3}$  \\ 
\hline
$[\Z1,[\Z1,\Z{23}]] - [\Z1,[\Z2,\Z{13}]] + [\Z1,[\Z3,\Z{12}]]$ & $\frac12 \XXJ{1,1,2,3}$ \\
$[\Z2,[\Z1,\Z{23}]] - [\Z2,[\Z2,\Z{13}]] + [\Z2,[\Z3,\Z{12}]]$ & $\XXJ{1,2,2,3}$ \\
$[\Z3,[\Z1,\Z{23}]] - [\Z3,[\Z2,\Z{13}]] + [\Z3,[\Z3,\Z{12}]]$ & $\XXJ{1,3,2,3}$ \\
	\hline
	\end{tabular}}
	\caption{$\Omega^{(4)}$'s coefficients in terms of the basis of \Cref{ex:dimensions}.
	The last three terms are in the kernel of the feedback.}
	\label{tab:general_case}
	\end{table}
	
	Here we use the short notation
	\begin{align*}
		\XXJ{i,j,k,\ell} &\coloneqq \int_0^{s_2} dq \int_0^{s_1} dr\ X^{(i)}_{r,q} X^{(j)}_{r,q} J^{(k\ell)}_{r,q} \\
		\AJ{i,j,k,\ell}  &\coloneqq \int_0^{s_2} dq \int_0^{s_1} dr\ \Area^{ij}(X_*(\HOOK{0,0}{r,q})) J^{(k\ell)}_{r,q} \\
		\JJ{i,j,k,\ell}  &\coloneqq \int_0^{s_2} dq_2 \int_0^{q_2} dq_1 \int_0^{s_1} dr\ J^{(ij)}_{r,q_1} \int_0^{s_1} dr'\ J^{(k\ell)}_{r',q_2}.
	\end{align*}

	We now proceed with explicit calculations and consistency checks.
	
	\textbf{Linear data.}
	Consider a linear map $X\colon \R^2\to\R^3$, for coefficients $M_{ij} \in \R$,
	\begin{align*}
		X_{r,q} =
		\begin{pmatrix}
			M_{11} r + M_{12} q \\
			M_{21} r + M_{22} q \\
			M_{31} r + M_{32} q
		\end{pmatrix}.
	\end{align*}
	Then, the 1- and 2-forms are given by
	\begin{align*}
		{\alphafree}
			&= (M_{11} \Z{1} + M_{21} \Z{2} + M_{31} \Z3)\ dt_1 + (M_{12} \Z{1} + M_{22} \Z{2} + M_{32} \Z3)\ dt_2 \\
		\betafree
			&=
			\left(
				\det(M^{(12)}) \Z{12}
				+
				\det(M^{(13)}) \Z{13}
				+
				\det(M^{(23)}) \Z{23} \right)
			\ dt_1 \wedge dt_2,
	\end{align*}
	with $M^{(ij)}$ the $2\times 2$ matrix obtained by restricting to rows $i,j$, i.e.~$J^{(ij)}_{t_1,t_2} = \det(M^{(ij)})$ is constant.
	Further
	\begin{align*}
		X^{(i)}_{r,q} = M_{i1} r + M_{i2} q, \qquad
		\Area^{ij}(X_*(\HOOK{0,0}{r,q})) = - \frac{rq}2 \det(M^{(ij)}),
	\end{align*}
	and
	\begin{align*}
		\int_0^{s_2} dq \int_0^{s_1} dr\ J^{(ij)}_{r,q}
		&=
		s_2 s_1 \det(M^{(ij)}) \\[0.5cm]
		\int_0^{s_2} dq \int_0^{s_1} dr\ X^{(i)}_{r,q} J^{(jk)}_{r,q}
		&=
		\int_0^{s_2} dq \int_0^{s_1} dr\ \left( M_{i1} r + M_{i2} q \right) \det(M^{(jk)}) \\
		&
		= ( \frac12 s_2 s_1^2 M_{i1} + \frac12 s_1 s_2^2 M_{i2} ) \det(M^{(jk)}) 
	\end{align*}
	\begin{align*}
		\lefteqn{\int_0^{s_2} dq \int_0^{s_1} dr\ X^{(i)}_{r,q} X^{(j)}_{r,q} J^{(k\ell)}_{r,q}}\\
		&=
		\int_0^{s_2} dq \int_0^{s_1} dr\ 
		(M_{i1} r + M_{i2} q)(M_{j1} r + M_{j2} q) \det(M^{(k\ell)}) \\
		&=
		\left(
			\frac13 s_2 s_1^3 M_{i1} M_{j1} 
			+
			\frac14 s_2^2 s_1^2 M_{i1} M_{j2}
			+
			\frac14 s_2^2 s_1^2 M_{i2} M_{j1}
			+
			\frac13 s_2^3 s_1 M_{i2} M_{j2}
		\right)
		\det(M^{(k\ell)}) 
	\end{align*}
	as well as
	\begin{align*}	
		\int_0^{s_2} dq \int_0^{s_1} dr\ 
			\Area^{ij}(X_*(\HOOK{0,0}{r,q})) J^{(k\ell)}_{r,q}
		&= -\int_0^{s_2} dq \int_0^{s_1} dr\ \frac{rq}2 \det(M^{(ij)}) \det(M^{(k\ell)}) \\
		&= -\frac18 s_2^2 s_1^2 \det(M^{(ij)}) \det(M^{(k\ell)}),
	\end{align*}
	and
	\begin{align*}
		\lefteqn{\int_0^{s_2} dq_2 \int_0^{q_2} dq_1 \int_0^{s_1} dr\ J^{(ij)}_{r,q_1} \int_0^{s_1} dr'\ J^{(k\ell)}_{r',q_2}}\\
		&=
		\int_0^{s_2} dq_2 \int_0^{q_2} dq_1 \int_0^{s_1} dr\ \det(M^{(ij)}) \det(M^{(k\ell)}) \\
		&=
		\frac12 s_2^2 s_1^2 \det(M^{(ij)}) \det(M^{(k\ell)}).
	\end{align*}



	\bigskip

	\textbf{Verifying Stokes (modulo feedback)}

	Looking at the coefficient of $\Z{ij}$, we get,
	using (classical) Stokes' formula
	\begin{align}
		\label{eq:level2_stokes}
		\begin{split}
		\int_0^{s_2} dq \int_0^{s_1} dr\ J^{(ij)}_{r,q} 
		&=
		\int_{[0,s_1]\times [0,s_2]} d[ X^* (\frac12 x_i \diff{x}^j - \frac12 x_j \diff{x}^i) ] \\
		&=
		\int_{ \partial [0,s_1]\times [0,s_2]} X^* (\frac12 x_i \diff{x}^j - \frac12 x_j \diff{x}^i), 
		\end{split}
	\end{align}
	where we recall $J^{(ij)}$ from \eqref{eq:J_def}.
	The final expression is
	the area integral of the path around the boundary, as expected.

	Now, regarding $\Omega^{(3)}$,
	\begin{align*}
		d[ X^*[ x_i x_j \diff{x}^k ] ]
		&=
		X^*[ d[ x_i x_j \diff{x}^k ] ] \\
		&=
		X^*[ x_j \diff{x}^i \wedge \diff{x}^k + x_i \diff{x}^j \wedge \diff{x}^k ] \\
		&=
		X^{(j)} \left( \partial_1 X^{(i)} \partial_2 X^{(k)}
					- \partial_2 X^{(i)} \partial_1 X^{(k)} \right) \diff{t}^1 \wedge \diff{t}^2 \\
					&\qquad
		+
		X^{(i)} \left( \partial_1 X^{(j)} \partial_2 X^{(k)}
					 - \partial_2 X^{(j)} \partial_1 X^{(k)} \right) \diff{t}^1 \wedge \diff{t}^2.
	\end{align*}
	For $i=j=1, k=3$, we get
	\begin{align*}
		\int_0^{s_2} dq \int_0^{s_1} dr\ X^{(1)}_{r,q} J^{(13)}_{r,q}
		&=
		\int_{[0,s_1]\times [0,s_2]} \frac12 d[ X^*( x_1^2 \diff{x}^3 ) ] \\
		&=
		\langle \IIS( X_*( \partial\square ) ), \Z1 \Z1 \Z3 \rangle \\
		&=
		\langle \IIS( X_*( \partial\square ) ), \frac16 (\Z1 \Z1 \Z3 - 2 \Z1\Z3\Z1 + \Z3\Z1\Z1) \rangle,
	\end{align*}
	where we used that the increment is zero, hence
	$\Z3\Z1\Z1 = -\frac12 \Z1\Z3\Z1$ (modulo testing against $\IIS$),
	and
	$\Z1\Z3\Z1 = - 2 \Z1\Z1\Z3$ (modulo testing against the iterated-integrals signature $\IIS$).
	Doing this analogously for all terms, but the first one,
	we see that we get exactly the terms in \Cref{tab:in_basis} on level~$3$.

	\bigskip


	\textbf{Verify horizontal Chen's identity for the first few terms.}

	Consider $\square_A = \SQUARE{0,0}{s_1,s_2},
	\square_B = \SQUARE{s_1,0}{t_1,s_2}$.
	We use the notation $\mag_\gamma \coloneqq \log \PD^{\alphafree}( \gamma )$.
	Then, using \Cref{lem:chen_in_lie_algebra}, \eqref{eq:chenHorizontal-LA},
	\begin{align*}
		&\BCH\left(
				\exp( \actionFreeCmAA_{ \mag(\RT{0,0}{s_1,0}) } ) \Omega_{\square_B},
				\Omega_{\square_A} \right)_2 \\
		&\quad=
		\BCH\left(
				\Omega_{\square_B},
				\Omega_{\square_A} \right)_2
		=
		\Omega^{(2)}_{\square_B} + \Omega^{(2)}_{\square_A} \\
		&\quad=
		\sum_{i<j} \int_0^{s_2} dq \int_{s_1}^{t_1} dr\ J^{(ij)}_{r,q} \Z{ij} 
		+
		\sum_{i<j} \int_0^{s_2} dq \int_0^{s_1} dr\ J^{(ij)}_{r,q} \Z{ij}  \\
		&\quad=
		\sum_{i<j} \int_0^{s_2} dq \int_0^{t_1} dr\ J^{(ij)}_{r,q} \Z{ij}
		=
		\Omega^{(2)}_{\square_A \square_B},
	\end{align*}
	as expected.
	
	Further
	\begin{align*}
		&\BCH\left(
				\exp( \actionFreeCmAA_{ \mag(\RT{0,0}{s_1,0}) } ) \Omega_{\square_B},
				\Omega_{\square_A} \right)_3 \\
		&\quad=
		\BCH\left(
				( \id + \actionFreeCmAA_{\mag(\RT{0,0}{s_1,0})} ) \Omega_{\square_B},
				\Omega_{\square_A} \right)_3
		=
		\Omega^{(3)}_{\square_B}
		+
		\actionFreeCmAA_{\mag(\RT{0,0}{s_1,0})}( \Omega^{(2)}_{\square_B} )
		+
		\Omega^{(3)}_{\square_A}\\
		&\quad=
		\sum_i \sum_{j<k} \int_0^{s_2} dq \int_{s_1}^{t_1} dr\ \left( X^{(i)}_{r,q} - X^{(i)}_{s_1,0} \right) J^{(jk)}_{r,q} [\Z{i}, \Z{jk}] \\
		&\qquad\quad
		+
		\sum_i \left( X^{(i)}_{s_1,0} - X^{(i)}_{0,0} \right)
			\sum_{j<k} \int_0^{s_2} dq \int_{s_1}^{t_1} dr\ J^{(jk)}_{r,q}
		\left[ \Z{i}, \Z{jk} \right] \\
		&\qquad\qquad
		+
		\sum_i \sum_{j<k} \int_0^{s_2} dq \int_0^{s_1} dr\ \left( X^{(i)}_{r,q} - X^{(i)}_{0,0} \right) J^{(jk)}_{r,q} [\Z{i}, \Z{jk}] \\
		&\quad=
		\sum_i \sum_{j<k} \int_0^{s_2} dq \int_0^{t_1} dr\ \left( X^{(i)}_{r,q} - X^{(i)}_{0,0} \right) J^{(jk)}_{r,q} [\Z{i}, \Z{jk}]
		=
		\Omega^{(3)}_{\square_A \square_B},
	\end{align*}
	as expected.

\section{Non-commutative sewing}
\label{sec:sewing}

In this section, we present our aforementioned non-commutative 2D sewing lemma, which is one of the main results of this article.
The general idea is that, starting from local descriptions (or germs), one can construct a global surface development that is close to the local description in a quantitative way.
This idea permeates rough path theory \cite{Lyons98, Gubinelli_04_controlling, Feyel_06_sewing}, and rough analysis in general \cite{Hairer_14,FH20,Harang_21_sewing}, and our contribution in particular generalises previous non-commutative sewing lemmas~\cite{feyel2008non,Bailleul15,gerasimovivcs2021non} to the 2D setting of crossed modules.
Our proof is novel on several points, using a special dyadic decomposition that gives an easy verification of Chen's identity and which appears new even in the 1D setting.

\subsection{1D sewing}
\label{ss:1d_sewing}

For comparison with the more involved 2D case,
we first give a statement of the 1D non-commutative sewing lemma.
Suppose $G$ is a group equipped with a complete metric $\rho$ such that there exists $Q>0$ for which $\rho(xy,xz)\leq Q\rho(y,z)$ and $\rho(yx,zx)\leq Q\rho(y,z)$ for all $x,y,z$ in a ball $\{a\,:\, \rho(a,e_G)<\delta\}$ with $\delta>0$ and where $e_G$ is the identity element of $G$.
(For instance, take $G$ to be a Lie group equipped with a Riemannian metric.)

\begin{lemma}\cite [Theorem 3.4]{gerasimovivcs2021non}\label{lem:1D_sewing}
	Assume that $\widehat\PD \colon [0,T] \times [0,T] \to G$ is a continuous function satisfying
	\begin{align*}
		\rho(\widehat\PD (s,t),\widehat\PD (s,u) \widehat\PD (u,t)) \lesssim |t-s|^\theta\;, \quad \theta > 1\;,
	\end{align*}
and $\widehat{\PD} (t,t)=e_G$.
	Then there exists a unique continuous map $(s,t)\mapsto \PD (s,t)$ satisfying $\PD (s,t) = \PD (s,u) \PD (u,t)$ and $\rho(\PD (s,t), \widehat \PD (s,t)) \lesssim |t-s|^\theta$.
\end{lemma}

\begin{remark}
The statement of \Cref{lem:1D_sewing} is less quantitative than \cite[Thm.~3.4]{gerasimovivcs2021non} (see also \cite[Thm.~6]{feyel2008non}).
Our statement furthermore differs from the latter because we do not assume the ``moderate growth'' condition of \cite[Def.~3.3]{gerasimovivcs2021non}.
The reason for this is that, by a simple inductive argument, the continuity assumption on $\widehat\PD$ in fact implies the moderate growth condition on sufficiently small intervals.
\end{remark}

\subsection{2D sewing}
\label{subsec:2D_sewing}

We now come to the main result of this section, which is a 2D analogue of \Cref{lem:1D_sewing} in the context of crossed modules.
Throughout this section, suppose $H \stackrel{\feedbackGG}\to G \stackrel{\actionGG}\to \Aut(H)$ is a crossed module.
We denote by $e_G,e$ the identities of $G,H$ respectively.
We suppose that $G,H$ are equipped with complete metrics $\rho_G$ and $\rho$
that satisfy the following conditions:
\begin{itemize}
	\item For all $\eps>0$ there exists $\delta>0$ such that for all $x_1,x_2,y,z \in B_\delta(e)\coloneqq\{a\in H\,:\, \rho(a,e)<\delta\}$
\begin{equ}[eq:rho_assump]
	\rho(x_1 y x_2,x_1 z x_2 )\leq \rho(y,z) (1+\eps)
\end{equ}
and for all $g_1,g_2,u,v \in B^G_\delta(e)\coloneqq\{a\in G\,:\, \rho_G(a,e_G)<\delta\}$
\begin{equ}[eq:rho_assump_G]
	\rho_G(g_1 u g_2,g_1 v g_2 )\leq \rho_G(u,v) (1+\eps).
\end{equ}

\item For all $\eps>0$ there exists $\delta>0$
such that, whenever $\rho_G(g,e_G) \vee \rho(h_1,e) \vee \rho(h_2,e) < \delta$,
\begin{equation}
\label{eq:G_act_lip}
	\rho(\actionGG_{g} h_1, \actionGG_{g} h_2) \leq \rho(h_1,h_2)(1+\eps)
\end{equation}
and, whenever $\rho_G(g,e_G) \vee \rho_G(g_1,e_G) \vee \rho_G(g_2,e_G) < \delta$,
\begin{equation}
\label{eq:G_Ad_lip}
	\rho_G(\Ad_{g} g_1, \Ad_{g} g_2) \leq \rho(g_1,g_2)(1+\eps).
\end{equation}

\item The maps $\feedbackGG\colon H\to G$ and $G\times H \ni (g,h) \mapsto \actionGG_g h\in H$ as well as the group operations $(h,h')\mapsto hh'$ and $h\mapsto h^{-1}$ (for both $G$ and $H$) are locally uniformly continuous.
\end{itemize}

We leave it as an exercise to show that these conditions are satisfied if $G,H$ are Lie groups equipped with Riemannian metrics.

\begin{convention}\label{conv:rect}
Throughout this section, all rectangles in $\RECT$ are assumed to be subrectangles of $[0,1]^2$.
Similarly, all paths in $\Paths \coloneqq \Paths_2$ take values in $[0,1]^2$.
\end{convention}

Unless otherwise stated, we let $\square$ denote a rectangle with height $\bfh$ and width $\bfw$.
We denote $\bfa = \max\{\bfh,\bfw\}$ and $\bfb = \min\{\bfh,\bfw\}$.
Define the \DEF{eccentricity} of $\square$ by 
\begin{equ}
E(\square) = \bfa/\bfb\;.
\end{equ}

\begin{definition}\label{def:subcontrol}
A two-element partition $\{\square_A,\square_B\}$ of a rectangle $\square\in \RECT$ is called \DEF{balanced} if $\square_A,\square_B\in\RECT$ are rectangles and
\begin{equ}
\max\{E(\square_A),E(\square_B) \}\leq \max\{3,2E(\square)/3\}\;.
\end{equ}

A function  $F\colon \RECT \to [0,\infty)$ is called a \DEF{subcontrol} if there exists $0<L<1$ such that, for any balanced partition $\{\square_A,\square_B\}$ of a rectangle $\square$,
\begin{equ}[eq:super_add]
F(\square_A) + F(\square_B) \leq L F(\square)
\end{equ}
and
\begin{equ}[eq:vanish_area]
\lim_{|\square| \to 0} F(\square) = 0\;.
\end{equ}
\end{definition}
\begin{remark}
		A (1D) \emph{control} (see e.g. \cite{Lyons98,StFlour07,FV10_RPs}) is a continuous functon
		$\omega\colon \R^2 \to [0,\infty)$, zero on the diagonal, and 
		superadditive, i.e., $\omega(s,t) + \omega(t,u) \leq \omega(s,u)$.
		Our notion of subcontrol plays the role of (a special case of) $\omega^\theta$, where $\omega$ is a 1D control and $\theta>1$,
		which frequently appears in sewing lemmas.
		Indeed, $(s,t)\mapsto |t-s|$ is a 1D control and, for any $\theta>1$, the function $(s,t)\mapsto |t-s|^\theta$
		satisfies $|u-s|^\theta + |t-u|^\theta < L|t-s|^\theta$ for some $L\in (0,1)$ whenever $s<u<t$ and $|s-u|/|t-s|$ and $|t-u|/|t-s|$ are bounded away from zero.
		In the 1D case, the existence of such $L$, while elementary, is a crucial ingredient in the proofs of the sewing lemmas of \cite{Feyel_06_sewing,feyel2008non,Bailleul15,FH20}.
\end{remark}

The following proposition gives a basic example of a subcontrol that we will later use.
\begin{proposition}\label{prop:subcontrol}
Define $F\colon\RECT \to [0,1]$ by
\begin{equ}
F(\square)=
\begin{cases}
	(\bfa\bfb)^\theta \quad &\text{if } E(\square)\leq \bfe \\
	\bfa^\lambda \bfb^\zeta  \quad &\text{if } E(\square) > \bfe
\end{cases}
\end{equ}
where $\bfe\geq 3$, $\theta,\lambda>1, \zeta\geq 0$ with $\lambda\geq \zeta$ and $2\theta \geq \lambda+\zeta$.
Then $F$ satsifies~\eqref{eq:super_add}.
If furthermore $\zeta>0$, then $F$ satisfies~\eqref{eq:vanish_area}, so $F$ is a subcontrol.
\end{proposition}

For the proof, we require the following lemmas.

\begin{lemma}\label{lem:long}
Suppose $E(\square)\geq 9/2$ and that $\bfw > \bfh$. Then any balanced partition $\square_A,\square_B$ is formed by a vertical cut and moreover, if $\bfw_A,\bfw_B,\bfh_A,\bfh_B$ are the widths and heights of $\square_A,\square_B$ respectively,
then $\bfw_i/\bfw \in [\frac13,\frac23]$ and $\bfw_i>\bfh_i$ for $i=A,B$.
\end{lemma}

\begin{proof}
Since $2E(\square)/3 \geq 3$, we have $\bfw_i \leq \bfh 2E(\square)/3 = \frac{2}{3}\bfw$.
Since $\bfw_A+\bfw_B = \bfw$, it follows that $\bfw_A/\bfw \in [\frac13,\frac23]$.
In particular, $\bfw_i \geq \frac13 \bfw = \frac13 E(\square) \bfh \geq \frac32 \bfh$.
\end{proof}

\begin{lemma}\label{lem:areas}
If $\square_A,\square_B$ is a balanced partition of $\square$, then $|\square_i| > \frac{2}{27}|\square|$.
\end{lemma}

\begin{proof}
If $E(\square)\geq 9/2$, it follows from \Cref{lem:long} that $|\square_i|\geq \frac13|\square|$.
On the other hand, if $E(\square) < 9/2$, then suppose without loss of generality that $\square_A,\square_B$ is formed by a vertical cut and that $\bfw_A,\bfw_B,\bfw$ are the widths of $\square_A,\square_B,\square$,
so that $\bfw_A+\bfw_B=\bfw$, and all three rectangles have common height $\bfh$.
Then $\frac \bfw \bfh < 9/2$ and
$\frac{\bfw_i}{\bfh} \geq 1/3$, which implies
$\bfw_i \geq \bfh/3 > \frac{2\bfw}{27}$
and thus $|\square_i| = \bfh_i \bfw > \frac{2\bfw \bfh}{27} = \frac{2}{27}|\square|$.
\end{proof}

\begin{proof}[Proof of \Cref{prop:subcontrol}]
To show~\eqref{eq:super_add},
note that, if $E(\square_A),E(\square_B)>\bfe$, then
$E(\square)>3\bfe/2 \geq 9/2$ and therefore, by \Cref{lem:long}, letting $\bfa_A,\bfa_B$ denote the long sides of $\square_A,\square_B$, we have $\bfa_A,\bfa_B \in \bfa[\frac13,\frac23]$ and $\bfa_A+\bfa_B=\bfa$.
Hence, since $\lambda>1$,
\begin{equ}
F(\square_A)
+
F(\square_B)
= \bfa_A^\lambda \bfb^\zeta
+ \bfa_B^\lambda \bfb^\zeta
\leq L \bfa^\lambda \bfb^\zeta
= F(\square)\;.
\end{equ}
Here and below, we let $L$ denote a constant depending only on $\lambda,\zeta,\theta$ such that $L<1$.
If $E(\square_A)>\bfe\geq E(\square_B)$ then again $E(\square)>9/2$ and thus
\begin{equ}
F(\square_A)
+
F(\square_B)
= \bfa_A^\lambda \bfb^\zeta
+ \bfa_B^\theta \bfb^\theta
\leq
\bfa_A^\lambda \bfb^\zeta 
+ \bfa_B^\lambda \bfb^\zeta
\leq L \bfa^\lambda \bfb^\zeta
= L F(\square)\;,
\end{equ}
where we used that $\bfa_B > \bfb$, $\lambda \geq \zeta$, and $2\theta \geq \lambda+\zeta$.
Finally, suppose $E(\square_A), E(\square_B) \leq \bfe$.
Then, because $\theta>1$ and $\square_A,\square_B$ partition $\square$ and $|\square_i|>\frac{2}{27}|\square|$ by \Cref{lem:areas},
\begin{equ}
F(\square_A)
+
F(\square_B)
= |\square_A|^\theta
+ |\square_B|^\theta
\leq
L|\square|^\theta \leq L F(\square)\;.
\end{equ}
Finally, if $\zeta>0$, then \eqref{eq:vanish_area} is clear. 
\end{proof}

\begin{definition}\label{def:dyadic}
A rectangle $\square$ is called \DEF{dyadic} if its four corners are in $2^{-N}\ZZ$ for some $N\geq 0$.
If $N$ is the smallest integer with this property, then we say that $\square$ is \DEF{$2^{-N}$-dyadic}.
Let $\RECT^N$ denote the set of all $2^{-N}$-dyadic rectangles.
If $\square\in\RECT^N$ has both height and width $2^{-N}$ (in particular $\square$ is a square), then we say that $\square$ is \DEF{elementary}.
We write $\DD^N\subset\RECT^N$ for the set of all elementary $2^{-N}$-dyadic squares.
\end{definition}

\begin{convention}\label{conv:metrics}
We equip $\Paths$ with the bounded variation metric $|\gamma-\eta|_{BV} = \int |\dot\gamma-\dot\eta|$
and $\Rect$ with the Hausdorff metric where we identify every $\square\in\RECT$ canonically with a closed subset of $\R^2$.
By \Cref{conv:rect}, note that $\Rect$ is a compact metric space.
\end{convention}

The following is our announced 2D non-commutative ``sewing lemma''.
For a function $f$, we frequently use the shorthand $f x = f(x)$.

\begin{theorem}[Sewing Lemma]\label{thm:2D_sewing}
	Suppose we are given locally uniformly continuous
	\begin{equs}[eq:PD_SD]
		\PD\colon \PATHS &\to G \\
		\widehat\SD \colon \RECT &\to H,
	\end{equs}
	such that $\widehat\SD (\square) = e$ if $|\square|=0$.
	Assume that $\PD$ is a $1$-cocycle, i.e.~$\PD$ satisfies the usual Chen's identity
	\begin{equ}[eq:PD_Chen]
	\PD( \gamma \sqcup \eta ) = \PD(\gamma ) \PD(\eta) \;.
	\end{equ}

	
	Assume there exists a subcontrol $F$ such that, for all $2^{-N}$-dyadic rectangles $\square\in\RECT^N$, $N\geq 0$, $\widehat\SD$ satisfies the following $2$-parameter approximate Chen identity:
	\begin{align}\label{eq:approx_Chen1}
		\rho\Big(\widehat\SD (\square_A\square_B) 
		,
		\actionGG_{\PD(\bottomboundary{1})}( \widehat\SD \square_B )\, \widehat\SD \square_A \Big)
		\leq F (\square)\;,
	\end{align}
	for all balanced horizontal partitions $\square_A\square_B=\square$ and
	\begin{align}\label{eq:approx_Chen2}
		\rho\Big(\widehat\SD( \genfrac{}{}{0pt}{1}{\square_B}{\square_A} ),
		\widehat\SD \square_A\, \actionGG_{\PD(\leftboundary{1})}( \widehat\SD \square_B )\Big)
		\leq F(\square)
	\end{align}
	for all balanced vertical partitions $\genfrac{}{}{0pt}{1}{\square_B}{\square_A}=\square$.
	
	Finally, assume there exists a subcontrol $\bar F$ such that $\widehat\SD$ satisfies approximate Stokes
	\begin{align}\label{eq:approx_Stokes}
		\rho_G(\feedbackGG(\widehat\SD  \square ) , \PD(\partial \square)) \leq \bar F(\square)	
	\end{align}
	for all elementary squares $\square\in\DD^N$.

	Then there exists a unique map
	\begin{align*}
		\SD \colon \RECT  \to H,
	\end{align*}
	which is a $2$-cocycle, i.e.~satisfies Chen's identities
	\eqref{eq:chenHorizontal}-\eqref{eq:chenVertical} and
	Stokes' identity \eqref{eq:stokes},
	and, for some $C>0$, for all $\square\in\RECT$
	\begin{align}\label{eq:approx_bound}
		\rho(\SD \square , \widehat\SD \square ) \leq C F(\square)\;.
	\end{align}
	Furthermore, $\SD$ is continuous and the map
	\begin{equ}[eq:sewing_map]
	\CC(\Paths,G)\times \CC(\Rect,H) \ni (\PD,\widehat\SD)  \mapsto \SD \in \CC(\Rect,H)
	\end{equ}
	is uniformly continuous on any set $\BB_1\times\BB_2$ equipped with the uniform metric, where $\BB_1\subset \CC(\Paths,G)$ is such that
	\begin{equ}[eq:BB1]
	\lim_{\eps\downarrow0} \sup_{|\gamma|_{BV}<\eps}\sup_{\PD \in \BB_1} \rho_G(\PD(\gamma),e_G) = 0
	\end{equ}
	and $\BB_2\subset \CC(\Rect,H)$ is such that 
	\begin{equ}[eq:BB2]
	\lim_{\eps\downarrow0} \sup_{|\square|<\eps}\sup_{\widehat\SD\in\BB_2} \rho(\widehat\SD \square,e) = 0
	\end{equ}
	and there exists a subcontrol $F$ such that all $\widehat\SD \in \BB_2$ satisfy the bounds~\eqref{eq:approx_Chen1}-\eqref{eq:approx_Chen2},
	and, for each $\widehat\SD\in\BB_2$, there exists $\bar F$ such that \eqref{eq:approx_Stokes} is satisfied.
\end{theorem}

\begin{remark}
The subcontrol $\bar F$ in \Cref{thm:2D_sewing} is used only to ensure that the constructed $\SD$ satisfies Chen and Stokes.
In particular, it does not feature in the bound~\eqref{eq:approx_bound}. 
\end{remark}

\begin{remark}\label{rem:rect_inc}
For another multiparameter sewing lemma in the case when $H$ is a linear space and $\PD, G$ play no role, see \cite[Lem.~14]{Harang_21_sewing}.
Restricting to this setting, our sewing lemma is different from (in fact weaker than) the latter.
The reason behind this difference is that we do not employ the notion of rectangular increments
and it remains open whether there is an analogue of \Cref{thm:2D_sewing} that utilizes rectangular increments.
A similar problem was noticed in \cite{lee2024_2Dsigs}.
\end{remark}

\begin{proof}[Proof of Uniqueness in \Cref{thm:2D_sewing}]
Suppose $\bar \SD$ is another map satisfying the given properties. In particular $\rho(\SD \square,\bar \SD \square) \lesssim F(\square)$. By Chen's identities for $\SD,\bar\SD$, it suffices to show that $\bar\SD (\square) = \SD$ for all $\square\in\RECT$ with sufficiently small diameter and $E(\square)\leq 2$.

Let us divide $\square$ in half $k$ times horizontally and vertically to arrive at rectangles $\square_i$, $1\leq i\leq 2^k$, with $E(\square_i)\leq 2$.
Chen's identity for $\SD,\bar\SD$ and the bounds~\eqref{eq:rho_assump} and~\eqref{eq:G_act_lip}
imply the telescoping bound
\begin{equ}
	\rho( \SD \square,\bar \SD \square) \lesssim \sum_{i=1}^{2^k} F(\square_i)\;,
\end{equ}
where the bound is uniform in $k$.
Since each division resulted in a balanced partition and since $F$ is a subcontrol,
we have $\sum_{i=1}^{2^k} F(\square_i) \leq L^k F(\square)$.
Taking $k\to\infty$, it follows that $\rho(\SD \square,\bar \SD \square)=0$.
\end{proof}

\begin{remark}\label{rem:weaker_F}
	The above proof of uniqueness in \Cref{thm:2D_sewing} reveals that $\SD$ is in fact the unique $2$-cocycle such that, for some $C>0$,
	\begin{equ}
	\rho(\SD {\square} , \widehat\SD {\square}) \leq C \tilde F(\square)
	\end{equ}
	where $\tilde F$ is any function satisfying~\eqref{eq:super_add} for some $L<1$ and $F\lesssim \tilde F$.
	In other words, for uniqueness, it is not necessary to assume the vanishing area condition~\eqref{eq:vanish_area}.
\end{remark}

For the proof of existence of \Cref{thm:2D_sewing}, we require the following definitions and lemma.
For simplicity, we write $\gamma\bar\gamma$ for the concatenation of two paths $\gamma,\bar\gamma$.

\begin{definition}\label{def:lassos}
For $N\geq 0$, let $\PATHS^N$ denote the set of paths $c \in \PATHS$ such that $c$ linearly interpolates points of $2^{-N}\ZZ^2$.

A \DEF{$2^{-N}$-dyadic lasso} is a loop of the form $\ell = c b c^{-1}$ where $b,c\in\PATHS^N$ and
$b$ is the boundary of an elementary square (recall \Cref{def:dyadic} and the boundary map \eqref{eq:boundary_rect_paths}).
We let $\Lasso^N$ denote the set of all $2^{-N}$-dyadic lassos.
\end{definition}

Our notion of a lasso is inspired by the works \cite{Gross85,Driver89,Levy03} on 2D gauge theory and is 
similar to that of a `kite' of Yekutieli \cite{yekutieli2016nonabelian}.
Given $\PD\colon \Paths \to G$ and $\SD\colon\Rect\to H$,
we extend the definition of $\SD$ to all lassos $\ell=c b c^{-1}$ by
\begin{equ}
\SD(\ell) = \actionGG_{\PD(c)}\SD(b)
\end{equ}
where $\SD(b) \coloneqq \SD(B)$ and $B$ is the elementary square such that $\partial B= b$.

\begin{lemma}\label{lem:Stokes_to_Chen}
Suppose $\PD\colon \Paths \to G$ and $\SD\colon\Rect\to H$ where $\PD$ satisfies the usual Chen's identity \eqref{eq:PD_Chen}
and Stokes identity for every $\square\in\DD^N$:
\begin{equ}[eq:Stokes]
\feedbackGG( \SD  \square ) = \PD(\partial \square)\;.
\end{equ}
Suppose $\ell_i = c_i b_i c_i^{-1}\in\Lasso^N$, $i=1,\ldots, k$, and $\bar \ell_i = \bar c_i \bar b_i {\bar c_i}^{-1}\in\Lasso^N$, $i=1,\ldots, \bar k$, are lassos
where $c_i,\bar c_i \in \PATHS^N$ all share a common basepoint.
Suppose that
\begin{equ}[eq:ell_bar_ell]
\prod_{i=1}^k \ell_i \sim_t \prod_{i=1}^{\bar k} \bar\ell_i\;.
\end{equ}
where $\sim_t$ denotes thin-homotopy equivalence.
Then
\[
	\prod_{i=1}^k \SD(\ell_i) = \prod_{i=1}^{\bar k} \SD(\bar \ell_i) \;.
\]
\end{lemma}

\begin{remark}
	In the language of double categories,
	which we avoid in this paper,
	the lemma states the following.
	Given a partition of a (big) rectangle into
	(small) rectangles,
	and an assignment
	of edges to $1$-cells and of (small) rectangles
	to $2$-cells both in some double category
	which is compatible with the boundary operations
	(this is where Stokes' identity comes in),
	then there is a unique ``aggregated''
	$2$-cell assigned to the big rectangle
	(\cite[Theorem 1.2]{dawson1993general}
		and
	\cite[Section 3]{fiore2008model}).
\end{remark}

\begin{proof}
We have the equalities
\begin{align*}
	\feedbackGG\Big( \prod_{i=1}^k \SD(\ell_i) \Big)
	&=
	\prod_{i=1}^k \feedbackGG(\SD(\ell_i))
	=
	\prod_{i=1}^k \feedbackGG(\actionGG_{\PD(c_i)}\SD (b_i))
	=
	\prod_{i=1}^k \Ad_{\PD(c_i)}( \feedbackGG(\SD (b_i))) \\
	&=
	\prod_{i=1}^k \Ad_{\PD(c_i)}( \PD (\partial b_i))
	=
	\prod_{i=1}^k \PD( \ell_i ),
\end{align*}
where we used Stokes identity in the 4th equality.
It follows from~\eqref{eq:ell_bar_ell} and Chen's identity for $\PD$ that
\begin{equ}[eq:T_to_P]
\feedbackGG\Big( \prod_{i=1}^k \SD(\ell_i) \Big)
=
\PD\Big( \prod_{i=1}^k \ell_i \Big)
=
\PD \Big( \prod_{i=1}^k \bar \ell_i \Big)
=
\feedbackGG\Big( \prod_{i=1}^{\bar k} \SD(\bar \ell_i) \Big)
\end{equ}
for all $\ell_i,\bar\ell_i$ as in the lemma statement.

In the rest of the proof, 
it suffices to consider the case when $\prod_{i=1}^k \ell_i$ is thin-homotopy equivalent to a constant path and show that $\prod_{i=1}^k \SD(\ell_i)=e$.

We treat $2^{-N}\ZZ^2$ as graph in which two vertices $x,y\in 2^{-N}\ZZ^2$ are connected if and only if $|x-y|=2^{-N}$ for the Euclidean norm $|\cdot|$.
Recalling that all rectangles are assumed to be subsets of $[0,1]^2$,
we may identify the group of loops modulo thin-homotopy with the fundamental group $\pi_1 (2^{-N}\ZZ^2\cap [0,1]^2)$.

Observe that $\pi_1 (2^{-N}\ZZ^2\cap [0,1]^2)$ is free and a set of generators is given by lassos $\{\hat\ell_i\}_{i=1}^{4^N}\subset \Lasso^N$ where the corresponding $\hat b_i$ ranges over all the squares in $2^{-N}\ZZ^2\cap [0,1]^2$.
(The corresponding $\hat c_i\in\PATHS^N$, which connects $0$ to the bottom-left corner of $\hat b_i$, is left unspecified.)

For every $\ell_i$, there exists $\hat\ell_{k(i)}$ such that $b_i = \hat b_{k(i)}$.
Let us write $\ell_i \sim_t c_i \hat c_i^{-1} \hat\ell_{k(i)} \hat c_i c_i^{-1}$.
Then, for the loop $c_i \hat c_i^{-1}$, there exists a word $w_i$ on the alphabet $\{1,\ldots, 4^N\}$ such that
\begin{equ}
c_i \hat c_i^{-1} \sim_t \prod_{w \in w_i} \hat\ell_{w}\;.
\end{equ}
It follows from Chen's identity \eqref{eq:PD_Chen}, the equalities \eqref{eq:T_to_P}, and identity \eqref{eq:peiffer} that
\begin{equ}
\SD(\ell_i) = \actionGG_{\PD(c_i)} \SD(b_i) = \actionGG_{\PD(c_i \hat c_i^{-1})} \SD(\hat \ell_{k(i)})
=
\actionGG_{\feedbackGG \prod_w \SD(\hat\ell_w)} \SD(\hat \ell_{k(i)})
=
\Ad_{\prod_w \SD(\hat\ell_w)} \SD(\hat \ell_{k(i)}) \;.
\end{equ}
Doing this for every $\ell_i$, it follows that $\prod_{i=1}^k \SD(\ell_i)$ can be written as $\prod_{w\in w^*} \SD(\hat\ell_w)$ for a word $w^*$.
Furthermore, $\prod_{w\in w^*} \hat\ell_w$ is thin-homotopy equivalent to $\prod_{i=1}^k\ell_i$ and thus to the constant path.
By freeness of $\pi_1 (2^{-N}\ZZ^2\cap [0,1]^2)$, it follows that the word $w^*$ is trivial and therefore $\prod_{w\in w^*} \SD(\hat\ell_w)=e$ as claimed.
\end{proof}

\begin{definition}\label{def:midway}
Suppose $\square\in\Rect^N$ is a $2^{-N}$-dyadic rectangle with lower-left corner $x\in 2^{-N}\ZZ^2$, height $\bfh=\bfH 2^{-N}$ and width $\bfw=\bfW 2^{-N}$, $\bfH,\bfW\geq 1$.
We define the \DEF{midway partition} $\{\square_A,\square_B\}$ of $\square$ as follows.

If $\square$ is an elementary square, i.e.~$\bfH=\bfW = 1$,
then $\square_A$ is the $2^{-N-1}$-dyadic rectangle with lower-left corner $x$, height $\bfh$ and width $\bfw/2$.

If $\square$ is not an elementary square then
\begin{enumerate}
	\item 
if $\bfW\geq \bfH$ and $\bfW$ is even, then $\square_A$ is the $2^{-N}$-dyadic rectangle with lower-left corner $x$, height $\bfh$ and width $\bfw/2$,

\item 
if $\bfW\geq \bfH$ and $\bfW$ is odd, then $\square_A$ is the $2^{-N}$-dyadic rectangle with lower-left corner $x$, height $\bfh$ and width $2^{-N}(\bfW+1)/2$,

\item 
if $\bfH>\bfW$ and $\bfH$ is even, then $\square_A$ is the $2^{-N}$-dyadic rectangle with lower-left corner $x$, height $\bfh/2$ and width $\bfw$,

\item 
if $\bfH>\bfW$ and $\bfH$ is odd, then $\square_A$ is the $2^{-N}$-dyadic rectangle with lower-left corner $x$, height $2^{-N}(\bfH+1)/2$ and width $\bfw$.
\end{enumerate}
Then $\square_B$ is the complement of $\square_A$ in $\square$.
See \Cref{fig:2D_midway} for examples.

\begin{figure}[ht]
	\centering
	\includegraphics[width=1\textwidth]{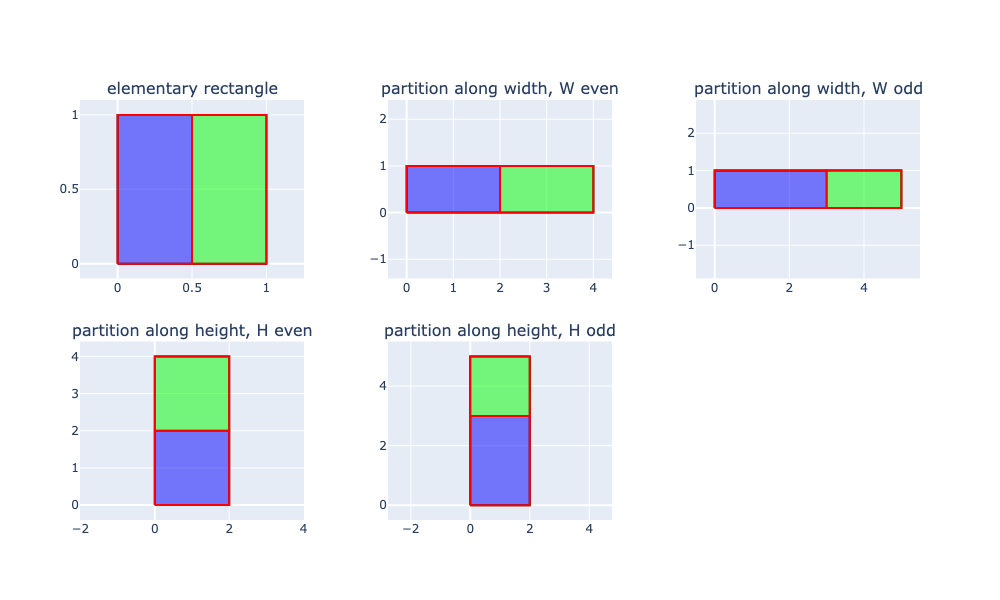}
	\caption{
	Examples of the five cases of the midway partition of a dyadic rectangle.
	The blue  rectangle is $\square_A$ and the green rectangle is $\square_B$.
	}
	\label{fig:2D_midway}
\end{figure}

\end{definition}

\begin{lemma}\label{lem:midway-balanced}
	Let $\square$ be a dyadic rectangle and $\{\square_A,\square_B\}$ its midway partition.
	Then $\{\square_A,\square_B\}$ is balanced.
\end{lemma}

\begin{proof}
If $\square$ is elementary, then the claim is clear.
Suppose henceforth that $\square$ is not elementary.

\textit{Case 1:} $E(\square)\leq 9/2$.

\textit{Case 1a:} $\bfW= \bfH$, so that $E(\square)=\bfW/\bfH =1$.
	\begin{itemize}
		\item If $\bfW$ is even, then $\square_A$ and $\square_B$ have width $\bfw/2$ and height $\bfw=\bfh$, so $E(\square_A)=E(\square_B)=2$.
		\item If $\bfW$ is odd, then  $\square_A$ has width $2^{-N}(\bfW+1)/2$ and height $\bfw=\bfh$, so $E(\square_A)=2\bfW/(\bfW+1) \leq 6/4$ and $E(\square_B) = 2\bfW/(\bfW-1)\leq 3$, where the inequalities follow from $\bfW\geq 3$ since $\square$ is not elementary.
	\end{itemize}

\textit{Case 1b:} $\bfW>\bfH$, so that $E(\square)=\bfw/\bfh \leq 9/2$.
\begin{itemize}
	\item If $\bfW$ is even, then $\square_A$ and $\square_B$ have width $\bfw/2$ and height $\bfh$.
Remark that $\bfw/(2\bfh) = \frac12 E(\square)$ and $2\bfh/\bfw \leq 2$, which implies $E(\square_A),E(\square_B) \leq 3$.

\item If $\bfW$ is odd, then $\square_A$ has width $2^{-N}(\bfW+1)/2$ and height $\bfh$.
Remark that
\begin{equ}
(\bfW+1)/(2\bfH)\leq \frac94+\frac{1}{2\bfH}\leq \frac{11}{4} < 3
\end{equ}
and $2\bfH/(\bfW+1) \leq 2$, which implies $E(\square_A) < 3$,
and likewise
$(\bfW-1)/(2\bfH) < 3$ and $2\bfH/(\bfW-1) \leq 2\bfW/(\bfW-1) \leq 3$ due to $\bfW\geq 3$ since $\square$ is not elementary, which implies $E(\square_B) \leq 3$.
\end{itemize}

The case $\bfH>\bfW$ follows by symmetry from the case $\bfW>\bfH$.
This completes the proof of \textit{Case 1}.

\textit{Case 2:} $E(\square)>9/2$.
Without loss of generality, suppose that $\bfW>\bfH$.
\begin{itemize}
\item If $\bfW$ is even, then $\square_A$ and $\square_B$ have width $\bfw/2$ and height $\bfh$ and thus, since $\bfw>9\bfh/2$,
\begin{equ}
E(\square_A) = E(\square_B) = \frac{\bfw}{2\bfh} = \frac12 E(\square)\;.
\end{equ}

\item If $\bfW$ is odd, then $\square_A$ has width $2^{-N}(\bfW+1)/2$ and height $\bfh$ and thus,
\begin{equ}
E(\square_A) = \frac{\bfW+1}{2\bfH} = \frac12 E(\square) + \frac{1}{2\bfH} < \frac23 E(\square)\;,
\end{equ}
where we used that $\bfH\geq 1$ and thus
\begin{equ}
\frac{1}{2\bfH} \leq \frac12 < \frac34 = \frac16 \times \frac{9}{2} < \frac{1}{6} E(\square)\;.
\end{equ}

Similarly $\square_B$ has width $2^{-N}(\bfW-1)/2$ and height $\bfh$ and thus, since $\bfW\geq 3$ and therefore $(\bfW-1)/2 \geq \bfW/3 \geq \bfH$,
we have
\begin{equ}
E(\square_B) = \frac{\bfW-1}{2\bfH} < E(\square_A) < \frac23 E(\square)\;.
\end{equ}
\end{itemize}

This completes the proof in \textit{Case 2}. 
\end{proof}

\begin{proof}[Proof of Existence and Continuity in \Cref{thm:2D_sewing}]
Consider $\square \in \RECT^N$ with midway partition $\{\square_A,\square_B\}$.
Define the scale $n$ approximation of $\widehat\SD$ on $\square$ iteratively by $\widehat\SD^0 \square = \widehat\SD \square$ and
\begin{equ}[eq:Wn_def]
\widehat\SD^n\square =
\begin{cases}
\actionGG_{\PD(\bottomboundary{A})}( \widehat\SD^{n-1}(\square_B) )\, \widehat\SD^{n-1}(\square_A)\quad &\text{ if } \square_A\square_B=\square
\\
\widehat\SD^{n-1}(\square_A)\, \actionGG_{\PD(\leftboundary{A})} \widehat\SD^{n-1}(\square_B) 
		\quad &\text{ if } \genfrac{}{}{0pt}{1}{\square_B}{\square_A}=\square\;.
\end{cases}
\end{equ}
To lighten notation, we let $\rightarrow$ and $\uparrow$ denote $\bottomboundary{A}$ and $\leftboundary{A}$ from \Cref{def:Stokes_Chen} respectively.

Take $\eps>0$ such that
\begin{equ}[eq:eps_gamma_def]
K \coloneqq L(1+\eps)^2 < 1
\end{equ}
where $L<1$ is as in \Cref{def:subcontrol} associated to $F$.
Consider $\delta>0$ such that~\eqref{eq:rho_assump} and~\eqref{eq:G_act_lip} hold.
We let $\diam(X)$ denote the diameter of a set $X\subset\R^2$.

\textbf{Claim 1:} There exists $\kappa>0$ with the following property. If $\diam (\square)<\kappa$ then for all $n\geq 0$
\begin{equ}\label{eq:Wn_Cauchy}
\rho(\widehat\SD^{n+1}\square, \widehat\SD^{n}\square) \leq F(\square) K^n \;.
\end{equ}
To prove the claim, we proceed by induction. The case $n=0$ follows by assumption~\eqref{eq:approx_Chen1}-\eqref{eq:approx_Chen2}.

For the inductive step, we will assume that $\square_A\square_B=\square$ (the case that $\genfrac{}{}{0pt}{1}{\square_B}{\square_A}=\square$ is identical).
We then have the telescoping bound
\begin{equs}[eq:telescope]
\rho(\widehat\SD^{n+1} \square, \widehat\SD^{n} \square)
&\leq
\rho(
\actionGG_{\PD(\rightarrow)}( \widehat\SD^n {\square_B} )\, \widehat\SD^{n} {\square_A},
\actionGG_{\PD(\rightarrow)}( \widehat\SD^n {\square_B} )\, \widehat\SD^{n-1} {\square_A})
\\
&\qquad +
\rho(
\actionGG_{\PD(\rightarrow)}( \widehat\SD^{n} {\square_B} )\, \widehat\SD^{n-1} {\square_A},
\actionGG_{\PD(\rightarrow)}( \widehat\SD^{n-1} {\square_B} )\, \widehat\SD^{n-1} {\square_A})\;.
\end{equs}
We also have the telescoping bound
\begin{equ}
\rho(
\actionGG_{\PD(\rightarrow)}( \widehat\SD^{n} {\square_B} ) , e)
\leq
\rho(\actionGG_{\PD(\rightarrow)}( \widehat\SD^{0} {\square_B} ), e) +
\sum_{i=1}^{n} \rho(\actionGG_{\PD(\rightarrow)} \widehat\SD^{i} {\square_B}, \actionGG_{\PD(\rightarrow)} \widehat\SD^{i-1} {\square_B} )\;.
\end{equ}
Since $\widehat\SD$ is continuous and $\widehat\SD (\bar\square)=e$ for any $\bar\square\in \Rect$ with $\diam(\bar\square)=0$, we have $\rho(\widehat\SD^{0} \square_B,e) \leq c\delta$ for $\diam(\square)$ sufficiently small, where $c$ is a small constant.
By telescoping and our inductive hypothesis,
we have for all $i=1,\ldots, n-1$,
\begin{equ}
\rho(\widehat\SD^{i} {\square_B},\widehat\SD^0 {\square_B}) \leq 
\sum_{j=1}^i \rho( \widehat\SD^{j} {\square_B}, \widehat\SD^{j-1} {\square_B} )
\leq F(\square)/(1-K)\;.
\end{equ}
Hence $\rho(\widehat\SD^{i} {\square_B},e)\leq c\delta$ whenever $\diam(\square)$ is sufficiently small since $\lim_{|\square|\to0}F(\square)=0$.
Furthermore, since $\PD$ is locally uniformly continuous and is a $1$-cocycle, $\rho_G(\PD(\rightarrow),e_G) \leq \delta$ for $\diam(\square)$ sufficiently small.
Therefore, by~\eqref{eq:G_act_lip}, for all $i=1,\ldots, n$,
\[
\rho(\actionGG_{\PD(\rightarrow)} \widehat\SD^{i} {\square_B}, \actionGG_{\PD(\rightarrow)} \widehat\SD^{i-1} {\square_B} ) \leq (1+\eps)\rho(\widehat\SD^{i} {\square_B},\widehat\SD^{i-1} {\square_B})\;.
\]
In conclusion, applying once more our inductive hypothesis,
\begin{equ}
\rho(
\actionGG_{\PD(\rightarrow)}( \widehat\SD^{n} {\square_B} ) , e)
\leq
\rho(\actionGG_{\PD(\rightarrow)}( \widehat\SD^{0} {\square_B} ), e)
+ (1+\eps) \sum_{i=1}^{n} K^{- i} F(\square) \leq \delta
\end{equ}
where the final bound holds once $\diam(\square)$ is sufficiently small.

The same conclusion applies with $\actionGG_{\PD(\rightarrow)}( \widehat\SD^{n} {\square_B} )$ replaced by
$\widehat\SD^{n-1} {\square_A}$.

Now applying the assumptions~\eqref{eq:rho_assump} and~\eqref{eq:G_act_lip} on $\rho$ and again the inductive hypothesis, we conclude that the right-hand side of \eqref{eq:telescope} is bounded by
\begin{equ}
(1+\eps)^2 K^{n-1} (F(\square_A) + F(\square_B))\;.
\end{equ}
(The factor $(1+\eps)^2$ appears from applying \eqref{eq:rho_assump} then \eqref{eq:G_act_lip} to the first term on the right-hand side of \eqref{eq:telescope}.)

Since $\square_A,\square_B$ is a balanced partition by \Cref{lem:midway-balanced}
and since $F$ is a subcontrol,
\begin{equ}
F(\square_A)+F(\square_B)\leq L F(\square)\;.
\end{equ}
Since $K = (1+\eps)^2 L$, it follows that the right-hand side of \eqref{eq:telescope} is bounded by
$
K^n F(\square)
$.
This concludes the inductive step and the proof of Claim 1.

(\textbf{Remark:} the proof of Claim 1 only uses the Lipschitz estimates \eqref{eq:rho_assump} and \eqref{eq:G_act_lip}, uniform continuity of $\PD,\widehat\SD$, and approximate Chen \eqref{eq:approx_Chen1}-\eqref{eq:approx_Chen2}.)

By completeness of $\rho$, the previous claim implies the existence of the limit
\begin{equ}
\SD (\square) = \lim_{n\to\infty } \widehat\SD^n (\square)\;.
\end{equ}
Furthermore, the claim shows the bound \eqref{eq:approx_bound} for all $\square \in \RECT^N$, $N\geq 0$.

\textbf{Claim 2.} For all $N$ sufficiently large, $\feedbackGG \SD (\square) = \PD(\partial\square)$ for all elementary squares $\square\in\DD^N$.

The proof of the claim is essentially the same as (even simpler than) that of Claim 1, except that we compare $\feedbackGG \widehat\SD^n (\square)$ to $\PD(\partial \square)$,
so we only highlight the differences.

Let $\square\in\DD^N$.
By definition \eqref{eq:Wn_def}, $\widehat\SD^{2n}(\square)$ satisfies the recursion $\widehat\SD^0 (\square)=\widehat\SD (\square)$ and
\begin{equ}
\widehat\SD^{2n} (\square) = \actionGG_{\PD(\rightarrow)}( \widehat\SD^{2(n-1)} \square_3 )
\actionGG_{\PD(\rightarrow\sqcup \uparrow)}( \widehat\SD^{2(n-1)} {\square_4} )
\widehat\SD^{2(n-1)} {\square_A}
\actionGG_{\PD( \uparrow)}( \widehat\SD^{2(n-1)} {\square_B} )  \;,
\end{equ}
where we decompose $\square$ into four squares as
	\begin{align*}
		\square = \genfrac{}{}{0pt}{1}{\square_B \square_4}{\square_A \square_3}
	\end{align*}
	and the arrows join the lower left corner of $\square_A$ with that of the other squares.

Suppose $\bar L \in (0,1)$ is as in \Cref{def:subcontrol} associated to $\bar F$ and take $\eps>0$ such that $\bar K \coloneqq \bar L^2 (1+\eps)^2 < 1$.
Consider $\delta>0$ such that~\eqref{eq:rho_assump_G} and~\eqref{eq:G_Ad_lip} hold.
We claim that, similar to \eqref{eq:telescope}, for all $n\geq 0$,
\begin{equ}[eq:W_PD_dist]
	\rho_G( \feedbackGG \widehat\SD^{2n} \square, \PD(\partial \square)) \leq \bar F(\square) \bar K^{n}
\end{equ}
whenever $\diam (\square)$ is sufficiently small.

Indeed, proceeding by induction, the base case $n=0$ is true by assumption \eqref{eq:approx_Stokes}.
For the inductive step,
we have, by definition of $\widehat\SD^n$, the identity \eqref{eq:equivariant} $\feedbackGG m_g h = \Ad_g \feedbackGG h$,
the morphism property of $\feedbackGG$,
and Chen's identity \eqref{eq:PD_Chen} for $\PD$,
\begin{equs}
\rho_G(\feedbackGG \widehat\SD^{n} \square, \PD (\partial \square))
&=
\rho_G
\Big\{
\Ad_{\PD(\rightarrow)} (\feedbackGG \widehat\SD^{2(n-1)} {\square_3} )
\Ad_{\PD(\rightarrow\sqcup \uparrow)} (\feedbackGG \widehat\SD^{2(n-1)} {\square_4} )
\\
&\qquad\qquad
\feedbackGG \widehat\SD^{2(n-1)} {\square_A}
\Ad_{\PD(\uparrow)}(\feedbackGG \widehat\SD^{2(n-1)} {\square_B} )
\;,
\\
&\qquad
\Ad_{\PD(\rightarrow)} (\PD ({\partial \square_3}) )
\Ad_{\PD(\rightarrow\sqcup \uparrow)} (\PD ({\partial \square_4}) )
\PD {(\partial \square_A)}
\Ad_{\PD(\uparrow)}(\PD ({\partial \square_B}) )
\Big\}
\\
&=: \rho_G \{\Delta_1,\Delta_2 \}\;.
\end{equs}
By the same argument as in the proof of Claim 1, combined with continuity of $\feedbackGG$ at $e$
and the facts that $\lim_{\diam(\square)\to 0}\widehat\SD (\square)=e$
and $\lim_{|\gamma|_{BV}\to 0} \PD(\gamma)= e_G$,
we can ensure that all products of terms comprising $\Delta_1$ and $\Delta_2$ are in the $\delta$-ball of $e_G$
once $\diam(\square)$ is sufficiently small.

It follows that, applying a similar telescoping sum as in~\eqref{eq:telescope}, the inductive hypothesis, and our assumptions on $\rho$~\eqref{eq:rho_assump_G} and~\eqref{eq:G_Ad_lip},
we obtain
\begin{equ}
\rho_G(\feedbackGG \widehat\SD^{n} \square, \PD (\partial \square)) \leq (1+\eps)^2 \bar K^{n-1} (\bar F(\square_A)+\cdots+\bar F(\square_4)) \;.
\end{equ}
Observe that $\bar F(\square_A)+\cdots+\bar F(\square_4)\leq \bar L^2 \bar F(\square)$ by the fact that $\bar F$ is a subcontrol, and thus, by our choice $\bar K = \bar L^2(1+\eps)^2$, we obtain
\begin{equ}
\rho_G(\feedbackGG \widehat\SD^{n} \square, \PD (\partial \square)) \leq \bar K^{n} \bar F(\square)\;.
\end{equ}
This completes the proof of \eqref{eq:W_PD_dist} and thus of Claim 2 by continuity of $\feedbackGG$.

(\textbf{Remark:} the proof of Claim 2 uses only the Lipschitz estimates~\eqref{eq:rho_assump_G} and~\eqref{eq:G_Ad_lip}, the morphism property and continuity of $\feedbackGG$, Chen's identity \eqref{eq:PD_Chen} for $\PD$, approximate Stokes \eqref{eq:approx_Stokes}, and uniform continuity of $\PD,\widehat\SD$.)

\textbf{Claim 3.} $\feedbackGG \SD (\square) = \PD(\partial\square)$ on any dyadic rectangle $\square\in\RECT^N$, $N\geq 0$.

Observe that the midway partition (\Cref{def:midway}) reduces $\square\in\RECT^N$ to two $2^{-N-1}$-dyadics $\square_A,\square_B$ if and only if $\square$ is an elementary square.
Moreover, in this case, the midway partitions of $\square_A$ and $\square_B$ are again elementary $2^{-N-1}$-dyadic squares.
Therefore, for any $\square\in\RECT^N$, there exists $m$ sufficiently large such that, for all $n \geq m$,
\begin{equ}[eq:pre-lassos]
\widehat\SD^{n} (\square) = \prod_{i=1}^{k} \actionGG_{\PD(\tilde c_i)} \widehat\SD^{n-m}(\tilde b_i)
\end{equ}
where $\tilde b_i, \tilde c_i \in\PATHS^{N_i}$ with $\tilde b_i$ the boundary of either an elementary square or a composition of two elementary squares
and with $\prod_{i=1}^{k} \tilde c_i \tilde b_i \tilde c_i^{-1} = \partial \square$.
(Note, that different $\tilde b_i$ may be boundaries of dyadic rectangles of different scales.)

Taking $n$ sufficiently large, we can furthermore find $r_i\geq 0$ for $i=1,\ldots, k$ such that
\begin{equ}[eq:lassos]
\widehat\SD^{n} (\square) = \prod_{i=1}^{k} \actionGG_{\PD(c_i)} \widehat\SD^{n-m-r_i}(b_i)
\end{equ}
where now $c_i b_i c_i^{-1} \in \Lasso^{N^*}$ are lassos (\Cref{def:lassos}) of the \emph{same} dyadic scale $N^*\geq N$ and with $\prod_{i=1}^{k}  c_i  b_i  c_i^{-1} = \partial \square$.
See \Cref{fig:2D_midway_partition_iteration} for an example.

\begin{figure}
	\centering
	\includegraphics[width=0.8\textwidth]{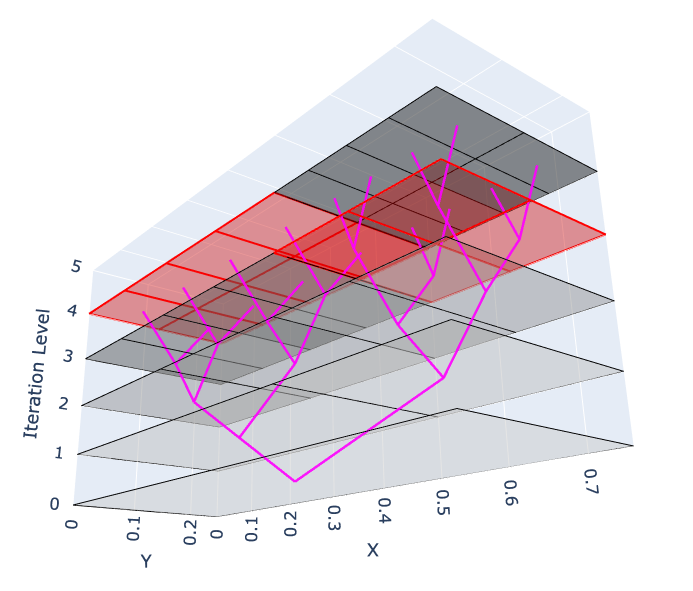}
	\caption{The first $4$ iterations of the midway partition of
	the dyadic rectangle $[0,\frac34]\times [0,\frac14] \in \Rect^0$.
	The rectangles colored red are elementary
	squares in $\Rect^1$ that partition the original rectangle.
	These correspond to $b_i$ in \eqref{eq:lassos}.
	}
	\label{fig:2D_midway_partition_iteration}
\end{figure}

Taking the limit $n\to\infty$, we obtain
\begin{equ}[eq:lasso_decomp]
	\SD (\square) = \prod_{i=1}^{k} \actionGG_{\PD(c_i)} \SD(b_i)\;.
\end{equ}
We can assume without loss of generality that $N^*$ is sufficiently large so that $\feedbackGG \SD(b_i) = \PD(b_i)$ by Claim 2.
Hence applying $\feedbackGG$ to both sides of \eqref{eq:lasso_decomp} and using \eqref{eq:equivariant} $\feedbackGG m_g h =\Ad_g \feedbackGG h$, Chen's identity for $\PD$, and $\prod_{i=1}^{k}  c_i  b_i  c_i^{-1} = \partial \square$,
it follows that $\feedbackGG \SD (\square) = \PD(\partial\square)$.

\textbf{Claim 4.} $\SD$ satisfies Chen for all dyadic $\square\in\RECT^N$, $N\geq 0$.

This claim is a consequence of Claim 3 and \Cref{lem:Stokes_to_Chen}.
More precisely, suppose $\square_A\square_B = \square$ is a subdivision of $\square\in\RECT^N$ into two dyadic rectangles.
Then, by the argument in the proof of Claim 3, there exist $N_1\geq N$ and lassos $c^1_i b^1_i (c^1_i)^{-1}\in\Lasso^{N_1}$ for $i=1,\ldots, k_1$, such that $\prod_{i=1}^{k_1} c_i^1 b_i^1 (c_i^1)^{-1} = \partial \square_A$
(these lassos are based at the lower left corner of $\square_A$, which is also the lower left corner of $\square$),
and
\begin{equ}[eq:SD_sq1]
\SD (\square_A) = \prod_{i=1}^{k_1} \actionGG_{\PD(c^1_i)} \SD(b^1_i)\;.
\end{equ}
Likewise there exist $N_2\geq N$ and lassos $c^2_i b^2_i (c^2_i)^{-1}\in \Lasso^{N_2}$ for $i=1,\ldots, k_2$, such that $\prod_{i=1}^{k_2}c_i^2 b_i^2 (c_i^2)^{-1} = \partial \square_B$
(these lassos are now based at the lower left corner of $\square_B$)
and
\begin{equ}[eq:SD_sq2]
	\SD (\square_B) = \prod_{i=1}^{k_2} \actionGG_{\PD(c^2_i)} \SD(b^2_i)\;.
\end{equ}
By taking further subdivisions, we may furthermore suppose that $N_1=N_2$, i.e.~the lassos are on the same dyadic scale.

Finally, consider the new lassos $\bar c^2_i b^2_i (\bar c^2_i)^{-1}$
where $\bar c^2_i = \rightarrow\sqcup c^2_i$ (these new lassos are based at the lower left corner of $\square$).
It follows that
\begin{equs}
\partial \square &= \rightarrow \sqcup (\partial \square_B) \sqcup (\leftarrow) \sqcup (\partial \square_A)
=
\Big(\prod_{i=1}^{k_2} \bar c_i^2  b_i^2 (\bar c_i^2)^{-1} \Big)
 \sqcup
 \Big(\prod_{i=1}^{k_1} c_i^1 b_i^1 (c_i^1)^{-1}\Big)\;.
\end{equs}

We thus obtain
\begin{equs}
\SD (\square)
= \prod_{i=1}^{k} \actionGG_{\PD(c_i)} \SD(b_i)
= \Big(\prod_{i=1}^{k_2} \actionGG_{\PD(\bar c^2_i)} \SD(b^2_i) \Big) \Big( \prod_{i=1}^{k_1} \actionGG_{\PD(c^1_i)} \SD(b^1_i) \Big)
= \actionGG_{\PD(\rightarrow)} (\SD \square_B) \SD \square_A\;,
\end{equs}
where we used \eqref{eq:lasso_decomp} in the first equality,
\Cref{lem:Stokes_to_Chen} in the second equality, and Chen's identity for $\PD$ and \eqref{eq:SD_sq1}-\eqref{eq:SD_sq2} in the final equality.
This proves horizontal Chen. The proof of vertical Chen is identical. This completes the proof of Claim 4.

To complete the proof of existence, it remains to extend $\SD$ to all rectangles such that \eqref{eq:approx_bound}, Chen, and Stokes are satisfied.
To this end, suppose $\square\in \RECT$ and let $\square_i\in\RECT^{N_i}$ be a sequence of dyadics that converges to $\square$.
Then the symmetric difference between $\square_i$ and $\square_j$ (as subsets of $\R^2$) can be written as a disjoint union of at most 8 dyadic rectangles whose area is vanishing as $i,j\to\infty$.
Therefore, by the approximation bound \eqref{eq:approx_bound} for dyadic rectangles, the fact that
\begin{equ}[eq:vanish_area_SD_F]
\lim_{|\square|\to 0} \rho(\widehat\SD \square,e)=\lim_{|\square|\to 0}F(\square)= 0\;,
\end{equ}
and the continuity of group and crossed module operations,
we obtain that $\SD({\square_i})\in H$ is a Cauchy sequence.
The limit $\SD(\square)\coloneqq\lim_{i\to\infty} \SD(\square_i)$ therefore exists.
It follows from continuity of $\PD$ and $\widehat\SD$ that $\SD$ satisfies \eqref{eq:approx_bound}, Chen, and Stokes,
concluding the proof of existence of $\SD$.

To prove the final claim, by again considering the symmetric difference between rectangles and using \eqref{eq:vanish_area_SD_F}, $\square\mapsto \SD (\square)$ is a continuous function on $\RECT$.
To prove the stated uniform continuity of $(\PD,\widehat\SD)\mapsto \SD$, consider $(\PD^{(1)},\widehat\SD^{(1)}),(\PD^{(2)},\widehat\SD^{(2)})\in\BB_1\times\BB_2$.
Suppose $\sup_{\square\in\RECT} \rho(\widehat\SD^{(1)}\square,\widehat\SD^{(2)} \square) < \eps$
and $\sup_{|\gamma|_{BV}\leq 1}\rho_G(\PD^{(1)}(\gamma),\PD^{(2)}(\gamma))<\eps$.
Then, for all $\square\in\RECT$ with sufficiently small diameter,
\begin{equ}
\rho(\actionGG_{\PD^{(1)}(\gamma)}\widehat\SD^{(1)} \square,\actionGG_{\PD^{(2)}(\gamma)}\widehat\SD^{(2)} \square) < \bar\eps
\end{equ}
where $\bar\eps\to0$ as $\eps\to0$ due to the joint (uniform) continuity of $(g,h)\mapsto \actionGG_g h$ around $(e_G,e)\in G\times H$.
It follows that,
for any dyadic $\square\in\Rect^N$, $\rho(\widehat\SD^{{(1)},n} \square,\widehat\SD^{{(2)},n} \square) \lesssim \bar\eps 2^{n}$
uniformly in $n,N\geq 0$.
On the other hand, if $\diam(\square)$ is sufficiently small, then by \eqref{eq:Wn_Cauchy},
$\rho(\widehat\SD^{{(1)},n} \square, \SD^{(1)} \square) \lesssim L^n F(\square)$
where $L<1$,
and likewise for $\widehat\SD^{(2)},\SD^{(2)}$.
Therefore $\rho(\SD^{(1)} \square,\SD^{(2)} \square)\lesssim L^n + \bar \eps 2^{n}$
uniformly in $n\geq 1$.
Taking $n$ as a small multiple of $\log \bar\eps$ we have $L^n + \bar \eps 2^{n}\to 0$.
This shows that $\rho(\widehat\SD^{{(1)}} \square, \SD^{(2)} \square) = o(1)$ as $\eps\downarrow0$ whenever $\diam(\square)$ is small.
The fact that $\rho(\widehat\SD^{{(1)}} \square, \SD^{(2)} \square) = o(1)$ for all $\square\in\RECT$ follows from Chen's identity
and local uniform continuity of $(g,h)\mapsto \actionGG_g h$.
\end{proof}

\begin{proof}[Proof of \Cref{thm:2cocycle_sewing}]
	By \Cref{thm:2D_sewing}, it suffices to show that $\widehat\SD$ satisfies the approximate Chen identities~\eqref{eq:approx_Chen1}-\eqref{eq:approx_Chen2} and approximate Stokes~\eqref{eq:approx_Stokes} identity for some subcontrols $F,\bar F$ respectively.
	Starting with Stokes, consider an elementary square $\square\in\DD^N$. Then
	\begin{equ}
	\feedbackGG \widehat\SD (\square) = \feedbackGG \exp\Big(\int_\square \beta \Big) = \exp \Big( \int_{\square} \feedbackAA \beta \Big) = \exp \Big( \int_\square \{d\alpha + \alpha\wedge \alpha\} \Big)\;,
	\end{equ}
	where the final equality is due to the vanishing fake curvature condition \eqref{eq:vanishfakecurvature}.
	By standard Euler estimates for ODEs and the BCH formula, $\rho_G(\PD^\alpha(\partial\square), e^{\int_\square \{d\alpha + \alpha\wedge \alpha\} }) \lesssim |\square|^{3/2}$.
	Therefore, approximate Stokes is satisfied with the subcontrol $\bar F$ taken as a sufficiently large multiple of the subcontrol from \Cref{prop:subcontrol} with $\theta=\lambda=\zeta=3/2$.

	To show approximate horizontal Chen~\eqref{eq:approx_Chen1}, consider $\square\in\Rect$ and a balanced partition $\square_A\square_B = \square$.
	Then, recalling the differential $\actionGA_g\in\Aut(\h)$ of $\actionGG_g$,
	\begin{equ}
	\actionGG_{\PD(\rightarrow)}( \widehat\SD {\square_B} ) = \exp \Big\{\actionGA_{\PD(\rightarrow)} \Big(\int_{\square_B}\beta \Big)\Big\}
	= \exp \Big\{\int_{\square_B}\beta + O(\bfw_A |\square_B|)\Big\}
	\end{equ}
	where the final equality is due to $\actionGA_g - \id = O(\bfw_A) \in \Aut(\h)$.
	Therefore
	\begin{equ}
	\actionGG_{\PD(\rightarrow)}( \widehat\SD {\square_B} ) \, \widehat\SD {\square_A} = \exp\Big(\int_{\square_A}\beta + \int_{\square_B}\beta + O(\bfw_A |\square_B|) + O(|\square|^2)\Big)
	\end{equ}
	where we used the BCH formula.
	Since $|\square| \lesssim |\square_B|$ by \Cref{lem:areas} and clearly $|\square|^{1/2} \lesssim \bfw_A$, we have $|\square|^2 = O(\bfw_A |\square_B|)$.
	Furthermore $\widehat\SD \square = \exp({\int_\square \beta}) = \exp({\int_{\square_A} \beta+\int_{\square_B} \beta})$ by additivity of integrals.
	It follows that
	\begin{equ}
	\rho(\widehat\SD \square , \actionGG_{\PD(\rightarrow)}( \widehat\SD {\square_B} ) \, \widehat\SD {\square_A}) \lesssim \bfw_A |\square_B|\;.
	\end{equ}
In conclusion, $\widehat\SD$ satisfies~\eqref{eq:approx_Chen1} with $F$ taken as a sufficiently large multiple of the subcontrol from \Cref{prop:subcontrol} with $\theta=\frac32$ and $\lambda=2$, $\zeta=1$, which shows \eqref{eq:strong_bound}.

Finally, it follows from \Cref{rem:weaker_F} and \Cref{prop:subcontrol} with $\theta=\frac32$ and $\lambda=3$, $\zeta=0$ that $\SD$ is furthermore the unique $2$-cocycle verifying the weaker bound
$\rho(\SD_\square,\widehat\SD_\square) \lesssim \bfa^{\lambda}$ for any $\lambda\in(2,3]$.
\end{proof}


\subsection{Extension theorem}

In addition to giving a proof of \Cref{thm:2cocycle_sewing}, another important application of the Sewing Lemma (\Cref{thm:2D_sewing}) is an extension theorem for (possibly rough) surfaces.
We phrase this extension theorem in the general setting of homogenous groups, of which the nilpotent (truncated) crossed module from \Cref{ex:free_nilpotent} is a special case.

We believe this general setting has potential to define the surface signature taking values in different algebraic structures, analogous to signatures of branched rough paths \cite{gubinelli_10_branched} in the Grossman--Larson algebra, planarly branched rough paths \cite{Curry_et_al_20_planarly_BRP} in the Munthe-Kaas--Wright algebra, multi-index rough paths \cite{Linares_23_multiindex} in the free Novikov algebra, or other graded Hopf algebras~\cite{BZ22,Manchon_25_RP_Hopf}.


\subsubsection{Homogenous groups}
\label{sec:hom_grps}

For details on homogenous groups, see~\cite{FollandStein82,HebischSikora90}.
Our notation is taken from~\cite{Chevyrev18}.

Throughout the section, we fix a homogenous group $G$, i.e.~a nilpotent, connected, and simply connected Lie group endowed with a one-parameter family of dilations (group automorphisms) $(\delta_\lambda)_{\lambda > 0}$, which, upon identifying $G$ with its Lie algebra $\g$ via the diffeomorphism $\exp \colon \g \to G$, is given by
\[
\delta_\lambda(u_i) = \lambda^{d_i} u_i
\]
for a basis $u_1,\ldots, u_m$ of $\g$ and real numbers $d_m \geq \cdots \geq d_1 \geq 1$ which we call the \emph{degrees} of $G$.
We equip $G$ with a ``homogeneous norm'', i.e.~a continuous map $\|\cdot\| \colon G\to [0,\infty)$ for which $\|\delta_\lambda x\| = \lambda\|x\|$ and $\|x\|=0$ if and only if $x=e_G$.
We assume that $\|\cdot\|$ is sub-additive, i.e.~$\|x y\|\leq \|x\|+\|y\|$ (see~\cite{HebischSikora90} for the existence of such a norm)
and equip $G$ with the corresponding left-invariant metric $d(x,y) = \|{x^{-1}y}\|$ on $G$.

Unless otherwise stated, we identify $G$ with $\g$ via the $\exp$-map and write $x = \sum x^{[i]} u_i$ for $x \in G$.
For a multi-index $\alpha = (\alpha^1,\ldots, \alpha^m)\in\ZZ^m$, $\alpha^i \geq 0$, we define $\deg(\alpha) = \sum_{i=1}^m \alpha^i d_i$, and for $x \in G$, write $x^{[\alpha]} = (x^{[1]})^{\alpha^1}\cdots(x^{[m]})^{\alpha^m}$.
By the BCH formula, for all $i \in \{1,\ldots, m\}$ there exist constants $C^i_{\alpha,\beta}$ such that
\begin{equation}
\label{eq:CBH}
	(xy)^{[i]} = x^{[i]} + y^{[i]} + \sum_{\alpha,\beta} C^i_{\alpha,\beta} x^{[\alpha]} y^{[\beta]},
\end{equation}
where the (finite) sum runs over all non-zero multi-indexes $\alpha,\beta$ such that $\deg(\alpha) + \deg(\beta) = d_i$.

\begin{example}[Graded groups]\label{ex:gradedLie}
Recall that a Lie group $G$ is called graded if its Lie algebra is endowed with a decomposition
\begin{equation}
\label{eq:LieAlgDecomp}
	\g = \g^1 \oplus \cdots \oplus \g^N
\end{equation}
such that $[\g^i,\g^j] \subset \g^{i+j}$, where $\g^k = \{0\}$ for $k > N$ (we allow the possibility that $\g^k = \{0\}$ for some $k \leq N$). Every graded Lie group can be equipped with a natural family of dilations $(\delta_{\lambda})_{\lambda > 0}$, and thus a homogeneous structure, for which $d_1,\ldots, d_m$ are integers determined by $\delta_\lambda(u) = \lambda^{k}u$ for all $u \in \g^k$.

Recall also that a graded Lie group $G$ is called a step-$N$ Carnot group (or stratified group in the terminology of~\cite{FollandStein82}) if the decomposition~\eqref{eq:LieAlgDecomp} further satisfies $[\g^i,\g^j] = \g^{i+j}$, where $\g^k = \{0\}$ for $k > N$.
Every Carnot group is a homogeneous group for which the metric $d$ can be taken as the Carnot--Carath{\'e}odory distance~\cite[p.~38]{Baudoin04}.
\end{example}

Due to~\eqref{eq:CBH}, for every $1\leq n \leq m$, the set $K_n\coloneqq \mathrm{span}\{u_{n+1},\ldots, u_m\}$ is a normal subgroup of $G$ and therefore the projection map
\begin{equ}[eq:pi_proj]
\pi_n\colon G \to G^{[n]} \coloneqq G/K_n
\end{equ}
is a group homomorphism.
Furthermore $G^{[n]}$ carries a canonical homogenous group structure and we identify $G^{[n]}$ with $\mathrm{span}\{u_1,\ldots, u_n\}$ equipped with the obvious group law given by the BCH formula.


\subsubsection{1D extension}

Consider a homogenous group $G$ as in \Cref{sec:hom_grps}.
For a function $X\colon[0,T]\to G$ and $\zeta \in [0,1]$, define the $\zeta$-H\"older ``norm'' $\|X\|_\zeta = \sup_{s\neq t} |t-s|^{-\zeta}d(X_s,X_t)$.

For comparison with the 2D (surface) case, we first provide an extension theorem for paths in $G^{[n]}$ of sufficient regularity.
The statement and proof is very similar to that of~\cite[Prop.~4.5]{BZ22}, although our setting of homogenous groups is more general (vs. the Hopf algebra setting of \cite{BZ22}).

\begin{lemma}[1D Extension]\label{lem:extension_1D}
Suppose $1\leq n< m$ and $\zeta\in (0,1]$ such that $d_{n+1}\zeta>1$.
Let $X\colon [0,T]\to G^{[n]}$ be a path with $\|X\|_\zeta < \infty$.

Then there exists a unique path $Y\colon [0,T]\to G$ such that $\pi_n(Y)=X$ and $\|Y\|_\zeta<\infty$.

Furthermore, $\|Y\|_\zeta\lesssim \|X\|_{\zeta}$,
where the proportionality constant depends only on $\zeta$ and $G$, and the map $X\mapsto Y$ is locally uniformly continuous in the $\zeta$-H\"older metric $d_\zeta(X,\bar X) \coloneqq \sup_{s,t} |t-s|^{-\zeta}d(X_{s,t},\bar X_{s,t})$ where we denote $X_{s,t} = X_s^{-1}X_t$.
\end{lemma}

\begin{proof}
By induction, it suffices to consider $n=m-1$.
By~\eqref{eq:CBH}, the kernel of $\pi_{n}$ is contained in the centre of $G$.


The proof is done once we show the existence of a unique two-parameter map $R_{s,t}\in\R$ such that $|R_{s,t}|\lesssim |t-s|^{\zeta d_m}$ and, for $s<u<t$,
\[
\delta R_{s,u,t} \coloneqq R_{s,u} + R_{u,t} - R_{s,t} =
-F_{s,u,t}
\coloneqq -\sum_{\alpha,\beta} C^m_{\alpha,\beta} X_{s,u}^{[\alpha]} X_{u,t}^{[\beta]}
= -(X_{s,u}\cdot_G X_{u,t})^{[m]}\;,
\]
where we write $X_{s,t} = X_{s}^{-1}X_t$ (the inverse and product are in $G^{[n]}$) and treat $X_{s,u},X_{u,t}\in G$ by setting the $u_m$-component to $0$,
and $\cdot_G$ denotes the product in $G$.
Indeed, setting $Y_{s,t} = X_{s,t} + R_{s,t}u_m$,
leads to
\[
	(Y_{s,u} \cdot_G Y_{u,t})^{[m]} 
	= R_{s,u} + R_{u,t} + \sum_{\alpha,\beta} C^m_{\alpha,\beta} X_{s,u}^{[\alpha]} X_{u,t}^{[\beta]} 
	= R_{s,t} 
	= Y_{s,t}^{[m]}\;.
\]
Therefore $Y_{s,t}$ satisfies Chen's identity and is therefore the increment of a path $Y$.

To find $R$, by the classical (linear) sewing lemma (see~\cite[Thm.~3.1]{BZ22}, also \cite{Gubinelli_04_controlling,Feyel_06_sewing}), since $|F_{s,u,t}|\lesssim |t-s|^{\zeta d_m}\|X\|_{\gamma}$,
it suffices to show that there exists $A\colon (s,t)\mapsto A_{s,t}\in \R$ such that $\delta A = F$, or equivalently, that
\begin{equ}
\delta F_{s,u,v,t}
\coloneqq F_{s,u,v} - F_{s,u,t} + F_{s,v,t} - F_{u,v,t} = 0
\end{equ}
for all $(s,u,v,t)$.
To this end, observe that
\[
(X_{s,u}\cdot_G X_{u,t}\cdot_G X_{t,s})^{[m]}
= (X_{s,u}\cdot_G X_{u,t})^{[m]} + X_{t,s}^{[m]}
+ \sum_{\alpha,\beta}C^m_{\alpha,\beta}X_{s,t}^{[\alpha]} X_{t,s}^{[\beta]}
=F_{s,u,t}
\]
since the sum in the middle term vanishes.
More generally, if $x_1,\ldots, x_k\in G^{[n]}$ with $x_1 \cdots x_k=e_{G^{[n]}}$, then
\[
	(x_1\cdot_G\ \cdots\ \cdot_G x_k)^{[m]} 
	= (x_1 \cdot_G\ \cdots\ \cdot_G x_{k-1})^{[m]}
	=(x_2\cdot_G\ \cdots\ \cdot_G x_k)^{[m]} \;,
\]
which means the expression is invariant under cyclic permutations.
Therefore
\begin{align*}
\delta F_{s,u,v,t}
&= F_{s,u,v} - F_{s,u,t} + F_{s,v,t} - F_{u,v,t}
\\
&= (X_{s,u}\cdot_G X_{u,v})^{[m]}
-(X_{s,u}\cdot_G X_{u,t})^{[m]}
+(X_{s,v}\cdot_G X_{v,t})^{[m]}
-(X_{u,v}\cdot_G X_{v,t})^{[m]}
\\
&= (X_{s,u}\cdot_G X_{u,v} \cdot_G X_{v,s})^{[m]}
+(X_{t,u}\cdot_G X_{u,s} \cdot_G X_{s,t})^{[m]}
\\
&\quad
+(X_{s,v}\cdot_G X_{v,t}\cdot_G X_{t,s})^{[m]}
+(X_{t,v}\cdot_G X_{v,u} \cdot_G X_{u,t} )^{[m]}\;.
\end{align*}
We now use the cyclic property of the products, and the fact that their products in $G^{[n]}$ are all $e_G$ as well as 
\[
	(x_1\cdot_G\ \cdots\ \cdot_G x_k)^{[m]} 
		+ (y_1\cdot_G\ \cdots\ \cdot_G y_\ell)^{[m]} 
			= (x_1\cdot_G\ \cdots\ \cdot_G y_\ell)^{[m]}
\]
whenever $x_1\cdots x_k = e_{G^{[n]}}$ and $y_1 \cdots y_\ell = e_{G^{[n]}}$, to see that the above sum is equal to
\begin{align*}
&(X_{s,u}\cdot_G X_{u,v} \cdot_G X_{v,t}\cdot_G X_{t,s})^{[m]}
+(X_{t,u}\cdot_G X_{u,s} \cdot_G X_{s,t})^{[m]}
+(X_{t,v}\cdot_G X_{v,u} \cdot_G X_{u,t} )^{[m]}
\\
&= (X_{s,u}\cdot_G X_{u,v} \cdot_G X_{v,t}\cdot_G X_{t,u}\cdot_G X_{u,s})^{[m]}
+(X_{t,v}\cdot_G X_{v,u} \cdot_G X_{u,t} )^{[m]}
\\
&= (X_{u,v} \cdot_G X_{v,t}\cdot_G X_{t,u}\cdot_G X_{u,t} \cdot_G X_{t,v} \cdot_G X_{v,u})^{[m]}
\\
&= 0\;.
\end{align*}
This shows that $\delta F=0$ and concludes the proof of existence and uniqueness of $Y$.
The claim about continuity of $X\mapsto Y$ follows from the boundedness of the linear sewing map $F\mapsto R$,
the details of which we leave to the reader.
\end{proof}

\begin{remark}
In the above proof, a (non-unique) $A$ such that $\delta A= F$ is
\[
A_{s,t} = \sum_{\alpha,\beta} C^m_{\alpha,\beta} X_{0,s}^{[\alpha]} X_{s,t}^{[\beta]} = (X_{0,s}\cdot_G X_{s,t})^{[m]}.
\]
\end{remark}


\subsubsection{2D extension}

Consider now two homogenous groups $G,H$.
We denote $\dim(H)=m$ and $\dim(G)=\ell$
and let $1\leq d_1\leq \cdots \leq d_m$ be the degrees of $H$
with corresponding basis $u_1,\ldots, u_m$ of $\h$,
and let $1\leq c_1 \leq \cdots \leq c_\ell$ be the degrees of $G$
with corresponding basis $v_1,\ldots, v_\ell$ of $\g$.
We follow the notation of \Cref{sec:hom_grps}
(with $H$ here playing the role of $G$ from \Cref{sec:hom_grps}).
In particular, we identify $G$ with $\g$ and $H$ with $\h$ via the $\exp$ maps.

We suppose that $H \stackrel{\feedbackGG}\to G \stackrel{\actionGG}\to \Aut(H)$ is a crossed module where $\feedbackGG$ and $\actionGG$ are continuous and graded in the following sense.
In addition to the BCH formula~\eqref{eq:CBH},
we assume that there exist constants $B^i_{\alpha,\beta}$
such that for every $i=1,\ldots, m$
\begin{equ}[eq:act_hom]
	(\actionGG_{g} h)^{[i]} = h^{[i]} + \sum_{\alpha,\beta} B^i_{\alpha,\beta} g^{[\alpha]} h^{[\beta]}
\end{equ}
where the sum is over all non-zero multi-indexes
$\alpha = \{\alpha^1,\ldots,\alpha^\ell\}$
and $\beta = \{\beta^1,\ldots,\beta^m\}$
such that
\[
	\deg(\alpha)+\deg(\beta) = d_i
\]
where $\deg(\alpha)=\sum_{j=1}^\ell \alpha^j c_j$
and $\deg(\beta) = \sum_{k=1}^m \beta^k d_k$.
We similarly suppose that
\begin{equ}[eq:feedback_GH_assump]
	\feedbackGG u_k \in \g_{d_k}
\end{equ}
where $\g_{\zeta}$ for $\zeta>0$ is the subspace of all $X\in G\simeq \g$ for which $\delta_\lambda X = \lambda^{\zeta} X$.
We let $\pi_n \colon H\to H^{[n]}$
denote the projection~\eqref{eq:pi_proj}
with kernel $K_n\coloneqq \mathrm{span}\{u_{n+1},\ldots, u_m\}$.
Assumption \eqref{eq:feedback_GH_assump} implies that, for all $1\leq n \leq m$,
$\feedbackGG K_n$ is a normal subgroup of $G$.
We define the quotient map
\begin{equ}
\proj_{\leq n} \colon G \to G^{\leq n} \coloneqq G / \feedbackGG K_n\;.
\end{equ}
Our assumptions on
$\actionGG,\feedbackGG$ imply that
\begin{equ}
H^{[n]} \stackrel{\feedbackGG}\to G^{\leq n} \stackrel{\actionGG}\to \Aut(H^{[n]})
\end{equ}
forms a crossed module for every $1\leq n \leq m$.


Recall the metrics we put on $\Paths$ and $\Rect$ from \Cref{conv:metrics}.

\begin{lemma}[2D Extension]\label{lem:extension_2D}
Consider a $1$-cocycle $\PD\colon \Paths \to G$ which is locally uniformly continuous. Suppose further that, for some $\zeta\in (0,1]$,
\begin{equ}[eq:PD_bound]
	|\PD(\gamma)^{[i]}| \lesssim |\gamma|^{\zeta c_i}_{BV} \quad \text{uniformly over $\gamma\in\Paths$ with $|\gamma|_{BV}<1$.}
\end{equ}
Let $1\leq n\leq m$ and $\widehat \SD \colon \Rect \to H^{[n]}$ be continuous. Suppose that $\widehat\SD$ is a $2$-cocyle up to level $n$, i.e, $\widehat\SD$ satisfies the level-$n$ Stokes identity
\begin{equ}
\feedbackGG \widehat \SD({\square}) = \proj_{\leq n} \PD(\partial \square)
\end{equ}
and similarly for level-$n$ Chen identities (for the crossed module $H^{[n]} \stackrel{\feedbackGG}\to G^{\leq n} \stackrel{\actionGG}\to \Aut(H^{[n]}$).

Finally, consider $\eps\in(0,\frac12]$ and suppose that
\begin{equ}[eq:dn_gamma]
\zeta d_{n+1} > 2
\end{equ}
and, for all $1\leq k\leq n$, uniformly over $\square\in\RECT$,
\begin{equ}[eq:d_k_2D]
| \widehat \SD(\square)^{[k]} | \lesssim
(\bfa^{d_k -\eps d_1}\bfb^{\eps d_1})^\zeta\;,
\end{equ}
where we recall that $\bfa,\bfb$ are the side lengths of $\square$ with $\bfa\geq \bfb$.


Then there exists a unique continuous $2$-cocycle $\SD \colon \Rect\to H$ such that $\pi_n\SD = \widehat\SD$ and which satisfies
\begin{equ}[eq:weak_bound_ext]
|\SD(\square)^{[k]}| \lesssim \bfa^{d_k \zeta}
\end{equ}
for all $k\leq m$.
Moreover, $\SD$ satisfies the stronger bound \eqref{eq:d_k_2D} for all $k\leq m$, i.e.~$| \SD(\square)^{[k]} | \lesssim
(\bfa^{d_k -\eps d_1}\bfb^{\eps d_1})^\zeta$.

Furthermore, the `extension map'
\begin{equ}[eq:ext_map]
\CC(\Paths,G)\times\CC(\Rect,H^{[n]}) \ni (\PD,\widehat\SD) \mapsto \SD \in \CC(\Rect,H)
\end{equ}
is uniformly continuous over $\BB_1\times\BB_2$
where $\BB_1$ is any set such that \eqref{eq:PD_bound} holds uniformly over $\PD\in\BB_1$ and
$\BB_2$ is any set such that \eqref{eq:d_k_2D} holds uniformly over $\widehat\SD\in\BB_2$.
\end{lemma}

\begin{remark}
Since $\bfa^{d_k-\eps d_1}\bfb^{\eps d_1}$ is decreasing in $\eps$,
the bound \eqref{eq:d_k_2D} becomes weaker as $\eps\to 0$.
Therefore, both the assumption and conclusion of \Cref{lem:extension_2D} weaken as $\eps\to0$. 
\end{remark}

\begin{proof}
By induction, it suffices to consider the case $n=m-1$.
Extend $\widehat\SD$ to a map $\widehat\SD\in\CC(\Rect,H)$ by $\widehat\SD(\square)^{[m]} = 0$ for all $\square$.
Since $\widehat\SD$ satisfies Chen on level $m-1$, we have
$\widehat\SD (\square_A\square_B)^{[i]} = ( \actionGG_{\PD(\rightarrow)} (\widehat\SD{\square_B})  \widehat\SD\square_A)^{[i]}$
for all $1\leq i\leq m-1$,
where $\actionGG_{\PD(\rightarrow)}$ is understood as the action of $\PD(\rightarrow)\in G$ on $\widehat\SD{\square_B}\in H$.
On the other hand, $\widehat\SD ( \square_A\square_B)^{[m]}=0$ and
\begin{equ}[eq:mSD_1]
		(\actionGG_{\PD(\rightarrow)} (\widehat\SD{\square_B} ) \widehat\SD{\square_A})^{[m]}
		=
		(\actionGG_{\PD(\rightarrow)} \widehat\SD{\square_B} )^{[m]}
		+ \sum_{\alpha,\beta} C_{\alpha,\beta}^m
		(\actionGG_{\PD(\rightarrow)} \widehat\SD {\square_B} )^{[\alpha]}
		\widehat\SD (\square_A)^{[\beta]}
	\end{equ}
	for all balanced partitions $\square_A\square_B=\square$
	where $\alpha,\beta$ are non-zero with $\deg(\alpha)+\deg(\beta) = d_m$.

Furthermore, by \eqref{eq:act_hom}, for all $i\leq m$,
\begin{equ}[eq:mSD_2]
(\actionGG_{\PD(\rightarrow)} \widehat\SD {\square_B} )^{[i]} =
\widehat\SD (\square_B)^{[i]} + \sum_{\alpha,\beta} B^i_{\alpha,\beta} \PD(\rightarrow)^{[\alpha]} \widehat\SD (\square_B)^{[\beta]}
\end{equ}
where $\alpha,\beta$ are non-zero with $\deg(\alpha)+\deg(\beta) = d_i$.
Recall that $\widehat\SD (\square_B)^{[m]}=0$.

For $i=1,\ldots, m$, define the function $F_i\colon\Rect\to [0,\infty)$
\begin{align*}
F_i(\square)	=
\begin{cases}
(\bfa\bfb)^{\zeta d_i/2} \quad &\text{ if } E(\square) \leq 3
\\
\bfa^{\zeta (d_i - \eps d_1)}\bfb^{\eps \zeta d_1}\quad &\text{ if } E(\square) > 3\;.
\end{cases}
\end{align*}

Since $|\PD(\rightarrow)^{[i]}| \lesssim \bfa^{\zeta c_i}$ due to \eqref{eq:PD_bound}
and $|\widehat\SD (\square_A)^{[i]}|\lesssim F_i(\square_A),|\widehat\SD (\square_B)^{[i]}| \lesssim F_i(\square_B)$ due to~\eqref{eq:d_k_2D},
and $F_i(\square_j) \leq F_i(\square)$ for $j=1,2$,
it follows that
\begin{equ}[eq:Wm_bound]
|(\actionGG_{\PD(\rightarrow)} (\widehat\SD {\square_B})  \widehat\SD {\square_A})^{[m]}|
\lesssim
\begin{cases}
(\bfa\bfb)^{\zeta d_m/2} \quad &\text{ if } E(\square) \leq 3
\\
\bfa^{\zeta (d_m-\eps d_1)} \bfb^{\eps d_1} \quad &\text{ if } E(\square) > 3
\end{cases}
\end{equ}
where we used that, if $E(\square)\leq 3$, then $\bfa \asymp (\bfa\bfb)^{1/2}$,
and that, if $E(\square) > 3$, then the `worst' term in \eqref{eq:mSD_1}-\eqref{eq:mSD_2} is $\PD(\rightarrow)^{[\alpha]} (\widehat\SD {\square_B})^{[\beta]}$ with multi-index $\beta=(1,0,\ldots,0)\in \ZZ^m$
for which we have the bounds $(\widehat\SD{\square_B})^{[\beta]} \lesssim \bfa^{\zeta(1-\eps)d_1}\bfb^{\zeta \eps d_1}$
and $\PD(\rightarrow)^{[\alpha]}\lesssim \bfa^{\zeta (d_m - d_1)}$ due to \eqref{eq:act_hom}.

The right-hand side of~\eqref{eq:Wm_bound} is precisely $F_m(\square)$.
Since $\zeta d_m/2 > 1$ by assumption \eqref{eq:dn_gamma} and $\eps \leq \frac12$, it follows that $d_m-\eps d_1 \geq d_m-\eps d_m \geq d_m/2$,
thus
\begin{equ}
\zeta (d_m - \eps d_1) \geq \zeta d_m/2 > 1\;.
\end{equ}
Therefore $F_m$ is a subcontrol by \Cref{prop:subcontrol}.


To summarize, for the Euclidean norm $| \cdot |_\h$ on $H\simeq \h$,
\begin{align*}
	| \widehat\SD (\square_A\square_B) - 
		\actionGG_{\PD(\rightarrow)}( \widehat\SD\square_B )\, \widehat\SD \square_A
		|
		\lesssim  F_m(\square)\;.
\end{align*}
The same bound holds for vertical balanced partitions $\genfrac{}{}{0pt}{1}{\square_B}{\square_A} = \square$.
We therefore satisfy approximate Chen in the sense of \Cref{thm:2D_sewing} wherein $\rho$ is the Euclidean metric.

To verify approximate Stoke in the sense of \Cref{thm:2D_sewing}, note that, since $\widehat\SD$ satisfies Stokes to level $m-1$, we have for every square $\square$ with side length $\bfa$
\begin{equ}
|\widehat\SD \square - \PD(\partial\square)| \leq |\PD(\partial\square)|_{\g_{d_m}} \lesssim
\bfa^{\zeta d_m} = |\square|^{\zeta d_m/2}\;.
\end{equ}
Since $\zeta d_m/2 > 1$ by assumption \eqref{eq:dn_gamma}, the map $\square\mapsto |\square|^{\zeta d_m/2}$ is a subcontrol.

Therefore all the conditions of \Cref{thm:2D_sewing} are satisfied and we conclude that there exists a unique $2$-cocycle $\SD\in\CC (\Rect, H)$ such that
\begin{equ}
|\widehat\SD \square - \SD \square|_\h \lesssim F_m(\square)\;.
\end{equ}
Since $\widehat\SD$ is already a $2$-cocycle on level $n$, the uniqueness part of \Cref{thm:2D_sewing} implies that $\pi_n\SD = \widehat\SD$.
Furthermore since $\widehat\SD^{[m]}=0$, it follows that $|(\SD\square)^{[m]}|_\h \lesssim F_m(\square)$ as desired.
The fact that the bound~\eqref{eq:weak_bound_ext} is enough to determine $\SD$ uniquely follows from \Cref{rem:weaker_F}.

Finally, uniform continuity of the extension map \eqref{eq:ext_map} follows from that of the sewing map \eqref{eq:sewing_map}
once we remark that the bound~\eqref{eq:PD_bound} uniform over $\BB_1$ implies \eqref{eq:BB1}
and likewise for \eqref{eq:d_k_2D} and \eqref{eq:BB2} with $\BB_2$.
\end{proof}

\subsection{Signature of a rough surface}

We now focus on the crossed module
\begin{equ}
\freeH^N \coloneqq \freeH^N(\R^n) \stackrel{\feedbackFreeCmGG}\to 
\freeG^N \coloneqq \freeG^{N}(\R^n)
\stackrel{\actionGG}\to \Aut(\freeH^N(\R^n))\;,
\end{equ}
which is the exponentials of the level-$N$ truncation of the ``free crossed module of Lie algebras''
as in \Cref{ex:free_nilpotent}.
That is, $\freeG^N$ is the exponential of $\g^N \coloneqq \g^0_{n;\leq N}$ (i.e., the level-$N$ free nilpotent Lie group over $\R^n$), and $\freeH^N$ is the exponential of $\h^N \coloneqq \g_{n;\leq N}^{-1}$,
where we apply \Cref{prop:crossedmodule} to build a crossed module of nilpotent Lie groups from a crossed module of nilpotent Lie algebra.
The degrees of $\freeG^N$ are integers $1=c_1 \leq \cdots \leq c_\ell = N$
and those of $\freeH^N$ are integers $2=d_1 \leq\cdots \leq d_m = N$.
As usual, we identify $\freeG^N,\freeH^N$ with $\g^N,\h^N$ via the $\exp$-maps.

For $1\leq k\leq N$, and $x\in \g^N$, we let $x^{(k)}\in \g_k \coloneqq \g_{k}^N$ denote the projection of $x$ onto the $k$-th homogenous level $\g_k$ (so that $x=x^{(1)}+\cdots+x^{(N)}$ and $\delta_\lambda x^k = \lambda^k x^k$),
and similarly for $x\in \h^N$ and $2\leq k \leq N$.

The extension \Cref{lem:extension_2D} (with $\eps=\frac12$), yields the following natural definition.

\begin{definition}[Rough surface]
\label{def:RS}
Consider $\zeta\in(0,1]$ and $N\geq 1$ such that $\zeta(N+1) > 2$.
A \DEF{level-$N$ $\zeta$-H\"older rough surface} is a pair of maps $\RS = (\PD,\SD)\in \CC(\Paths,\freeG^N)\times \CC(\Rect,\freeH^N)$, 
locally uniformly continuous,
such that $\PD$ is a $1$-cocycle and $\SD$ is a $2$-cocycle and there exists $C>0$ such that
\begin{equ}
|\PD(\gamma)^{(k)}| \leq C|\gamma|_{BV}^{\zeta k}
\end{equ}
for all $\gamma\in\Paths$ with $|\gamma|<1$ and $k=1,\ldots, N$, and
\begin{equ}
|\SD(\square)^{(k)}| \leq C (\bfa^{k-1}\bfb)^\zeta
\end{equ}
for all $\square\in\RECT$ and $k=2,\ldots, N$,
and where we recall that $\bfa,\bfb$ are the side length of $\square$ with $\bfa\geq\bfb$.
The smallest constant $C$ for which these bounds hold is denoted $\|\RS\|_\zeta$ and is called the \DEF{$\zeta$-H\"older norm} of $\RS$.
\end{definition}

\begin{remark}
Although $\PD$ in \Cref{def:RS} is, a priori, an arbitrary $1$-cocycle, the facts that $(\feedbackFreeCmGG h)^{(1)}=0$ for all $h\in \freeH^N$ and that there exists a $2$-cocycle $\SD$ above $\PD$ implies that $\PD^{(1)}(\gamma)=0$ whenever $\gamma$ is a loop.
It follows that there exists a (unique up to translation) function $X\colon [0,1]^2\to \R^n$ such that $\PD(\gamma)^{(1)}=X(y)-X(x)$ for any path $\gamma$ starting at $x$ and ending at $y$.
This $X$ is the `surface' aluded to in the above definition.
\end{remark}

We can now state the extension theorem for rough surfaces, which is a corollary of the Extension \Cref{lem:extension_1D,lem:extension_2D}.
For $M\geq N$, as earlier let $\pi_{\leq N}\colon \freeG^M(\R^n)\to \freeG^N(\R^n)$ and $\pi_{\leq N}\colon H^M \to H^N$ denote the canonical projections (the fact that both maps are denoted $\pi_{\leq N}$ will not cause confusion).

\begin{theorem}[Extension of rough surfaces]\label{thm:RS_ext}
Suppose $\RS = (\PD,\SD)$ is a level-$N$ $\zeta$-H\"older rough surface.
For every $M\geq N$, there exists a unique level-$M$ $\zeta$-H\"older rough surface $\YY$ such that $\pi_{\leq N}\YY = \RS$.
Furthermore, the map $\RS\mapsto\YY$ is uniformly continuous in the uniform metric on every set of the form $\{\RS\,:\,\|\RS\|_\zeta \leq C\}$ for $C>0$.
\end{theorem}

\begin{proof}
The facts that there exists a unique level-$M$ $\zeta$-H\"older rough `$1$-cocycle' $\bar \PD$ such that $\pi_{\leq N}\bar\PD = \PD$ and that the map $\PD\mapsto \bar\PD$ is locally uniformly continuous follows from \Cref{lem:extension_1D}.
(In \Cref{lem:extension_1D}, continuity is stated for the $\zeta$-H\"older norm, from which continuity under the uniform metric on sets bounded in $\zeta$-H\"older norm follows by interpolation, see \cite[Lem.~8.16]{FV10_RPs}.)
In turn,
the existence of the extension of $\SD$ (and its continuity) to a level-$M$ $2$-cocycle
is given by \Cref{lem:extension_2D}.
\end{proof}

\begin{remark}
In \cite{lee2024_2Dsigs} another notion of rough surface and extension theorem is given,
and a decay estimate as $M\to\infty$ is shown that is similar to the original extension theorem of \cite{Lyons98} for rough paths.
\end{remark}

\subsection{Young case}
\label{subsec:Young}

In this section, we give several conditions under which a `surface' $X\colon [0,1]^2 \to \R^n$ admits a canonical (and/or unique) lift to a rough surface in the sense of \Cref{def:RS}.
We first give a condition based on $\zeta$-H\"older norms for $X$ to admit a unique lift for $\zeta > \frac23$ (and demonstrate non-uniqueness for $\zeta=\frac23$),
and then give a simple condition based on `rectangular' H\"older norms for $X$ to admit a canonical (but not unique) lift.

Consider a $X\colon [0,1]^2 \to \R^n$.
Assume that, for some $\zeta>\frac23$,
\begin{equ}[eq:X_Hol]
\|X\|_{\Hol\zeta} \coloneqq \sup_{x\neq y} \frac{|X(y)-X(x)|}{|x-y|^{\zeta}} < \infty\;.
\end{equ}
For any path $\gamma \colon [0,1]\to [0,1]^2$,
it holds that $\|X\circ\gamma\|_{\Hol\zeta} \leq \|X\|_{\Hol\zeta}|\gamma|_{\Hol1}^\zeta$ and thus, if $\gamma$ is parametrized at constant speed, $\|X\circ\gamma\|_{\Hol\zeta} \leq \|X\|_{\Hol\zeta}|\gamma|_{BV}^\zeta$.
It follows therefore from \Cref{lem:extension_1D} (or equivalently, from classical Young integration) that there exists a unique continuous $1$-cocycle $\PD\colon\Paths \to \freeG^2$ such that
$\PD(\gamma)^{(1)} = X(\gamma(1)) - X(\gamma(0))$ and,
for $k=1,2$,
\begin{equ}
|\PD(\gamma)^{(k)}| \lesssim |\gamma|_{BV}^{k\zeta}\;.
\end{equ}

Since the kernel of $\feedbackFreeCmGG$ is trivial on level $2$, there exists a unique continous $2$-cocycle $\SD\colon \Rect\to \freeH^2$ such that
\begin{equ}
\feedbackFreeCmGG \SD(\square) = \PD(\partial \square)\;.
\end{equ}
(The fact that $\SD$ is indeed a $2$-cocycle over $\PD$ is a consequence of the classical Chen identity.)

We now \emph{assume} the extra regularity condition
\begin{equ}[eq:P_Hol]
|\PD (\partial \square)^{(2)}| \lesssim |\square|^\zeta = \bfa^\zeta\bfb^\zeta\;.
\end{equ}
Since $3\zeta>2$, $\RS \coloneqq (\PD,\SD)\in\CC(\Paths,\freeG^2)\times\CC(\Rect,\freeH^2)$ is a level-$2$ $\zeta$-H\"older rough surface.
\Cref{thm:RS_ext} now implies the following result. 

\begin{proposition}
Assume the two regularity assumptions \eqref{eq:X_Hol} and \eqref{eq:P_Hol} for some $\zeta\in(\frac23,1]$ and let notation be as above.
Then $\RS$
admits a unique extension to all levels $N>2$ that satisfies the bound
\begin{equ}[eq:SD_Young]
|\SD(\square)^{(k)}| \lesssim \bfa^{\zeta k}\;,
\end{equ}
and this extension furthermore satisfies the stronger bound
\begin{equ}[eq:SD_Young_strong]
|\SD(\square)^{(k)}| \lesssim \bfa^{\zeta(k-1)}\bfb^\zeta\;.
\end{equ}
\end{proposition}

\begin{example}[Optimality of $\zeta>\frac23$ --- pure level-$3$ rough surface]
It may appear strange that $\zeta>\frac23$ is the threshold for the `Young regime' to extend $\SD$ uniquely.
This is because $\zeta>\frac12$ is the threshold for the Young regime to extend $\PD$ uniquely.

The following example demonstrates a `pure level-$3$' $\frac23$-H\"older rough surface and thus the necessity of $\zeta>\frac23$ to obtain unique lifts in our setting.
This is analogous to the `pure area' $\frac12$-H\"older lift of the constant path in rough path theory.

Consider the constant surface $X\equiv 0$.
The canonical level-$N$ extension of $X$, given by the above construction, is simply $\PD\equiv e_{\freeG^N}$ and $\SD \equiv e_{\freeH^N}$,
and this is the unique lift which satisfies \eqref{eq:SD_Young} for all $\zeta\in(\frac23,1]$.

We now demonstrate a \emph{different} extension of $X$ to level-$3$ which satisfies
for $\zeta = \frac23$
\begin{equ}[eq:SD_k_23]
|\SD(\square)^{{(k)}}| \lesssim |\square|^{\zeta k / 2}\;.
\end{equ}
(which is even stronger than \eqref{eq:SD_Young_strong}).
Let $\PD \equiv e_{\freeG^3}$ as before (here we have no choice but to use the canonical Young lift),
but now define $\SD\colon \Rect\to \freeH^3$ by
\[
\SD(\square)^{(2)} = 0\;,\qquad \SD(\square)^{(3)} = \omega \int_\square Y\;,
\]
where $\omega\in \h^3$ is any non-zero element such that $\feedbackFreeCmAA\omega=0$ and $Y \in L^\infty ([0,1]^2,\R)$.
Then $\SD$ is a $2$-cocyle above $\PD$ by additivity of integration.
Furthermore
\begin{equ}
|\SD(\square)^{(3)}| \leq |Y|_{\infty} |\square|
\end{equ}
showing that \eqref{eq:SD_k_23} is satisfied for $k=3$ (the case $k=2$ is trivial).
Therefore $\RS = (\PD,\SD)$ is a level-$3$ $\frac23$-H\"older rough surface lift of $X$.

(Applying \Cref{thm:RS_ext} to $\RS$, there is a unique level-$N$ $\zeta$-H\"older rough surface lift of $\RS$ for all $N\geq 3$ and $\zeta\in (\frac12,\frac23]$,
where the lower bound $\zeta>\frac12$ comes from the condition $(N+1)\zeta>2$ from \Cref{def:RS} with $N=3$.)
\end{example}

\subsubsection{Rectangular increment regularity}

We now show that one can still obtain a `canonical' level-$3$ rough surface lift for $\zeta>\frac12$ provided that we control the rectangular increments of $X$.
\textit{En passant}, we provide a sufficient (linear) condition on $X$ for \eqref{eq:P_Hol} to hold.

Assume that, now for $\zeta\in(\frac12,1]$, $X$ satisfies the H\"older bound~\eqref{eq:X_Hol} and
\begin{equ}[eq:X_Hol_rect_inc]
\|X\|_{\Hol\zeta;\RECT} \coloneqq \sup_{\square\in\RECT} \frac{|\delta_\square X|}{|\square|^{\zeta} } < \infty\;,
\end{equ}
where, for $\square = (e^-,e^+)$ with lower left corner $e^- = (e^-_1,e^-_2)$ and top right corner $e^+ = (e^+_1,e^+_2)$,
\begin{equ}
\delta_\square X = X(e^+) - X(e^{-}_1, e^+_2) - X(e^+_1,e^-_2) + X(e^-)
\end{equ}
is the rectangular increment of $X$.
As before, the $1$-cocycle $\PD\colon \Paths\to \freeG^3$ is well-defined since $\zeta>\frac12$.

We claim that \eqref{eq:P_Hol} holds, or more precisely, that
\begin{equ}[eq:P_Hol_rect]
|\PD(\partial\square)^{(2)} | \lesssim |\square|^\zeta \|X\|_{\Hol\zeta} ( \|X\|_{\Hol\zeta} + \|X\|_{\Hol\zeta;\Rect})\;.
\end{equ}
Indeed, for two parallel lines $\gamma,\bar\gamma$ forming the long sides of $\square$, we have
\begin{equ}
\|X\circ\gamma\|_{\Hol\zeta} \vee \|X\circ\bar\gamma\|_{\Hol\zeta} \leq \bfa^\zeta \|X\|_{\Hol\zeta}\;,
\end{equ}
and likewise for the short sides.
Furthermore
\begin{equ}
\|X\circ\gamma - X\circ\bar\gamma\|_{\Hol\zeta} \leq \bfb^\zeta \|X\|_{\Hol\zeta;\Rect}\;.
\end{equ}
Together with Chen's identity, these bounds readily imply \eqref{eq:P_Hol_rect}.

Recall from \Cref{ex:general_case} that the level-$3$ lift of $X$ takes the form
\begin{equ}
\SD(\square)^{(3)} \coloneqq \Omega^{(3)} = \sum_i \sum_{j<k} \int_\square dq dr\ X^{(i)}_{r,q} J^{(jk)}_{r,q} [\Z{i}, \Z{jk}]\;.
\end{equ}
The identity \eqref{eq:level2_stokes}
(extrapolated to general rectangles) and the bound \eqref{eq:P_Hol_rect} imply that $|\int_\square J^{(jk)}_{r,q}|\lesssim |\square|^\zeta$,
which, by 2D Young integration (see \cite[Thm.~16]{Harang_21_sewing} or \cite[Thm.~6.18]{FV10_RPs}),
implies that $\SD(\square)^{(3)} = \Omega^{(3)}$ is well-defined and satisfies
\begin{equ}
|\SD(\square)^{(3)} | \lesssim \bfa^{2\zeta}\bfb^{\zeta}\|X\|_{\Hol\zeta} ( \|X\|_{\Hol\zeta} + \|X\|_{\Hol\zeta;\Rect})^2 \;,
\end{equ}
where we supposed, without loss of generality, that $X(e^-)=0$.
This is precisely the bound required for $\RS = (\PD,\SD)$ to be a level-$3$ $\zeta$-H\"older rough surface.
This demonstrates a canonical (but not unique) level-$3$ lift of $X$ for $\zeta>\frac12$
under the H\"older regularity conditions \eqref{eq:X_Hol} and \eqref{eq:X_Hol_rect_inc}.


\appendix

\section{Some proofs of the main text}
\label{ape:proofs}

\begin{proof}[{Proof of \Cref{prop:quotient}}]
  Let us define \(\feedbackAA'(a+I)=\feedbackAA(a)+J\) and \(\actionAA'_{X+J}(a+I) = \actionAA_X(a) + I\).
  Note first that \(\feedbackAA'\) is well defined, since if \(a-a'\in I\) then property 2 implies that \(\feedbackAA(a)-\feedbackAA(a')\in J\), so that
  \begin{align*}
    \feedbackAA'(a+I) &= \feedbackAA(a)+J\\
    &= \feedbackAA(a') + J \\
    &= \feedbackAA'(a'+I).
  \end{align*}
  For a similar reason \(\feedbackAA'\) is a Lie algebra morphism. Indeed, if \(a,a'\in\mathfrak{h}\) then
  \begin{align*}
    \feedbackAA'([a+I,a'+I]) &= \feedbackAA'([a,a']+I) \\
                             &= \feedbackAA([a,a'])+J \\
                             &= [\feedbackAA(a),\feedbackAA(a')] +J\\
                             &= [\feedbackAA(a)+J,\feedbackAA(a')+J] \\
                             &= [\feedbackAA'(a+I),\feedbackAA'(a'+I)].
  \end{align*}

  Let us now verify that \(\actionAA'\) is well defined.
  Suppose that \(X,X'\in\mathfrak{g}\) and \(a,a'\in\mathfrak{h}\) are such that \(X'-X\in J\) and \(a'-a\in I\).
  Then
  \begin{align*}
    \actionAA'_{X'+J}(a'+I) &= \actionAA_{X'}(a') + I \\
                            &= \actionAA_X(a') + \actionAA_{X'-X}(a')  + I \\
                            &= \actionAA_X(a) + \actionAA_X(a'-a) + \actionAA_{X'-X}(a') + I
  \end{align*}
  However, since by property 4 we have that \(\actionAA_{X'-X}(a')\in I\) and by property 3 that \(\actionAA_X(a'-a)\in I\), we get that
  \begin{align*}
    \actionAA'_{X'+J}(a'+I) &= \actionAA_X(a) + I \\
                            &= \actionAA'_{X+J}(a+I).
  \end{align*}
  Moreover, for every \(X\in\mathfrak{g}\), \(\actionAA'_{X+J}\in \operatorname{Der}(\mathfrak{h}/I)\).
  Indeed, if \(a,a'\in\mathfrak{h}\) then
  \begin{align*}
    \actionAA'_{X+J}([a + I, a' + I]) &= \actionAA_X([a,a']) + I \\
                                      &= [\actionAA_X(a),a'] + [a,\actionAA_X(a')] + I \\
                                      &= [\actionAA_X(a) + I, a' + I] + [a + I,\actionAA_X(a') + I] \\
                                      &= [\actionAA'_{X+J}(a+I),a'+I] + [a+I,\actionAA'_{X+J}(a'+I)].
  \end{align*}
  Finally, let us verify that \(\actionAA'\colon\mathfrak{g}/J\to \operatorname{Der}(\mathfrak{h}/I)\) is a Lie morphism.
  Take \(X,X'\in\mathfrak{g}\) and \(a\in\mathfrak{h}\), then
  \begin{align*}
    \actionAA'_{[X+J,X'+J]}(a + I) &= \actionAA'_{[X,X']+J}(a+I) \\
                                   &= \actionAA_{[X,X']}(a) + I \\
                                   &= \actionAA_X(\actionAA_{X'}(a)) - \actionAA_{X'}(\actionAA_X(a)) + I \\
                                   &= \actionAA'_{X+J}(\actionAA'_{X'_J}(a)) - \actionAA'_{X'+J}(\actionAA_{X+J}(a)).
  \end{align*}

  Next, we verify that \(\feedbackAA'\) and \(\actionAA'\) satisfy \cref{eq:equiDiff,eq:peifferDiff}.
  Let \(X,X'\in\mathfrak{g}\) and \(a,a'\in\mathfrak{h}\). Then
  \begin{align*}
    \feedbackAA'(\actionAA'_{X+J}(a+I)) &= \feedbackAA'(\actionAA_X(a) + I) \\
    &= \feedbackAA(\actionAA_X(a))+J \\
    &= [X,\feedbackAA(a)] + J \\
    &= [X+J,\feedbackAA(a)+J] \\
    &= [X+J,\feedbackAA'(a + I)]
  \end{align*}
  so that \cref{eq:equiDiff} holds.
  Similarly,
  \begin{align*}
    \actionAA'_{\feedbackAA'(a+I)}(a'+I) &= \actionAA'_{\feedbackAA(a)+J}(a'+I) \\
    &= \actionAA_{\feedbackAA(a)}(a')+I \\
    &= [a,a']+I \\
    &= [a+I,a'+I].\qedhere
  \end{align*}
\end{proof}

\begin{proof}[{Proof of \Cref{lem:exp}}]
  The Leibniz rule yields
  \[
    \delta^n[X,Y] = \sum_{k=0}^n\binom{n}{k}[\delta^kX,\delta^{n-k}Y].
  \]
  Thus
  \begin{align*}
    \exp(\delta)[X,Y] &= \sum_{n=0}^{\infty}\frac{1}{n!}\delta^n[X,Y] \\
                      &= \sum_{n=0}^{\infty}\frac{1}{n!}\sum_{k=0}^n\binom{n}{k}[\delta^kX,\delta^{n-k}Y] \\
                      &= \sum_{k=0}^{\infty}\sum_{n=k}^\infty\frac{1}{k!}\frac{1}{(n-k)!}[\delta^kX,\delta^{n-k}Y] \\
                      &= [\exp(\delta)X,\exp(\delta)Y].
  \end{align*}
  The last statement is evident.
\end{proof}

\begin{proof}[{Proof of \Cref{prop:crossedmodule}}]
  We begin by noting that by the previous Lemma, indeed we have \(\actionGG_g\in \operatorname{Aut}(H)\).
  Next, we check that \(\actionGG\) is a morphism of groups.
  Since \(\actionAA\) is a morphism of Lie algebras we have
  \begin{align*}
    \actionGG_{gg'} &= \exp(\actionAA_{\operatorname{BCH}_{\mathfrak{g}}(g,g')}) \\
                    &= \exp\Bigl(\operatorname{BCH}_{\operatorname{Der}(\mathfrak{h})}(\actionAA_g,\actionAA_{g'})\Bigr) \\
                    &= \exp(\actionAA_{g})\circ\exp(\actionAA_{g'}).
  \end{align*}
  That \(\mathfrak{T}\) is a morphism of Lie groups is direct from its definition.

  Now we verify that \(\actionGG\) and \(\feedbackGG\) satisfy \cref{eq:equivariant,eq:peiffer}.
  For \cref{eq:equivariant}, it suffices to note that \cref{eq:equiDiff} implies that for any \(k\ge 1\),
  \[
    \feedbackAA(\actionAA^k_x(y'))=\ad_x^k(\feedbackAA(y')),
  \]
  so that
  \begin{align*}
    \feedbackGG(\actionGG_g(h')) &= \feedbackAA(\exp(\actionAA_g)(h')) \\
    &= \exp(\ad_g)(\feedbackAA(h)) \\
    &= g\feedbackGG(h')g^{-1}
  \end{align*}
  since \(\exp(\ad_x)=\Ad_{\exp(x)}\) (and the Lie exponential map is the identity in our case).
  \Cref{eq:peiffer} follows similarly, by simply noting that \cref{eq:peifferDiff} implies that for any \(k\ge 1\),
  \[
    \actionAA^k_{\feedbackAA(y)}(y')=\ad_y^k(y'),
  \]
  hence
  \begin{align*}
    \actionGG_{\feedbackGG(h)}(h') &= \exp(\actionAA_{\feedbackAA(h)})(h') \\
    &= \exp(\ad_h)(h') \\
    &= hh'h^{-1}.\qedhere
  \end{align*}
\end{proof}


\section{Coordinates of the first kind}
\label{s:appendix_coordinates}

\newcommand\w[1]{\mathtt{#1}} 


When expressing the logarithm of the iterated-integrals signature
in terms of a basis of the free Lie algebra,
the dual element which picks out the coefficient
of a basis element in the signature is called the \emph{coordinate of the first kind}
(see \cite[Section 3]{diehl2020areas} for a recent overview of this classical topic).

For the sanity checks in \Cref{ex:general_case} we need
the coordinates of the first kind on level $3$ and $4$,
for the basis used in \Cref{ex:dimensions} (for $n=3$).
It is given in \Cref{tab:in_basis}.
\begin{table}
\centering
\newcommand{\specialcell}[2][c]{\begin{tabular}[#1]{@{}l@{}}#2\end{tabular}}
\begin{tabular}{ll}
	Lie basis & Coefficient (``coordinates of the first kind'') \\
\hline
$[\Z1,[\Z1,\Z2]]$ & $\frac16 ( \w{112} - 2 \ \w{121} + \w{211})$ \\
$[\Z1,[\Z1,\Z3]]$ & $\frac16 (\w{113} - 2 \ \w{131} + \w{311})$ \\
$[\Z2,[\Z1,\Z2]]$ & $\frac16  (-\w{122} + 2\ \w{212} - \w{221})$ \\
$[\Z2,[\Z1,\Z3]]$ & $\frac16( \w{123} - 2 \  \w{132} + \w{213} - 2 \  \w{231} + \w{312} + \w{321} )$ \\
$[\Z2,[\Z2,\Z3]]$ & $\frac16 ( \w{223} - 2 \  \w{232} + \w{322} )$\\
$[\Z3,[\Z1,\Z2]]$ & $\frac16 ( -2 \ \w{123} + \w{132} + \w{213} + \w{231} +  \w{312} - 2 \  \w{321} )$ \\
$[\Z3,[\Z1,\Z3]]$ & $\frac16 ( -\w{133} + 2\ \w{313} - \w{331} )$\\
$[\Z3,[\Z2,\Z3]]$ & $\frac16 ( -\w{233} + 2\ \w{323} - \w{332} )$\\
	$[\Z1, [\Z1, [\Z1, \Z2]]]$ & $-\tfrac{1}{6} \w{1121} +\tfrac{1}{6} \w{1211}$ \\
	$[\Z2, [\Z1, [\Z1, \Z2]]]$ & $\tfrac{1}{6} \w{2211} -\tfrac{1}{6} \w{2121} +\tfrac{1}{6} \w{1212} -\tfrac{1}{6} \w{1122}$ \\
	$[\Z2, [\Z2, [\Z1, \Z2]]]$ & $-\tfrac{1}{6} \w{2122} +\tfrac{1}{6} \w{2212}$ \\
	$[\Z3, [\Z1, [\Z1, \Z2]]]$ & $-\tfrac{1}{6} \w{1123} +\tfrac{1}{6} \w{1213} +\tfrac{1}{6} \w{1231} -\tfrac{1}{6} \w{1321} -\tfrac{1}{6} \w{3121} +\tfrac{1}{6} \w{3211}$ \\
	$[\Z3, [\Z2, [\Z1, \Z2]]]$ & $\tfrac{1}{6} \w{3212} -\tfrac{1}{6} \w{2123} -\tfrac{1}{6} \w{3122} -\tfrac{1}{6} \w{1322} +\tfrac{1}{6} \w{1232} -\tfrac{1}{6} \w{2321} +\tfrac{1}{6} \w{2231} +\tfrac{1}{6} \w{2213}$ \\
	$[\Z3, [\Z3, [\Z1, \Z2]]]$ & $\tfrac{1}{6} \w{1233} -\tfrac{1}{6} \w{1323} -\tfrac{1}{6} \w{3123} +\tfrac{1}{6} \w{3213} +\tfrac{1}{6} \w{3231} -\tfrac{1}{6} \w{3321}$ \\
	$[\Z1, [\Z1, [\Z1, \Z3]]]$ & $-\tfrac{1}{6} \w{1131} +\tfrac{1}{6} \w{1311}$ \\
	$[\Z2, [\Z1, [\Z1, \Z3]]]$ & $-\tfrac{1}{6} \w{1132} +\tfrac{1}{6} \w{1312} +\tfrac{1}{6} \w{1321} -\tfrac{1}{6} \w{1231} -\tfrac{1}{6} \w{2131} +\tfrac{1}{6} \w{2311}$ \\
	$[\Z2, [\Z2, [\Z1, \Z3]]]$ & $\tfrac{1}{6} \w{1322} -\tfrac{1}{6} \w{1232} -\tfrac{1}{6} \w{2132} +\tfrac{1}{6} \w{2312} +\tfrac{1}{6} \w{2321} -\tfrac{1}{6} \w{2231}$ \\
	$[\Z3, [\Z1, [\Z1, \Z3]]]$ & $\tfrac{1}{6} \w{3311} -\tfrac{1}{6} \w{3131} +\tfrac{1}{6} \w{1313} -\tfrac{1}{6} \w{1133}$ \\
	$[\Z3, [\Z2, [\Z1, \Z3]]]$ & $\tfrac{1}{6} \w{1323} -\tfrac{1}{6} \w{3231} -\tfrac{1}{6} \w{1233} -\tfrac{1}{6} \w{2133} +\tfrac{1}{6} \w{2313} -\tfrac{1}{6} \w{3132} +\tfrac{1}{6} \w{3312} +\tfrac{1}{6} \w{3321}$ \\
	$[\Z3, [\Z3, [\Z1, \Z3]]]$ & $-\tfrac{1}{6} \w{3133} +\tfrac{1}{6} \w{3313}$ \\
	$[\Z2, [\Z2, [\Z2, \Z3]]]$ & $-\tfrac{1}{6} \w{2232} +\tfrac{1}{6} \w{2322}$ \\
	$[\Z3, [\Z2, [\Z2, \Z3]]]$ & $\tfrac{1}{6} \w{3322} -\tfrac{1}{6} \w{3232} +\tfrac{1}{6} \w{2323} -\tfrac{1}{6} \w{2233}$ \\
	$[\Z3, [\Z3, [\Z2, \Z3]]]$ & $-\tfrac{1}{6} \w{3233} +\tfrac{1}{6} \w{3323}$ \\
	$[[\Z1, \Z2], [\Z1, \Z3]]$ & $-\tfrac{1}{3} \w{1231} +\tfrac{1}{6} \w{2311} +\tfrac{1}{3} \w{1321} +\tfrac{1}{6} \w{3112} +\tfrac{1}{6} \w{1123} -\tfrac{1}{6} \w{3211} -\tfrac{1}{6} \w{2113} -\tfrac{1}{6} \w{1132}$ \\
	$[[\Z1, \Z2], [\Z2, \Z3]]$ & $-\tfrac{1}{3} \w{1232} +\tfrac{1}{3} \w{2321} +\tfrac{1}{6} \w{1322} +\tfrac{1}{6} \w{3122} +\tfrac{1}{6} \w{1223} -\tfrac{1}{6} \w{3221} -\tfrac{1}{6} \w{2213} -\tfrac{1}{6} \w{2231}$ \\
	$[[\Z1, \Z3], [\Z2, \Z3]]$ & $-\tfrac{1}{6} \w{1332} +\tfrac{1}{6} \w{3132} +\tfrac{1}{6} \w{1323} -\tfrac{1}{6} \w{3231} -\tfrac{1}{6} \w{2313} +\tfrac{1}{6} \w{2331}$ \\
\end{tabular}
\caption{Magnus in a basis (i.e.~coordinates of the first)}
\label{tab:in_basis}
\end{table}



\section{Why crossed modules?}
\label{ss:why}

In this section, we provide an intuitive description of the conceptual transition from 0-dimensional (0D) objects (points) to 1-dimensional (1D) objects (paths), to motivate the subsequent transition from 1D to 2-dimensional (2D) objects (surfaces). We will place particular emphasis on how the structure of crossed module is motivated by the latter transition.

This section aims to inspire, avoiding intricate details or strict mathematical rigor (such as, for instance,~questions of regularity).


\subsection{Paths}

The basic 1-dimensional (1D) objects that we aim to study are paths $X^{(1)} \colon [0,1]\to\R^n$, the set of which, only in this subsection, we denote $\CO^{(1)}$.
In this subsection, we use the superscript $(d)$ to emphasise the parameter dimension $d$ of the associated objects.
This space comes equipped with a natural notion of splitting / concatenation.

Given paths $X^{(1)}, X_A^{(1)}, X_B^{(1)} \in\CO^{(1)}$, we write $X^{(1)} = X_A^{(1)}\sqcup X_B^{(1)}$ if there exist $t_0\in [0,1]$ and increasing bijections
\begin{equ}[eq:psis]
	\psi_A\colon [0,1] \to [0,t_0]\;,\qquad \psi_B\colon [0,1] \to [t_0,1]
\end{equ}
such that for all $t\in[0,1]$,
\begin{equ}
	X^{(1)} (\psi_A(t)) = X_A^{(1)}(t)\;,
	\qquad 
	X^{(1)} (\psi_B(t)) = X_B^{(1)}(t)\;.
\end{equ}
\textbf{NB.} The equality $X^{(1)} = X_A^{(1)}\sqcup X_B^{(1)}$ should not be interpreted as an equality in the strict sense since we never define $X_A^{(1)}\sqcup X_B^{(1)}$ as an element of $\CO^{(1)}$.

Consider now a group $G$. We aim to study the path space $\CO^{(1)}$ via maps $M^{(1)}\colon \CO^{(1)}\to G$, which are \emph{multiplicative} in the sense that
\begin{equ}
	M^{(1)}(X^{(1)}) = M^{(1)}(X_A^{(1)})M^{(1)}(X_B^{(1)})
\end{equ}
whenever $X^{(1)} = X_A^{(1)}\sqcup X_B^{(1)}$. One of the motivations to consider such maps is computational: if we view $M^{(1)}$ as a feature map on $\CO^{(1)}$, then the multiplicativity property
implies that it suffices to compute $M^{(1)}$ on the `sub-path' over $[0,t_0]$, for any $t_0\in[0,1]$, then compute $M^{(1)}$ on the `sub-path' over $[t_0,1]$, and then combine the results directly in the group $G$.
This can be much more efficient than computing $M^{(1)}$ directly over the full path on $[0,1]$.

The following is a simple and natural way to build multiplicative maps.
Consider a linear map $m^{(1)}\colon \R^n\to \g$ where $\g$ is a vector space.
Then the map
\begin{equ}
m^{(1)}(X^{(1)} ) \coloneqq m^{(1)}(X^{(1)}(1)) - m^{(1)}(X^{(1)} (0))
= \int_0^1 m^{(1)}(\dot X^{(1)}(t)) \, dt
\end{equ}
is multiplicative for the abelian group structure of $\g$.

It is important to remark that $m^{(1)}$ factors through a boundary operator on $\CO^{(1)}$.
More precisely, there is a natural boundary operator $\partial^{(1)} \colon \CO^{(1)} \to \CO^{(0)}$,
where $\CO^{(0)}$ is the space of all maps $X^{(0)} \colon \{0,1\}\to\R^n$, defined by
\begin{equ}
	(\partial^{(1)} X^{(1)} ) = X^{(1)} \circ\gamma^{(0)} \colon \{0,1\} \to \R^n\;,
\end{equ}
where $\gamma^{(0)} \colon \{0,1\} \to [0,1]$ is the canonical injection, i.e.~$(\partial^{(1)} X^{(1)} )(0) = X^{(1)} (0)$ and $(\partial^{(1)} X^{(1)} )(1)= X^{(1)} (1)$.
Then we define $m^{(0)} \colon \CO^{(0)}\to\g$ by $m^{(0)}(X^{(0)} ) = m^{(1)}(X^{(0)} (1))-m^{(1)}( X^{(0)}(0))$
and remark that, trivially, for $X^{(1)} \in\CO^{(1)}$,
\begin{equ}
m^{(1)}(X^{(1)}) = m^{(0)}(\partial^{(1)} X^{(1)})\;.
\end{equ}


Now suppose $\g$ is the Lie algebra of a Lie group $G$ (not necessarily abelian).
Equip $G$ with a suitable (e.g.~Riemannian) metric $\rho$.
Then $m^{(1)}$ induces a unique map $M^{(1)}\colon \CO^{(1)} \to G$ which is multiplicative
and satisfies the local increment property\footnote{Here \(\lesssim\) means up to a constant factor.}
\begin{equ}[eq:LI1]
	\rho(M^{(1)}(X^{(1)} ), e^{m^{(1)}(X^{(1)} )}) \lesssim \|X^{(1)} \|_{BV}^\theta
\end{equ}
for some $\theta>1$, where $BV$ stands for 1D bounded variation.
It is given as the time-$1$ solution of the ODE
\begin{equ}
	dY(t) = Y(t) m^{(1)} (dX^{(1)}(t))\;,\qquad Y(0) = 1 \in G\;.
\end{equ}

We call any multiplicative map $M^{(1)}$ constructed as above \textbf{linearly generated}.
If $G$ is abelian, then $M^{(1)}(X^{(1)} )=e^{m^{(1)}(X^{(1)} )} = e^{m^{(0)}(\partial^{(1)} X^{(1)} )}$, so that again $M^{(1)}$ factors through the boundary.

There is a universal linearly generated map $S^{(1)}\colon \CO^{(1)}\to \freeG((\R^n))$, called the signature, that every other linearly generated map factors through.
It is given by applying the above construction to the linear map denoted $s^{(1)} \colon\R^n\to \R^n$ (instead of $m^{(1)}$),
which we take to be the identity and where $\R^n$ on the right-hand side is identified canonically with a subspace of $\FL((\R^n))$, the Lie algebra corresponding to the group-like elements $\freeG((\R^n))$.
That is, $S^{(1)}(X^{(1)} )$ is the time-$1$ solution of the ODE
\begin{equ}
	dY(t) = Y(t) \otimes s^{(1)}(dX^{(1)}(t)) \;,\qquad Y(0)=1
\end{equ}
taking values in the algebra of tensor series.
Then, for every linear $m^{(1)}\colon \R^n\to\g$ there is a unique group morphism $\widetilde M^{(1)}\colon \freeG((\R^n))\to G$ such that
\begin{equ}
	M^{(1)}(X^{(1)} ) = \widetilde M^{(1)}\circ S^{(1)}(X^{(1)} )\;.
\end{equ}
Note that $\freeG((\R^n))$ is not a bona fide finite-dimensional Lie group, but a suitable subgroup of it can be turned into a Polish group with the stated universality property, see reference \cite{CL16_Char} for details.

It is known that the signature of $X^{(1)}$ determines $X^{(1)}$ up to tree-like (i.e.~thin-homotopy) equivalence \cite{Chen58_faithful,HL10, BGLY2016}.



\subsection{Surfaces}

We now turn to our 2D objects of interest, which are maps $X^{(2)} \colon [0,1]^2 \to \R^n$, the set of which, again only in this subsection, we denote $\CO^{(2)}$.

We first remark that, as in the 1D case, there is a boundary operator $\partial^{(2)} \colon \CO^{(2)} \to \CO^{(1)}$ defined by
\begin{equ}
	(\partial^{(2)} X^{(2)}) = X^{(2)}\circ\gamma^{(1)} \colon [0,1] \to \R^n\;,
\end{equ}
where $\gamma^{(1)} \colon [0,1] \to [0,1]^2$ is a fixed loop with $\gamma^{(1)}(0)=\gamma^{(1)}(1)=0$ which runs around the boundary of $[0,1]^2$, by convention, counter-clockwise.

Moreover, again as in the 1D case, $\CO^{(2)}$ comes equipped with a natural notion of splitting / concatenation, but the description of this is now more involved.

\begin{definition}\label{def:splitting}
Given $X^{(2)},X_A^{(2)},X_B^{(2)}\in\CO^{(2)}$, we write $X^{(2)} = X_A^{(2)}\sqcup X_B^{(2)}$ if
\begin{itemize}
\item there exists a point $x\in \partial [0,1]^2$ and a path $\zeta \colon [0,1] \to [0,1]^2$ such that $\zeta(0)=0$ and $\zeta(1)=x$,

\item there exist maps $\psi_A,\psi_B\colon [0,1]^2\to [0,1]^2$ with $\psi_i(0)=0$ and such that the boundary of $\partial \psi_A$ is given by first running counter-clockwise along $\partial [0,1]^2$ until we hit $x$, and then running backwards along $\zeta$, while $\partial \psi_B$ is given by first running along $\zeta$ and then running from $x$ counter-clockwise along $\partial [0,1]^2$ until we return to $0$.

\item one has for all $z\in [0,1]^2$
\begin{equ}
	X^{(2)}(\psi_A(z)) = X_A^{(2)}(z)\;,\quad X^{(2)}(\psi_B(z)) = X_B^{(2)}(z)\;.
\end{equ}
\end{itemize}
\end{definition}

See \Cref{fig:splitting} for illustrations.
Remark that we do not require $\psi_i$ to be injective: in the right drawing in \Cref{fig:splitting}, a triangle with non-zero area is collapsed to a line.
This example especially motivates our definition of $2$-cocycle in \Cref{ss:cocycles_lie_group}.

\begin{figure}
\centering
\begin{tikzpicture}[scale=1.65]
    \draw[thick] (0,1.5) rectangle (1,2.5);
    
    \draw[thick, blue, -{Latex[length=2mm]}] (0,1.5) to[out=55, in=200] node[midway, above] {$\zeta$} (1,2.3);
	\draw[thick, blue,dashed, {Latex[length=2mm]}-] (0.1,1.5) to[out=55, in=200] (1,2.2);
    
    \node[below left] at (0.1,1.5) {$0$};
    
	\draw[thick, {Latex[length=2mm]}-, purple] (-1.5,-0.5) --++ (0,1) --++ (1,0) --++ (0,-1);
	\draw[thick, -{Latex[length=2mm]}, blue] (-1.5,-0.5) --++ (1,0);

	\draw[thick, {Latex[length=2mm]}-, orange] (1.5,0.5) -- (2.5,0.5) -- (2.5,-0.5) -- (1.5,-0.5);
	\draw[thick, {Latex[length=2mm]}-, dashed, blue] (1.5,-0.5) --++ (0,1);
    
    \coordinate (center_left) at (0.2, 2.0); 
    \coordinate (center_right) at (0.8, 1.85); 
    
    \draw[-{Latex[length=2mm]}, thick, bend right=20] (-1,0.5) to[out=45, in=150] node[midway, left] {$\psi_B$} (center_left);
    \draw[-{Latex[length=2mm]}, thick, bend left=20] (2,0.5) to[out=-45, in=210]  node[midway, right] {$\psi_A$} (center_right);
    
\end{tikzpicture}
\quad\quad
\begin{tikzpicture}[scale=1.65]
    \draw[thick] (0,1.5) rectangle (1,2.5);
    
	\colorlet{darkgreen}{green!70!black}
    \draw[thick, darkgreen, -{Latex[length=2mm]}] (0,1.5) to (0,2);
	\draw[thick, blue, -{Latex[length=2mm]}] (0,2) to node[midway, above] {$\zeta$} (1,2);
	\draw[thick, blue, dashed, {Latex[length=2mm]}-] (0,1.95) to (1,1.95);
	\draw[thick, darkgreen, dashed, {Latex[length=2mm]}-] (0.05,1.5) to (0.05,1.95);
    
    \node[below left] at (0.1,1.5) {$0$};
    
	\draw[thick, {Latex[length=2mm]}-, purple] (-1.5,0.5) --++ (1,0) --++ (0,-1);
	\draw[thick, {Latex[length=2mm]}-, blue] (-0.5,0.5) --++ (0,-1);
	\draw[thick, {Latex[length=2mm]}-, darkgreen, dashed] (-1.5,-0.5) --++ (0,1);
	\draw[thick, -{Latex[length=2mm]}, darkgreen] (-1.5,-0.5) --++ (1,0);
	\draw[thick, brown, dashed] (-1.5,0.5) --++ (1,-1);

	\draw[thick, {Latex[length=2mm]}-, orange] (2.5,0.5) -- (2.5,-0.5) -- (1.5,-0.5);
	\draw[thick, {Latex[length=2mm]}-, dashed, blue] (1.5,0.5) -- (2.5,0.5);
	\draw[thick, {Latex[length=2mm]}-, dashed, darkgreen] (1.5,-0.5) --++ (0,1);
    
    \coordinate (center_left) at (0.2, 2.3); 
    \coordinate (center_right) at (0.6, 1.7); 
    
    \draw[-{Latex[length=2mm]}, thick, bend right=20] (-1,0.5) to[out=45, in=150] node[midway, left] {$\psi_B$} (center_left);
    \draw[-{Latex[length=2mm]}, thick, bend left=20] (2,0.5) to[out=-45, in=210]  node[midway, right] {$\psi_A$} (center_right);
    
\end{tikzpicture}

\caption{Examples of surface splitting.
In the right figure, the mapping $\psi_B$ collapses the bottom left triangle of the left square to the green line on the boundary of the top square;
in particular, the entire brown diagonal line in the left square is mapped to the intersection point of the green and blue lines in the top square.
}
\label{fig:splitting}
\end{figure}

The path $\zeta$ should be interpreted as a `splitting' or `cut' of the surface and plays the role of the point $t_0$ in \eqref{eq:psis} from the 1D setting.
Again, the notation $X^{(2)}=X_A^{(2)}\sqcup X_B^{(2)}$ should not be interpreted as a strict equality.

We remark that the boundary map is multiplicative in the following sense: if $X^{(2)}\in\CO^{(2)}$ with $X^{(2)}=X_A^{(2)}\sqcup X_B^{(2)}$, then
\begin{equ}[eq:Chen12]
	M^{(1)} (\partial^{(2)} X^{(2)}) = M^{(1)}( \partial^{(2)} X_A^{(2)}) M^{(1)} (\partial^{(2)} X_B^{(2)})
\end{equ}
for any linearly generated multiplicative map $M^{(1)}\colon\CO^{(1)}\to G$, and in particular for $M^{(1)}=S^{(1)}$ the (universal linearly generated) path signature map.

As in the 1D case, we wish to study maps $M^{(2)}\colon \CO^{(2)}\to H$, where $H$ is a group, that are multiplicative in the sense that 
\begin{equ}
	M^{(2)} (X^{(2)}) = M^{(2)}(X_A^{(2)}) M^{(2)}(X_B^{(2)})
\end{equ}
whenever $X^{(2)} = X_A^{(2)}\sqcup X_B^{(2)}$.
The computational motivation for such maps is the same as before: we can compute $M^{(2)}$ over subdomains first, and then combine the answers in the group $H$.

A natural way to build multiplicative maps is the following.
Consider a $2$-form $m^{(2)}\colon \R^n\wedge\R^n \to \h$, for a vector space $\h$.
Then the map
\begin{equ}
	m^{(2)}(X^{(2)}) \coloneqq \int_{[0,1]^2} (X^{(2)})^* m^{(2)} = \int_{[0,1]^2} m^{(2)}(\partial_1 X^{(2)}\wedge \partial_2 X^{(2)}) \, dt_1 \, dt_2 
\end{equ}
is multiplicative for the abelian group $\h$.

Now suppose $\h$ is the Lie algebra of a Lie group $H$.
Equip $H$ with a suitable (e.g. Riemannian) metric $\rho$.
Then, as in 1D, we aim to find a multiplicative map $M^{(2)}\colon \CO^{(2)} \to H$ which satisfies the local increment property
\begin{equ}[eq:LI2]
	\rho(M^{(2)}(X^{(2)}), e^{m^{(2)}(X^{(2)})}) \lesssim |X^{(2)}|^\theta
\end{equ}
for some $\theta>1$ and $|X| \coloneqq \|X\|_{BV}$, where $BV$ now stands for 2D bounded variation.
In analogy with the 1D case, 
such maps $M^{(2)}$, whenever they exist, will be called \textbf{linearly generated}.

An important case where we can find such a map $M^{(2)}$ is the so-called \emph{Schlesinger case}.
We take $\h$ to be the Lie algebra of a Lie group $H$ and $m^{(2)}(x , y) = [m^{(1)}(x),m^{(1)}(y)]_\h$ for $x,y\in\R^n$ and some linear $m^{(1)}\colon \R^n \to \h$.
Then $M^{(2)}(X^{(2)}) \coloneqq M^{(1)}(\partial^{(2)} X^{(2)})$, with $M^{(1)}\colon\CO^{(1)}\to H$ on the right determined as in \eqref{eq:LI1},
is indeed multiplicative; it is the unique map that satisfies \eqref{eq:LI2}.

To see this, observe that $2(X^{(2)})^* m^{(2)} = 2m^{(2)}(\partial_1 X^{(2)}, \partial_2 X^{(2)})\diff{t}^1 \wedge \diff{t}^2 $
is the exterior derivative of the $1$-form\footnote{To make sense of $\alpha$ as a $1$-form, we recall that $X^{(2)}$ is $\R^n$-valued.}
\begin{equ}
	\alpha = m^{(2)}(X^{(2)},\partial_1 X^{(2)})\diff{t}^1 + m^{(2)}(X^{(2)},\partial_2 X^{(2)}) \diff{t}^2\;.
\end{equ}
So, by Stokes' theorem,
\begin{equ}[eq:Levy_area]
	\int_{[0,1]^2}(X^{(2)})^* m^{(2)} = \frac12\int_{ [0,1]^2}\diff{\alpha} = \frac12\int_{\partial [0,1]^2} \alpha
  = \int_0^1 [m^{(1)} ( Y(t) ),m^{(1)} (\dot Y(t))]_{\mathfrak{h}}\diff{t}\;,
\end{equ}
where $Y(t) = X^{(2)}\circ \gamma^{(1)}(t)$ and where the final equality follows from a direct
computation\footnote{We note that \eqref{eq:Levy_area} can be identified with the so-called L{\'e}vy area of $X^{(2)}\circ\gamma^{(1)}$ contracted with $m^{(1)}\otimes m^{(1)}$.}  using $Y(1)=Y(0)=0$.
The term \eqref{eq:Levy_area} is the
first non-zero approximation of $M^{(1)}(\partial^{(2)} X^{(2)})$ so we indeed have, denoting $e_H$ the identity element of $H$,
\begin{equ}
	M^{(2)}(X^{(2)}) = e_H  + \int_{[0,1]^2}(X^{(2)})^* m^{(2)} + O(\|X\|_{BV}^{3/2}) = e^{m^{(2)}(X^{(2)})} +  O(\|X\|_{BV}^{3/2}) \;.
\end{equ}
where we implicitly treat $H$ and $\h$ as embedded in a matrix algebra to make sense of the final two expressions.

Any map built in the Schlesinger case will, by construction, factor through the boundary and thus not reveal any `surface' information.

In the rest of this subsection, we show how to build a class of linearly generated multiplicative maps that do reveal surface information.
A crucial step turns out to be a generalisation of an algebraic identity in the Schlesinger case (but which appears to have no analogue in 1D)
which naturally leads to a crossed module.

To describe this identity, note that $M^{(2)}(X^{(2)})$, in the Schlesinger case, is sensitive to the location of the basepoint $X^{(2)}(0)\in\R^n$.
That is, if
\begin{equ}
	X^{(2)}= \tilde X^{(2)}\circ \psi
\end{equ}
where $\psi\colon [0,1]^2\to[0,1]^2$ is a diffeomorphism such that $\psi(x)=0$ for some $x \in \partial [0,1]^2$ and which preserves the boundary and its orientation, then in general $M^{(2)}(X^{(2)}) \neq M^{(2)}(\tilde X^{(2)})$.
However, there is a \textbf{basepoint shift} identity given by
\begin{equ}[eq:base_point_shift]
	\Ad_{M^{(1)}(\tilde X^{(2)} \circ\eta)} M^{(2)}(\tilde X^{(2)}) = M^{(2)}(X^{(2)})
\end{equ}
where $\eta\colon [0,1]\to\partial [0,1]^2$ is the curve from $x$ to $0$ along $\partial [0,1]^2$ moving clockwise.

We now suitably generalise this identity.
Suppose we are given two groups, $G$ and $H$, with a group homomorphism $\triangleright \colon G \to \Aut(H)$,  $\triangleright\colon g \mapsto \triangleright_g$,
and multiplicative maps $M^{(1)}\colon \CO^{(1)}\to G$
and $M^{(2)}\colon \CO^{(2)}\to H$
that are linearly generated.

The map $M^{(1)}$ encodes boundary information of $X^{(2)}$ via the multiplicative map $X^{(2)}\mapsto M^{(1)}(X^{(2)}\circ \gamma^{(1)})$,
while $M^{(2)}$ encodes the `surface' information.
(As far as surface information is concerned, there is no loss of generality in assuming $G\subset \Aut(H)$ and $\triangleright$ is the inclusion map.)
The Schlesinger case would correspond to $G = H$ and $\triangleright_g = \Ad_g$.

We now ask that $M^{(1)}$ and $M^{(2)}$ satisfy the following generalisation of \eqref{eq:base_point_shift}:
\begin{equ}
\label{eq:bsc}
	\triangleright_{M^{(1)}(\tilde X^{(2)}\circ \eta)} M^{(2)}(\tilde X^{(2)}) = M^{(2)}(X^{(2)})\;.
\end{equ}
We now explain how this condition naturally leads to the structure of a crossed module.

Note that if $X^{(2)} = X^{(2)}_A\sqcup X^{(2)}_B$ with corresponding maps $\zeta,\psi_A,\psi_B$, then we can find a (non-unique) $Y^{(2)}$ such that $X^{(2)} = Y^{(2)}\sqcup X^{(2)}_A$ with corresponding maps $\tilde\zeta,\tilde\psi_A,\tilde\psi_B$ where $\tilde\psi_B=\psi_A$.
For example, as in \Cref{fig:splitting2},
take $\tilde\zeta$ as $\partial\psi_A$ run up until some point $\tilde x \in \partial \psi_B \cap \partial [0,1]^2$ (we could simply take $\tilde x=0$ in which case we run around the entire boundary $\partial \psi_A$).
This also determines $\partial \tilde\psi_A$ and $\partial \tilde\psi_B$ (the latter, by construction, agrees with $\partial \psi_A$).
From this, is it easy to build $\tilde \psi_A$ with the claimed property.

\begin{figure}
\centering
\begin{tikzpicture}[scale=1.65]
	\draw[thick, orange, -{Latex[length=2mm]}] (0,1.5) to (1,1.5) to (1,2);
	\draw[thick, purple, -{Latex[length=2mm]}] (1,2) to (1,2.5) to (0,2.5) to (0,2);
    
	\colorlet{darkgreen}{green!70!black}
    \draw[thick, darkgreen, -{Latex[length=2mm]}] (0,1.5) to (0,2);
	\draw[thick, blue, -{Latex[length=2mm]}] (0,2) to node[midway, above] {$\zeta$} (1,2);
	\draw[thick, blue, dashed, {Latex[length=2mm]}-] (0,1.95) to (1,1.95);
	\draw[thick, darkgreen, dashed, {Latex[length=2mm]}-] (0.05,1.5) to (0.05,1.95);
    
    \node[below left] at (0.1,1.5) {$0$};
    
	\draw[thick, {Latex[length=2mm]}-, purple] (-1.5,0.5) --++ (1,0) --++ (0,-1);
	\draw[thick, {Latex[length=2mm]}-, blue] (-0.5,0.5) --++ (0,-1);
	\draw[thick, {Latex[length=2mm]}-, darkgreen, dashed] (-1.5,-0.5) --++ (0,1);
	\draw[thick, -{Latex[length=2mm]}, darkgreen] (-1.5,-0.5) --++ (1,0);
	\draw[thick, brown, dashed] (-1.5,0.5) --++ (1,-1);

	\draw[thick, {Latex[length=2mm]}-, orange] (2.5,0.5) -- (2.5,-0.5) -- (1.5,-0.5);
	\draw[thick, {Latex[length=2mm]}-, dashed, blue] (1.5,0.5) -- (2.5,0.5);
	\draw[thick, {Latex[length=2mm]}-, dashed, darkgreen] (1.5,-0.5) --++ (0,1);
    
    \coordinate (center_left) at (0.2, 2.3); 
    \coordinate (center_right) at (0.6, 1.7); 
    
    \draw[-{Latex[length=2mm]}, thick, bend right=20] (-1,0.5) to[out=45, in=150] node[midway, left] {$\psi_B$} (center_left);
    \draw[-{Latex[length=2mm]}, thick, bend left=20] (2,0.5) to[out=-45, in=210]  node[midway, right] {$\psi_A$} (center_right);
    
\end{tikzpicture}
\quad\quad
\begin{tikzpicture}[scale=1.65]
	\draw[thick, orange, -{Latex[length=2mm]}] (0,2) to (0,1.5);
	\draw[thick, purple, -{Latex[length=2mm]}] (1,2) to (1,2.5) to (0,2.5) to (0,2);
    
	\colorlet{darkgreen}{green!70!black}
    \draw[thick, blue, -{Latex[length=2mm]}] (0,1.5) to (1,1.5) to (1,2);
	\draw[thick, blue, -{Latex[length=2mm]}] (1,2) to node[midway, above] {$\tilde \zeta$} (0,2);
	\draw[thick, blue, dashed, -{Latex[length=2mm]}] (0,1.95) to (0.95,1.95);
	\draw[thick, blue, dashed] (0.95,1.95) to (0.95,1.55);
	\draw[thick, blue, dashed, -{Latex[length=2mm]}] (0.95,1.55) to (0,1.55);
    
    \node[below left] at (0.1,1.5) {$0$};
	\node[below left] at (0.05,2.05) {$\tilde x$};
    
	\draw[thick, {Latex[length=2mm]}-, blue, dashed] (-1.5,0.5) --++ (1,0);
	\draw[thick, purple, {Latex[length=2mm]}-] (-0.5,0.5) --++ (0,-1);
	\draw[thick, {Latex[length=2mm]}-, blue, dashed] (-1.5,-0.5) --++ (0,1);
	\draw[thick, -{Latex[length=2mm]}, blue] (-1.5,-0.5) --++ (1,0);
	\draw[thick, brown, dashed] (-1.5,0.5) --++ (1,-1);

	\draw[thick, {Latex[length=2mm]}-, blue] (1.5,0.5) -- (2.5,0.5) ;
	\draw[thick, {Latex[length=2mm]}-, blue] (2.5,0.5) -- (2.5,-0.5) -- (1.5,-0.5);
	\draw[thick, {Latex[length=2mm]}-, orange] (1.5,-0.5) --++ (0,1);
    
    \coordinate (center_left) at (0.2, 2.3); 
    \coordinate (center_right) at (0.6, 1.7); 
    
    \draw[-{Latex[length=2mm]}, thick, bend right=20] (-1,0.5) to[out=45, in=150] node[midway, left] {$\tilde \psi_A$} (center_left);
    \draw[-{Latex[length=2mm]}, thick, bend left=20] (2,0.5) to[out=-45, in=210]  node[midway, right] {$\tilde\psi_B = \psi_A$} (center_right);
    
\end{tikzpicture}

\caption{
Two ways of subdividing the same square.
}
\label{fig:splitting2}
\end{figure}

Since $M^{(2)}$ is linearly generated, it is simple to see that $M^{(2)}(Z^{(2)})=\unit_H$
whenever $Z^{(2)}$ has zero area (like the maps sending the lower left triangles in \Cref{fig:splitting2} to curves).
In particular, taking above $\tilde x=0$
and suitably shifting the basepoint of $Y^{(2)}$, we obtain from \eqref{eq:bsc} that
\begin{equ}
\triangleright_{M^{(1)}(X^{(2)}_A\circ \gamma^{(1)})}
M^{(2)}(X^{(2)}_B) = M^{(2)}(Y^{(2)})\;.
\end{equ}
On the other hand, by multiplicativity, $M^{(2)}(X^{(2)}_A)M^{(2)}(X^{(2)}_B) = M^{(2)}(Y^{(2)})M^{(2)}(X^{(2)}_A)$.
It follows that
\begin{equ}[eq:pre_stokes]
	\triangleright_{M^{(1)}(X^{(2)}\circ \gamma^{(1)})}  M^{(2)}(X^{(2)}_B) = \Ad_{M^{(2)}(X^{(2)})} M^{(2)}(X^{(2)}_B)\;.
\end{equ}

It is now natural to assume that there exists a group morphism $\feedbackGG \colon H \to G$ that satisfies
\begin{equ}[eq:Peiffer1]
\triangleright_{\feedbackGG(h)} = \Ad_h
\end{equ}
and
\begin{equ}[eq:Peiffer2]
\feedbackGG\circ \triangleright_g =\Ad_g\circ  \feedbackGG\;.
\end{equ}
One should think of $\feedbackGG$ as the abstract analogue of $\Ad\colon H \to \Aut(H)$.
(Again, there is no loss of generality in taking $G\subset \Aut(H)$ and $\feedbackGG\colon h\mapsto \Ad_h$, but it is crucial that $G$ is larger than $\feedbackGG(H)$ to capture non-boundary information as otherwise we are in the Schlesinger case.)
The tuple $(G,H,\feedbackGG,\triangleright)$ is the structure of a \textbf{crossed module}.

Due to \eqref{eq:Peiffer1}, property \eqref{eq:pre_stokes} reads simply as
\begin{equ}
\feedbackGG(M^{(2)}X^{(2)}) = M^{(1)}(X^{(2)}\circ \gamma^{(1)})\;.
\end{equ}
This last condition is called (non-commutative) \textbf{Stokes' identity} and will play an important role in the sequel.
In light of \eqref{eq:bsc},
condition \eqref{eq:Peiffer1} follows naturally from the demand that $M^{(2)}(X^{(2)})=M^{(2)}(X^{(2)}_A)M^{(2)}(X^{(2)}_B)$
for \textit{all} $X^{(2)}_A, X^{(2)}_B$ such that $X^{(2)}=X^{(2)}_A\sqcup X^{(2)}_B$,
while \eqref{eq:Peiffer2} follows from the Stokes' identity.

It turns out that, starting with a crossed module $(G,H,\feedbackGG,\triangleright)$ and linear maps $m^{(2)}\colon \R^n\wedge\R^n \to\h$ and $m^{(1)}\colon \R^n\to\g$ that are suitably compatible with the differentials of $\feedbackGG$ and $\triangleright$ (see \Cref{thm:2cocycleIntegration,thm:2cocycle_sewing} for details),
we can build a unique multiplicative map $M^{(2)}\colon\CO^{(2)}\to H$ that satisfies the local increment property \eqref{eq:LI2} (and,
additionally, \eqref{eq:bsc}).
The existence and characterization of this map $M^{(2)}$ via the germ $e^{m^{(2)}(X^{(2)})}$,
as well as the study of the universal linearly generated map, i.e. the `surface signature', is one of the main contributions of this paper.

\section{Kapranov's proofs}
\label{app:kapranov}

The representation theory of $\GL(\R^n)$  will play a
central role for both of these tasks, so we briefly recall it
(see for example \cite{fulton1997young,fulton2013representation,sam2012introduction} for more background).
%
%
The finite-dimensional irreducible representations (\DEF{irreps}) of $\GL(\R^n)$ 
are indexed by partitions $\lambda$ of integers $p \in \N$,
that is $\lambda_1 \ge \lambda_2 \ge \dots \ge \lambda_k > 0$ with $\sum_{i=1}^k \lambda_i = p$.
A basis for the irrep corresponding to the \DEF{shape} $\lambda$ is indexed
by \DEF{semistandard} fillings (weakly increasing along rows
and strictly increasing along columns), from $[n]$,
of the Young diagram of $\lambda$.
For example, for $\lambda = (2,1)$ and $n=3$,
the semistandard fillings are
\begin{align*}
	\begin{ytableau}
		1 & 1 \\
		2
	\end{ytableau},
	\begin{ytableau}
		1 & 2 \\
		2
	\end{ytableau},
	\begin{ytableau}
		1 & 3 \\
		2
	\end{ytableau},
	\begin{ytableau}
		1 & 1 \\
		3
	\end{ytableau},
	\begin{ytableau}
		1 & 2 \\
		3
	\end{ytableau},
	\begin{ytableau}
		1 & 3 \\
		3
	\end{ytableau},
	\begin{ytableau}
		2 & 2 \\
		3
	\end{ytableau},
	\begin{ytableau}
		2 & 3 \\
		3
	\end{ytableau},
\end{align*}
so that particular irrep has dimension $8$.

Note that $\I_{n,p}$ and $\J_{n,p}$
index semistandard Young tableaux of shapes $(p-1,1)$ and $(p-1,1,1)$, respectively.

%

We start with well-known results, encapsulated in the following lemma.

\newcommand\dRham{d_{\scalebox{0.6}{$\Omega$}}}
\begin{lemma}
	\label{lem:cochain}
	Consider the cochain complex of polynomial forms on $V=\R^n$,
	\begin{align*}
		\Omega^\bullet(V) \coloneqq \bigoplus_{m \ge 0} \Omega^m(V), \qquad
		\Omega^m(V)
		\coloneqq  \Lambda^m( V^* ) \otimes S(V^*), 
	\end{align*}
	where $S(V^*)$ is the space of polynomials on $V$.
	We use the usual exterior/de Rham differential $\dRham$ on $\Omega^\bullet(V)$.
	As a vector space we declare the space $\Lambda^m(V^*) \otimes S^d(V^*)$
	of $m$-forms with polynomial coefficients of homogeneity $d$ to have degree $d+m$.
	With respect to the (shifted polynomial) grading, consider the graded dual
	\begin{align*}
		\Gamma_m(V) = \Lambda^m( V ) \otimes S(V),
	\end{align*}
	with analogous grading convention.
	Dualize $\dRham$ to a chain map $\partial$ on $\Gamma_\bullet$,
	$\partial_m \coloneqq (\dRham^{m-1})^*$.
	For $m \ge 1$ consider
	\begin{align*}
		\Gamma_m^\cl \coloneqq \ker \{ \partial_m\colon \Gamma_m \to \Gamma_{m-1} \}.
	\end{align*}
	Then (\Cref{fig:cochain})
	\begin{enumerate}[label=\roman*.]

	\item 
		$\Gamma_m^\cl \cong \coker \{\partial_{m+2}\colon \Gamma_{m+2} \to \Gamma_{m+1}\}$

	\item
	\begin{align*}
		\Gamma_m^\cl \cong (\ker \dRham^{m+1})^0,
	\end{align*}
	the graded dual (monomial grading).

	\item
	\label{it:gamma_is_irrep}
	%
	As $\GL(\R^n)$-representations (recall the ``shifted'' grading on $\Gamma_m$),
	\begin{align*}
		(\Gamma_m^\cl)_\ell \cong \text{ the irrep of $\GL(\R^n)$ corresponding to the shape $(\ell-m,\underbrace{1,..,1}_{m})$}.
	\end{align*}

	\end{enumerate}
\end{lemma}

\begin{figure}
	\centering
	\begin{tikzcd}
		& \Gamma_m^{\operatorname{cl}}=\ker \partial_m \arrow[d, hook] \arrow[r, "\cong" description, no head, phantom] & \operatorname{coker} \partial_{m+2} \arrow[r, "\cong" description, no head, phantom]                  & (\ker d^{m+1})^0 \arrow[lddd, "\text{\tiny{(graded) dual}}" description, no head, dotted, bend left=100] \\
	\Gamma_{m-1}                      & \Gamma_m \arrow[l, "\partial_m"']                                                                    & \Gamma_{m+1} \arrow[l, "\partial_{m+1}"'] \arrow[u, two heads] & \Gamma_{m+2} \arrow[l, "\partial_{m+2}"']                    \\
	\Omega^{m-1} \arrow[r, "d^{m-1}"] & \Omega^m \arrow[r, "d^m"]                                                                            & \Omega^{m+1} \arrow[r, "d^{m+1}"]                              & \Omega^{m+2}                                                 \\
		&                                                                                                      & \ker d^{m+1} \arrow[u, hook]                                   &                                                             
	\end{tikzcd}
	\vspace{-2em}
	\caption{$\Gamma^\cl_m$ etc.}
	\label{fig:cochain}
\end{figure}

%
%
%
%
%
%

Before proving \Cref{lem:cochain}, we introduce notation.
On $\Omega^m$ consider the basis
\begin{align*}
	x_1^{\alpha_1} \dots x_n^{\alpha_n} \diff{x}^{i_1} \wedge \dots \wedge \diff{x}^{i_m}.
\end{align*}
Write the dual basis (recall, we are working in the \emph{graded} dual) as
\begin{align*}
	y_1^{\alpha_1} \dots y_n^{\alpha_n} \diff{y}^{i_1} \wedge \dots \wedge \diff{y}^{i_m}.
\end{align*}
Then
\begin{align*}
	\partial_1 \left[ y_1^{\alpha_1} \dots y_n^{\alpha_n} \diff{y}^j \right] &= (\alpha_j+1) y_1^{\alpha_1} \dots y_j^{\alpha_j+1} \dots y_n^{\alpha_n} \\
	\partial_2 \left[ y^\alpha \diff{y}^j \wedge \diff{y}^k \right] &= (\alpha_j+1) y^{\alpha+e_j} \diff{y}^k - (\alpha_k+1) y^{\alpha+e_k} \diff{y}^j,
\end{align*}

We also consider on $\Omega^m$ the basis
\begin{align*}
	r_1^{\alpha_1} \dots r_n^{\alpha_n} \diff{r}^{i_1} \wedge \dots \wedge \diff{r}^{i_m}
	\coloneqq
	\frac1{\alpha_1!\dots \alpha_n!} x_1^{\alpha_1} \dots x_n^{\alpha_n} \diff{x}^{i_1} \wedge \dots \wedge \diff{x}^{i_m},
\end{align*}
and write the dual basis as
\begin{align*}
	s_1^{\alpha_1} \dots s_n^{\alpha_n} \diff{s}^{i_1} \wedge \dots \wedge \diff{s}^{i_m}.
\end{align*}
In this basis
\begin{align*}
	\dRham r^\beta = \sum_{\ell=1}^n r^{\beta-e_\ell} \diff{r}^\ell, \qquad\qquad
	\dRham [ r^\beta \diff{r}^i ] = \sum_{\ell=1}^n r^{\beta-e_\ell} \diff{r}^\ell \wedge \diff{r}^i.
\end{align*}

and
\begin{align*}
	\partial_1 \left[ s^\alpha \diff{s}^{j} \right] = s^{\alpha+e_j}, \qquad
	\partial_2 s^\alpha \diff{s}^j \wedge \diff{s}^k = s^{\alpha + e_j} \diff{s}^k - s^{\alpha + e_k} \diff{s}^j.
\end{align*}

\begin{proof}[Proof of \Cref{lem:cochain}]~

	\begin{enumerate}

	\item 
	Since $d$ is exact so is $\partial$ (except at $0$).
	Then
	\begin{align*}
		\Gamma_m^\cl = \ker \partial_m = \im \partial_{m+1} \cong \Gamma_{m+1} / \ker \partial_{m+1} = \Gamma_{m+1} / \im \partial_{m+2} = \coker \partial_{m+2}.
	\end{align*}

	\item 
	We establish
	\begin{align*}
		\coker \partial_{m+2} \cong (\ker \dRham^{m+1})^0.
	\end{align*}
	This is straightforward linear algebra.
	For $[f] \in \coker \partial_{m+2}$ define $T([f]) \in (\ker \dRham^{m+1})^0$
	via
	\begin{align*}
		T( [f] )( \omega ) \coloneqq f( \omega ).
	\end{align*}
	This is well-defined: if $f = \partial_{m+2} h \in \im \partial_{m+2}$ then
	for any $\omega \in \ker \dRham^{m+1}$:
	\begin{align*}
		f(\omega) = (\partial_{m+2} h)(\omega) = h( \dRham^{m+1} \omega ) = h(0) = 0,
	\end{align*}
	so we get a well-defined (linear) map $T\colon  \coker \partial_{m+2} \to (\ker \dRham^{m+1})^0$.

	$T$ is injective:
	assume
	\begin{align*}
		0 = T([f])(\omega) = f(\omega) \qquad \forall \omega \in \ker \dRham^{m+1}.
	\end{align*}
	Then (Fredholm alternative), there is $h \in \Gamma_{m+2}$
	with $\partial_{m+2} h = f$ and then $[f] = 0$.

	$T$ is surjective:
	let $\nu\colon \ker \dRham^{m+1} \to \R$ be linear.
	Pick $S$, a subspace of $\Omega^{m+1}$ complementary to $\ker \dRham^{m+1}$.
	Then, there exists a unique linear map $\hat \mu\colon \Omega^{m+1} \to \R$
	with 
	\begin{align*}
		\hat \mu(\omega) &= \mu(\omega), \omega \in \ker \dRham^{m+1} \\
		\hat \mu(s) &= 0, s \in S.
	\end{align*}
	Then, for $\omega \in \ker \dRham^{m+1}$,
	\begin{align*}
		T( [\hat \mu] )(\omega) = \hat \mu(\omega) = \mu(\omega),
	\end{align*}
	which gives surjectivity.

	\item
	A proof of this statement, using highest weight vectors, is sketched in \cite[Lemma 2.6]{kapranov1992hyperdeterminants},
	\cite[Proposition 2.2]{gelfand1994discriminants}.
	
	The isomorphism is given by the map
	\begin{align*}
		\phi:
		\mathbb S^{(\ell-m,1,\dots,1)}(V)
		&\to (\Gamma_m^\cl)_\ell \\
		\begin{ytableau}
			\scalebox{0.5}{$i_{\ell-m}$} & \dots & \scalebox{0.5}{$i_1$} \\
			\dots \\
			\scalebox{0.5}{$i_\ell$}
		\end{ytableau}
		&\mapsto
		\partial_2 ( s_{i_1} \dots s_{i_{\ell-m-1}} \diff{s}^{i_{\ell-m}} \wedge \dots \wedge \diff{s}^{i_{\ell}} ).
	\end{align*}
	where $\mathbb S^\lambda(V)$ is the Schur module corresponding to the shape $\lambda$.
	The map is well-defined, since 
	\begin{align*}
		\partial_2 ( s_{i_1} \dots s_{i_{\ell-m-1}} \diff{s}^{i_{\ell-m}} \wedge \dots \wedge \diff{s}^{i_{\ell}} )
		=
		0,
	\end{align*}
	if any of the $i_{\ell-m}, \dots, i_{\ell}$ are equal,
	and since for any $r = 1, \dots, \ell-m-1$
	\begin{align*}
		&\partial_2 ( s_{i_1} \dots s_{i_{\ell-m-1}} \diff{s}^{i_{\ell-m}} \wedge \dots \wedge \diff{s}^{i_{\ell}} )\\
		&\quad=
		\sum_{q=\ell-m,\dots,\ell}
		\partial_2 ( s_{i_1} \dots s_{i_{r-1}} s_{i_{q}} s_{i_{r+1}} \dots s_{i_{\ell-m-1}}    \diff{s}^{i_{\ell-m}} \wedge \dots \wedge \hbox{\sout{$\diff{s}^{i_{q}}$}} \diff{s}^{i_r} \wedge \dots \wedge \diff{s}^{i_{\ell}} ).
	\end{align*}

	\end{enumerate}

\end{proof}

\begin{theorem}
	\label{thm:reutenauers_map}
	For 
		$I = (i_1, i_2, \dots, i_p) \in [n]^p$,
	define the following element
	$\nu_I \in (\Omega^2)^0$,
	\begin{align*}
		\nu_I \coloneqq s^\alpha \diff{s}^{i_{p-1}} \wedge \diff{s}^{i_p}, \qquad \text{ where } \quad \alpha_\ell \coloneqq \#\{ k \le p-2 \mid i_k = \ell \},\ \ell=1,\dots,n.
	\end{align*}

	Then:
	\begin{enumerate}[label=\roman*.]

		\item
		Endow $\Gamma^\cl_1$ with the trivial Lie algebra structure (i.e.~with the zero bracket).
		Consider $\partial_2 \nu_I \in \Gamma^\cl_1$.
		Then: there exists a unique Lie algebra morphism
		$\rho_0\colon [\freecm{n}^0, \freecm{n}^0] \to \Gamma^\cl_1$ such that
		\begin{align*}
			\rho_0( [\Z{i_1}, [\Z{i_2}, \dots, [\Z{i_{p-1}}, \Z{i_p}]\dots]] ) = \partial_2 \nu_I, \quad \text{ $I=(i_1, \dots, i_p) \in \I_{n,p}$}.
		\end{align*}

		\item
		\label{it:any}
		$\rho_0$ induces a linear morphism
		\begin{align*}
			\rho^{\ab}_0\colon \frac{ [\freecm{n}^0, \freecm{n}^0] }{ \left[[\freecm{n}^0,\freecm{n}^0], [\freecm{n}^0,\freecm{n}^0]\right] } \to \Gamma^\cl_1,
		\end{align*}
		and for \emph{any} $J = (j_1,\dots,j_p) \in [n]^p$ we have
		\begin{align*}
			\rho_0( [\Z{j_1}, [\Z{j_2}, \dots, [\Z{j_{p-1}}, \Z{j_p}]\dots]] ) = \partial_2 \nu_J.
		\end{align*}

		\item
		\label{it:rho0iso}
		$\rho^\ab_0$ is in fact a linear isomorphism of $\GL(\R^n)$-representations,
		\begin{align}
			\label{eq:iso}
			\rho^{\ab}_0: \frac{ [\freecm{n}^0, \freecm{n}^0] }{ \left[[\freecm{n}^0,\freecm{n}^0], [\freecm{n}^0,\freecm{n}^0]\right] } \xrightarrow{\sim} \Gamma^\cl_1.
		\end{align}

	\end{enumerate}

\end{theorem}
\begin{figure}
	\centering
	\begin{tikzcd}
		\f^{-1} \arrow[d, "\sab"] \arrow[ddd, "\rho_{-1}"', bend right=80] \arrow[rr, "d"]    &  & {[\f^0,\f^0]} \arrow[d, "\id"] \arrow[ddd, "\rho_0", bend left=80] \\
		{\f^{-1}_{\sab} = \frac{\f^{-1}}{d[\f^{-1},\f^{-1}]}} \arrow[d, "\ab"] \arrow[rr, "d^{\sab}"]  \arrow[dd, "\rho^\sab_{-1}"', bend right=60]   &  & {[\f^0,\f^0]} \arrow[d, "\ab"]  \arrow[dd, "\rho^\sab_0", bend left=50]                                   \\
		{\f^{-1}_{\ab}=\frac{\f^{-1}}{[ [\f^0, \f^0],\f^{-1}]}} \arrow[d, "\rho^{\ab}_{-1}"] \arrow[rr, "d^{\ab}"] &  & {\frac{[\f^0,\f^0]}{[[\f^0,\f^0],[\f^0,\f^0]]}} \arrow[d, "\rho^{\ab}_{0}"]    \\
		\Gamma_2 \arrow[rr, "\partial"]                                                    &  & \Gamma^{\cl}_1                                                  
	\end{tikzcd}
	\vspace{-1em}
	\caption{
			 $\rho_0$ is defined in \Cref{thm:reutenauers_map}.
			 $\rho_{-1}$ is defined in \Cref{lem:rho_minus_1}
			 $\rho^\ab_{-1}, \rho^\ab_0$ are Lie algebra isomorphisms (since domain and codomain are abelian, these are just linear isomorphisms).
			 $\rho^\sab_{-1}, \rho^\sab_0$ are Lie algebra isomorphisms (with abelian codomain).
	}
	\label{fig:rho}
\end{figure}
\begin{remark}
	We have
	\begin{align*}
		\nu_I( \omega ) 
		=
		\partial_{i_1} \partial_{i_2} \dots \partial_{i_{p-2}} \omega_{e_{i_{p-1}}, e_{i_{p}} }(0),
	\end{align*}
	which is the definition used in \cite{kapranov2015membranes}.
	Indeed,
	\begin{align*}
		\nu_I( r^\alpha \diff{r}^j \wedge \diff{r}^k )
		&=
		\nu_I( \frac1{\alpha!} x^\alpha \diff{x}^j \wedge \diff{x}^k ) \\
		&=
		\begin{cases}
			\hphantom{-}1 & \text{ if } \alpha_\ell = \#\{ k \le p-2 \mid i_k = \ell \}, j=i_{p-1}, k=i_p \\
			-1 & \text{ if } \alpha_\ell = \#\{ k \le p-2 \mid i_k = \ell \}, j=i_{p}, k=i_{p-1} \\
			\hphantom{-}0 & \text{ else}.
		\end{cases}
	\end{align*}
\end{remark}
\begin{proof}[{Proof of \Cref{thm:reutenauers_map}}]~
	\begin{enumerate}[label=\roman*.]

		\item 
		By \Cref{thm:free_generators},
		the unique existence of $\rho_0$ follows from the universal property of free Lie algebras.

		\item
		Let $W \coloneqq [\freecm{n}^0,\freecm{n}^0]$.
		First, $\rho_{0}$ factors through the abelianization $W/[W,W]$, since,
		by construction, $\rho_0$ is a Lie morphism with abelian codomain.

		Further, for all $\sigma \in S_{p-2}$,
		\begin{align}
			\label{eq:first}
			&[\Z{i_1}, [\Z{i_2}, \dots [\Z{i_{p-2}}, [\Z{i_{p-1}}, \Z{i_p}]\dots]]] \notag \\
			&\qquad\quad 
			= [\Z{i_{\sigma(1)}}, [\Z{i_{\sigma(2)}}, \dots [\Z{i_{\sigma(p-2)}}, [\Z{i_{p-1}}, \Z{i_p}]\dots]]], \qquad \text{ mod } [W,W].
		\end{align}
		Indeed, for any $k, \ell \in [n]$ and any $I = (i_1, \dots, i_p) \in [n]^p$, $p\ge 2$, we have, mod $[W,W]$,
		\begin{align*}
			0
			&= \left[ [\Z{i_1}, [\Z{i_2}, \dots,[\Z{i_{p-1}}, \Z{i_p}]]], [\Z{k}, \Z{\ell}] \right] \\
			&= \left[ [[\Z{i_1}, [\Z{i_2}, \dots,[\Z{i_{p-1}}, \Z{i_p}]]], \Z{k}], \Z{\ell} \right]
			+
			\left[ \Z{k}, [[\Z{i_1}, [\Z{i_2}, \dots,[\Z{i_{p-1}}, \Z{i_p}]]], \Z{\ell}] \right] \\
			&=
			\left[ \Z{\ell}, [\Z{k}, [\Z{i_1}, [\Z{i_2}, \dots,[\Z{i_{p-1}}, \Z{i_p}]]]] \right]
			-
			\left[ \Z{k}, [\Z{\ell}, [\Z{i_1}, [\Z{i_2}, \dots,[\Z{i_{p-1}}, \Z{i_p}]]]] \right].
		\end{align*}
		Since $[W,W]$ is a Lie ideal, we have, mod $[W,W]$,
		\begin{align*}
			0 &= 
			[\Z{i_1}, [\Z{i_2}, \dots, [\Z{i_r}, [\Z{i_{r+1}}, \dots, [\Z{i_{p-1}}, \Z{i_p}]\dots]]]] \\
			&\qquad\quad 
			-
			[\Z{i_1}, [\Z{i_2}, \dots, [\Z{i_{r+1}}, [\Z{i_{r}}, \dots, [\Z{i_{p-1}}, \Z{i_p}]\dots]]]],
		\end{align*}
		for any $r=1,\dots,p-3$.
		Since $S_{p-2}$ is generated by transpositions, this shows \eqref{eq:first}.

		Now let
		\begin{align*}
			[\Z{i_1}, [\Z{i_2}, \dots [\Z{i_{p-2}}, [\Z{i_{p-1}}, \Z{i_p}]\dots]]] \in W. 
		\end{align*}
		By what we just showed, we can assume $i_1 \ge i_2 \ge \dots \ge i_{p-2}$,
		since re-ordering these indices does not change the image under $\rho_0 = \rho^\ab_0 \circ \pi^\ab$
		(where $\pi^\ab$ is the projection to the abelianization).
		By changing the sign (which is compatible with $\rho_0$), we can assume $i_{p-1}<i_p$.

		If $i_{p-2} \ge i_{p-1}$, we are in the setting of the previous point i.
		Otherwise, we apply the Jacobi identity to obtain
		\begin{align*}
			&[\Z{i_1}, [\Z{i_2}, \dots, [\Z{i_{p-2}}, [\Z{i_{p-1}}, \Z{i_p}]\dots]]] \\
			&\quad=
			[\Z{i_1}, [\Z{i_2}, \dots, [[\Z{i_{p-2}}, \Z{i_{p-1}}], \Z{i_p}]\dots]]
			+
			[\Z{i_1}, [\Z{i_2}, \dots, [\Z{i_{p-1}}, [\Z{i_{p-2}}, \Z{i_p}]]\dots]] \\
			&\quad=
			-
			[\Z{i_1}, [\Z{i_2}, \dots, [\Z{i_{p}}, [\Z{i_{p-2}}, \Z{i_{p-1}}]\dots]]]
			+
			[\Z{i_1}, [\Z{i_2}, \dots, [\Z{i_{p-1}}, [\Z{i_{p-2}}, \Z{i_p}]\dots]]],
		\end{align*}
		which under $\rho_0$ is mapped to
		\begin{align*}
		 	&\partial_2( - s_{I'} s_{i_p} \diff{s}^{i_{p-2}} \wedge \diff{s}^{i_{p-1}} + s_{I'} s_{i_{p-1}} \diff{s}^{i_{p-2}} \wedge \diff{s}^{i_p} ) \\
			&\quad=
			- s_{I'} s_{i_p} s_{i_{p-2}} \diff{s}^{i_{p-1}} + s_{I'} s_{i_{p-1}} s_{i_{p-2}} \diff{s}^{i_p} \\
			&\quad=
			\partial_2 \nu_I,
		\end{align*}
		where $I = (i_1, \dots, i_{p-2}, i_{p-1}, i_p)$ and $I' = (i_1, \dots, i_{p-2})$.
		This proves the claim.

		\item
		First, $\rho^\ab_0$ is a morphism of $\GL(\R^n)$-representations.
		Indeed, let $A \in \GL(\R^n)$.
		Then, for all $J \in [n]^p$, (using Einstein summation)
		\begin{align*}
			&\rho^\ab_0( A \cdot [\Z{j_1}, [\Z{j_2}, \dots, [\Z{j_{p-1}}, \Z{j_p}]\dots]] ) \\
			&\quad=
			\rho^\ab_0( [A_{r_1 j_1} \Z{r_1}, [A_{r_2 j_2} \Z{r_2}, \dots, [A_{r_{p-1} j_{p-1}}\Z{r_{p-1}}, A_{r_p j_p} \Z{r_p}]\dots]] ) \\
			&\quad=
			A_{r_1 j_1} A_{r_2 j_2} \dots A_{r_{p-1}j_{p-1}} A_{r_p j_p}  \rho^\ab_0( [\Z{r_1}, [\Z{r_2}, \dots, [\Z{r_{p-1}}, \Z{r_p}]\dots]] ) \\
			&\quad=
			A_{r_1 j_1} A_{r_2 j_2} \dots A_{r_{p-1}j_{p-1}} A_{r_p j_p} 
			\partial_2( x_{r_1} \dots x_{r_{p-2}} \diff{x}^{r_{p-1}} \wedge \diff{x}^{r_p} ) \\
			&\quad=
			A \cdot \partial_2( x_{j_1} \dots x_{j_{p-2}} \diff{x}^{j_{p-1}} \wedge \diff{x}^{j_p} ).
		\end{align*}


		Let $W \coloneqq [\freecm{n}^0, \freecm{n}^0]$.
		We consider the graded piece of degree $\ell$
		of both $W/[W,W]$ and $\Gamma^\cl_1$.

		By \Cref{lem:cochain} \Cref{it:gamma_is_irrep},
		$(\Gamma^\cl_1)_\ell$ is the irrep of $\GL(\R^n)$
		corresponding to the shape $(\ell-1,1)$.

		By \Cref{thm:free_generators},
		$(W/[W,W])_\ell$ is the irrep of $\GL(\R^n)$
		corresponding to the shape $(\ell-1,1)$.

		By Schur's lemma (\cite[Lemma 1.7]{fulton2013representation})
		it remains to show that $\rho^\ab_0$
		is not the zero map.
		And indeed, for $I = (i_1, \dots, i_\ell) \in \I_{n,\ell}$,
		\begin{align*}
			\rho^\ab_0( [ \Z{i_1}, [\Z{i_2}, \dots, [\Z{i_{\ell-1}}, \Z{i_\ell}]\dots]] )
			=
			\partial_2 \nu_I
			\neq 0.
		\end{align*}
		%









	\end{enumerate}
\end{proof}

\medskip

Next we consider $\Gamma^{\cl}_2(\R^n)$, which is indexed by $\J_{n,p}$, $p\ge 2$.

\begin{theorem}
	\label{lem:rho_minus_1}
	Recall $\nu_I \in (\Omega^2(V))^0 \cong \Gamma^2(V)$ defined in \Cref{thm:reutenauers_map}.
	There is a unique linear map
	\begin{align*}
		\rho_{-1}\colon \f^{-1} \to \Gamma_2(\R^n),
	\end{align*}
	satisfying, for $p \in \N_{\ge2},\ i_j \in [n],\ j=1,...,p$ 
	\begin{align*}
		\rho_{-1}\left( [\Z{i_1}, [\Z{i_2}, \dots, [\Z{i_{p-2}}, \Z{i_{p-1}i_p}]\dots]] \right)
		=
		\nu_I.
	\end{align*}
	It is a morphism of $\GL(\R^n)$-representations.
	Moreover,
	\begin{enumerate}[label=\roman*.]
		\item
		$\rho_{-1}$ factors through a map $\rho^\sab_{-1}\colon \f^{-1}_\sab \to \Gamma_2(\R^n)$ of $\GL(\R^n)$-representations.
		Moreover, for any $x \in \f^{-1}, y \in [\f^0,\f^0]$,
		\begin{align*}
			\rho^\sab_{-1}([x,y]) = 0.
		\end{align*}

		\item
		$\rho^\sab_{-1}$ induces a linear isomorphism of $\GL(\R^n)$-representations
		\begin{align*}
			\rho^\ab_{-1}: \frac{ \f^{-1} }{ \left[[\f^0,\f^0], \f^{-1}\right] } \xrightarrow{\sim} \Gamma_2(\R^n).
		\end{align*}

		\item 
		\begin{align*}
			\rho^\sab_{-1}\colon \ker d^{\sab} \to \Gamma^\cl_2 = \ker \{\partial_2\colon \Gamma_2 \to \Gamma_1\},
		\end{align*}
		is a linear isomorphism of $\GL(\R^n)$-representations.
	\end{enumerate}

\end{theorem}
\begin{proof}
	Define $\rho_{-1}$ on $\mathfrak f^{-1}(\R^n)$ as
	\begin{align*}
		\rho_{-1}\left( [\Z{i_1}, [\Z{i_2}, \dots, [\Z{i_{p-2}}, \Z{i_{p-1}i_p}]\dots]] \right)
		=
		\nu_I, 
	\end{align*}
	for $p \in \N_{\ge2},\ i_j \in [n],\ j=1,...,p,\ i_{p-1} < i_p$.
	This is well-defined by \Cref{lem:basis_f_minus_1}.
	Changing the order of $i_{p-1}$ and $i_p$ changes the sign on both sides,
	so the identity is true for \emph{all} $I \in [n]^p$.
	As in the proof of \Cref{thm:reutenauers_map},
	this implies that $\rho_{-1}$ is a morphism of $\GL(\R^n)$-representations.

	To verify that $\rho_{-1}$ factors through the quotient of the semi-abelianization,
	we need to verify $d [\f^{-1}(\R^n),\f^{-1}(\R^n)] \subset \ker \rho_{-1}$.
	We show the stronger statement: 
	\begin{align}
		\label{eq:stronger}
		\left[ \f^{-1}(\R^n),[\f^0(\R^n),\f^0(\R^n)] \right] \subset \ker \rho_{-1}.
	\end{align}	

	To show this, it is enough that, for any $i_1, ..., i_p \in [n], i_{p-1}<i_p$, $p \ge 2$,
	and any $j_1, ..., j_q \in [n]$, $q \ge 2$,
	\begin{align*}
		\rho_{-1}\left(
			\Big[ [\Z{i_1}, [\Z{i_2}, \dots, [\Z{i_{p-2}}, \Z{i_{p-1}i_p}]\dots]],
			  [\Z{j_1}, [\Z{j_2}, \dots, [\Z{j_{q-2}}, [\Z{j_{q-1}}, \Z{j_q}]]\dots]] \Big]
		 \right) = 0.
	\end{align*}
	First, it is true for $q=2$:
	\begin{align*}
		&\rho_{-1}
		\left(
			\Big[ [\Z{i_1}, [\Z{i_2}, \dots, [\Z{i_{p-2}}, \Z{i_{p-1}i_p}]\dots]],
			  [\Z{j_1}, \Z{j_2}] \Big]
		 \right) \\
		 &\quad=
		 \rho_{-1}
		\Big(
			[ [ [\Z{i_1}, [\Z{i_2}, \dots, [\Z{i_{p-2}}, \Z{i_{p-1}i_p}]\dots]], \Z{j_1} ],
			    \Z{j_2} ] \\
				&\qquad\qquad
			+
			[ \Z{j_1}, [ [\Z{i_1}, [\Z{i_2}, \dots, [\Z{i_{p-2}}, \Z{i_{p-1}i_p}]\dots]], \Z{j_2} ] ]
		 \Big) \\
		 &\quad=
		 \rho_{-1}
		\Big(
			[ \Z{j_2}, [ \Z{j_1}, [\Z{i_1}, [\Z{i_2}, \dots, [\Z{i_{p-2}}, \Z{i_{p-1}i_p}]\dots]] ] ] \\
			&\qquad\qquad
			-
			[ \Z{j_1}, [ \Z{j_2}, [\Z{i_1}, [\Z{i_2}, \dots, [\Z{i_{p-2}}, \Z{i_{p-1}i_p}]\dots]] ] ]
		 \Big) \\
		 &\quad= 0.
	\end{align*}
	Now let it be true for arbitrary $q-1 \ge 2$.
	Let $I=(i_1, \dots, i_p)$, $J=(J',j_q)=(j_1, \dots, j_q)$, be given.
	Define
	\begin{align*}
		Z^{-1}_I   &\coloneqq [\Z{i_1}, [\Z{i_2}, \dots, [\Z{i_{p-2}}, \Z{i_{p-1}i_p}]\dots]], \quad
		Z^{0}_{J'} &\coloneqq [\Z{j_1}, [\Z{j_2}, \dots, [\Z{j_{q-1}}, \Z{j_q}]\dots]].
	\end{align*}
	Then
	\begin{align*}
		\rho_{-1}( [Z^{-1}_I, [\Z{j_1},Z^0_{J'}]] )
		&=
		\rho_{-1}(  [ [Z^{-1}_I,\Z{j_1}], Z^0_{J'} ] )
		+
		\rho_{-1}(  [ \Z{j_1}, [Z^{-1}_I,Z^0_{J'}] ] )
	\end{align*}
	Now
	\begin{align*}
		\rho_{-1}(  [ [Z^{-1}_I,\Z{j_1}], Z^0_{J'} ] )
		=
		\rho_{-1}(  [ Z^0_{J'}, [\Z{j_1}, Z^{-1}_I] ] )
		= 0,
	\end{align*}
	by induction hypothesis, since $|J'|=q-1$.
	Lastly, by \Cref{lem:basis_f_minus_1}, we have
	\begin{align*}
		[Z^{-1}_I, Z^0_{J'}] = \sum_K \alpha_K Z^{-1}_K,
	\end{align*}
	for some $\alpha_K \in \R$ and the sum is over all $K=(k_1, \dots, k_{p+q-1})$ with 
	$k_{p+q-2} < k_{p+q-1}$.
	Then
	\begin{align*}
		\rho_{-1}(  [ \Z{j_1}, [Z^{-1}_I,Z^0_{J'}] ] )
		=
		\sum_K \alpha_K \rho_{-1}( [\Z{j_1}, Z_K] )
		=
		s_{j_1} \rho_{-1}(\sum_K \alpha_K Z_K ) = 0.
	\end{align*}

	\begin{enumerate}[label=\roman*.]

		\item
			The map $\rho_{-1}$ factors through the semi-abelianization, by \eqref{eq:stronger}.
			Since $\ker \rho_{-1}$ is a $\GL(\R^n)$-submodule of $\f^{-1}$,
			we get a map of $\GL(\R^n)$-representations.

		\item
		 	By \Cref{lem:basis_f_minus_1},
			the elements
			\begin{align*}
				z_I \coloneqq [\Z{i_1}, [\Z{i_2}, \dots, [\Z{i_{p-2}}, \Z{i_{p-1}i_p}]\dots]], \quad p \in \N_{\ge2}, i_j \in [n], j=1,...,p,\ i_{p-1} < i_p,
			\end{align*}
			form a basis of $\f^{-1}$.
			Its image under abelianization $[z_I]_{\ab}$ hence span $\f^{-1}_\ab$.
			Restricting to $I \in \I_{n,p}$, 
			they are mapped under $\rho^\ab_{-1}$ to the 
			basis $\nu_I$ of $\Gamma_2(\R^n)$.
			Hence $[z_I]_{\ab}$, $I \in \I_{n,p}$, form a basis of $\f^{-1}_\ab$
			and $\rho^\ab_{-1}$ is an isomorphism.
			
		\item
			Consider the diagram in \Cref{fig:rho}.
			The outermost square commutes i.e.~$\rho_{0} \circ d = \partial \circ \rho_{-1}$.
			Indeed, this follows from
			\Cref{thm:reutenauers_map},\Cref{it:any}:
			for an $J \in [n]^p$
			\begin{align*}
				\rho_0 d [\Z{j_1}, [\Z{j_2}, \dots, [\Z{j_{p-2}}, \Z{j_{p-1}j_p}]\dots]]
				&=
				\rho_0 [\Z{j_1}, [\Z{j_2}, \dots, [\Z{j_{p-2}}, [\Z{j_{p-1}},\Z{j_p}]]\dots]] \\
				&=
				\partial_2 \nu_J \\
				&=
				\partial_2 \rho_{-1} [\Z{j_1}, [\Z{j_2}, \dots, [\Z{j_{p-2}}, \Z{j_{p-1}j_p}]\dots]].
			\end{align*}

			Now $\rho_{-1}$ factors through the abelianization, \Cref{eq:stronger}.
			Again, since $\ker \rho^\sab_{-1}$ is a $\GL(\R^n)$-submodule of $\f^{-1}$, we get a map of $\GL(\R^n)$-representations.
			Further $\rho_0$ factors through the abelianization, by construction.
			Further, $d[\f^{-1},\f^{-1}] \subset \ker d$, so that $d^\sab$ is well-defined.
			Lastly, $d[[\f^0,\f^0],\f^{-1}] \subset [[\f^0,\f^0],[\f^0,\f^0]]$, so that $d^\ab$ is well-defined.
			Since everything is well-defined and the outer square commutes,
			all squares commute.

			$\rho_0$ is surjective and then so are $\rho^\sab_0, \rho^\ab_0$.
			By \Cref{thm:free_differential_crossed_module}, as $\GL(\R^n)$-representations,
			\begin{align*}
				\f_\sab^{-1} \cong \ker d^\sab \oplus [\f^0,\f^0].
			\end{align*}

			Now the action of $[\f^0,\f^0]$ on $\f_\sab^{-1}$ acts on the direct summands as follows
			\begin{align*}
				[[\f^0,\f^0], \ker d^\sab] = 0 \subset \ker d^\sab, \qquad
				[[\f^0,\f^0], [\f^0,\f^0]] \subset [\f^0,\f^0].
			\end{align*}
			Indeed, only the first statement is non-trivial, so let $x \in \ker d^\sab$, $y \in [\f^0,\f^0]$.
			Then, by \Cref{prop:dgLie} there is $b \in \f^{-1}$ such that $d b = y$.
			Then 
			\begin{align*}
				0 = d [b,x] = [db,x] + [b,dx] = [y,x] \text{ in } \f^{-1}_\sab.
			\end{align*}

			Then, as $\GL(\R^n)$-representations (and hence, as abelian Lie algebras),
			\begin{align*}
				\f^{-1}_\ab
				= \f^{-1}_\sab / \left[ [\f^0,\f^0], \f^{-1}_\sab \right]
				\cong \ker d^\sab \oplus \frac{[\f^0,\f^0]}{\left[ [\f^0,\f^0], [\f^0,\f^0] \right]}
				\cong \ker d^\ab  \oplus \frac{[\f^0,\f^0]}{\left[ [\f^0,\f^0], [\f^0,\f^0] \right]}.
			\end{align*}
			Since the lower square in \Cref{fig:rho} commutes and $\rho^\ab_{-1}$ is an isomorphism
			we thus have
			\begin{align*}
				\rho^\ab_{-1}\colon \ker d^\ab \to \ker \partial,
			\end{align*}
			is an isomorphism of $\GL(\R^n)$-representations (and hence, as abelian Lie algebras).
	\end{enumerate}
\end{proof}


\section{Symbolic index}

We collect in this appendix commonly used symbols of the article and a page of reference where they occur.
Groups are denoted using capital letters, e.g., $G$, $H$, etc. Lie algebras are denoted using lowercase fraktur letters, e.g., $\mathfrak{g}$, \(\mathfrak{h}\), and so on.

 \begin{center}
\renewcommand{\arraystretch}{1.1}
\begin{longtable}{lll}
\toprule
Symbol & Meaning & Page\\
\midrule
\endfirsthead
\toprule
Symbol & Meaning & Page\\
\midrule
\endhead
\bottomrule
\endfoot
\bottomrule
\endlastfoot

		$\alpha \wedge \alpha'$ & Lie algebra valued $2$-form & \pageref{rmk:wedge}\\
	
		\(\overline{\mathfrak{g}}\) & Lie series over Lie algebra $\mathfrak{g}$ & \pageref{thm:universality}\\
	
		$\mathrm{G}((\mathbb{R}^n))$ & group-like elements in $\text{T}((\mathbb{R}^n))$ & \pageref{eq:Chen1}\\ 

		$J^{(ij)}$ & Jacobian minor & \pageref{eq:J_def}\\

		$\FL(\mathbb{R}^{n})$/$\FL(( \mathbb{R}^n))$ & free Lie algebra/series & \pageref{sec:prelim}\\

		$\actionGG$ 		& (Lie) group action& \pageref{def:xmod}\\

		$\actionAA$		& Lie algebra action& \pageref{def:diffxmod}\\

		$n \in \mathbb \N$ & ambient dimension, i.e.~we will consider functions $X\colon [0,1]^2 \to \R^n$ & \pageref{first-appearance-n}\\ 

		$\mag$ & 1D Magnus expansion (log-signature)& \pageref{thm:magnus_expansion_1cocycle}\\
	
		$\Omega$ & 2D Magnus expansion (2-log-signature)& \pageref{eq:magnus_general}\\

		$\Omega^{k}(\R^n,\mathfrak g)$, $k \ge 1$  & smooth\footnote{In fact, the only thing we use is that our 1-forms are Lipschitz and our 2-forms are bounded.} $k$-forms on $\R^n$ with values in Lie algebra $\mathfrak g$ & \\

		$\PD$ & path development, i.e., 1D signature (1-cocycle)& \pageref{thm:1cocycleIntegration}\\

		$\widehat\PD$ &1D germ (see \Cref{ss:1d_sewing}) & \pageref{ss:1d_sewing} \\
	
		$\SD$ & surface development, i.e., 2D signature (2-cocycle)& \pageref{thm:2cocycleIntegration}\\

		$\widehat\SD$ &2D germ (see \Cref{thm:2D_sewing}) & \pageref{eq:PD_SD}\\

		$\text{T}(\mathbb{R}^{n})$/$\text{T}(( \mathbb{R}^n))$ & tensor algebra/series over $\mathbb{R}^{n}$ & \pageref{sec:prelim}\\

		$\feedbackGG$		& feedback (or boundary) map in a crossed module of (Lie) groups& \pageref{def:xmod}\\
		
		$\mathfrak{t}$		& feedback (or boundary) map in a crossed module of Lie algebras & \pageref{def:diffxmod}\\

		$\actionFreeCmAA$ & Lie algebra action in  free differential crossed module& \pageref{thm:free_differential_crossed_module}
\end{longtable}
 \end{center}

\bibliographystyle{alpha}
\bibliography{twosig}


\end{document}